\documentclass[11pt]{article} 
\usepackage{amssymb} 
\usepackage{cite}
\usepackage{dsfont}
\usepackage{epsfig} 
\usepackage{rotate}
\usepackage{rotating}

\newtheorem{theorem}{Theorem}
\newtheorem{proposition}{Proposition}
\newtheorem{lemma}{Lemma}

\newcommand{\Q}{\mathbb{Q}}
\newcommand{\R}{\mathbb{R}}
\newcommand{\Z}{\mathbb{Z}}
\newcommand{\cA}{{\mathcal A}}
\newcommand{\cB}{{\mathcal B}}
\newcommand{\cD}{{\mathcal D}}

\newcommand{\cP}{{\mathcal P}}

\newcommand{\cT}{{\mathcal T}}
\newcommand{\beq}{\begin{equation}} 
\newcommand{\eeq}{\end{equation}}
\newcommand{\Leq}[1]{\label{#1}\end{equation}}
\newcommand{\bdm}{\begin{displaymath}} 
\newcommand{\edm}{\end{displaymath}}
\title{Renormalization of a one-parameter family of piecewise isometries}
 \author{J. H. Lowenstein and F. Vivaldi\dag}
\date{
{\small\it
Dept.~of Physics, New York University, 2 Washington Place, New York, NY 10003, USA
\\
\dag School of Mathematical Sciences, Queen Mary, University of London, London E1 4NS, UK
}
}
\begin{document} 
\maketitle
\begin{abstract}
We consider a one-parameter family of piecewise isometries of a rhombus.
The rotational component is fixed, and its coefficients belong to the 
quadratic number field $K=\mathbb{Q}(\sqrt{2})$. 
The translations depend on a parameter $s$ which is allowed to vary in an interval. 
We investigate renormalizability. We show that recursive constructions of 
first-return maps on a suitable sub-domain eventually produce a scaled-down 
replica of this domain, but with a renormalized parameter $r(s)$. 
The renormalization map $r$ is the second iterate of a map $f$ of the generalised 
L\"uroth type (a piecewise-affine version of Gauss' map). We show that
exact self-similarity corresponds to the eventually periodic points of $f$,
and that such parameter values are precisely the elements of the field $K$ 
that lie in the given interval. 

The renormalization process is organized by a graph analogous to those used
to construct renormalizable interval-exchange transformations. There are ten
distinct renormalization scenarios corresponding to as many closed circuits in
the graph. The process of induction along some of these circuits involves
intermediate maps undergoing, as the parameter varies, infinitely many bifurcations.

Our proofs rely on computer-assistance.
\end{abstract}
\def\JHLabstract{
We consider a one-parameter family of piecewise isometries (PWIs) of
a rhombus, with rotation angle $-\pi/4$ and a fixed point on the short diagonal, 
controlled by a parameter $s\in I=[0,\sqrt{2}]$.  
We investigate the renormalizability of the model, focusing on the scaling 
properties of a triangular sub-domain whose induced return map $\rho$ has five 
atoms which are deformed continuously for $s$ ranging over $I$.  
Recursive, multi-level construction of induced return maps produces, for all 
$s\in I$, a scaled-down replica of the triangle and its PWI, but with a renormalized 
parameter $r(s)$.  
We show that $r$ is the second iterate of a piecewise affine function $f$ of 
generalized L\"{u}roth type (a piecewise affine version of Gauss's map).  
Exact dynamical self-similarity corresponds to the eventually periodic points of $f$, 
which, according to one of our main theorems, are precisely those elements of $I$ 
in the algebraic number field $\Q(\sqrt{2})$.

The underlying return-map analysis involves a variety of intermediate PWIs evolving 
with $s$ with infinitely many bifurcations. The renormalization process over the 
entire interval $I$ is efficiently organized by a graph analogous to those used 
to construct renormalizable interval exchange transformations.  
The distinct multi-step renormalization scenarios correspond to the closed circuits 
on the graph, of which there are ten.  
The return-map analysis of each induction step is handled by computer-assisted 
iteration of a PWI along non-branching orbits, often supplemented by application 
of certain lemmas to treat efficiently the sequences of bifurcations.
}

\centerline{\small\it\today}

%-------------------------------------------------------------------------------
\section{Introduction}

In piecewise isometries (PWI), renormalizability is a key for a complete 
description of the dynamics.
The phase space of these systems is partitioned into domains, 
called \textit{atoms}, over each of which the dynamics is an isometry.
By choosing a sub-domain of the original space ---typically an atom or 
a union of atoms--- and considering the first-return map to it, one 
constructs a new system, the \textit{induced} PWI on the chosen domain.  
If this process is repeated, then it may happen that an induced system 
be conjugate to the original system. This circumstance usually leads 
to a detailed understanding of the dynamics.

In one-dimension (interval-exchange transformations ---IET)
there is a satisfactory theory of renormalization. The Rauzy induction 
gives a criterion for selecting an interval (not one of the atoms)
over which to induce, resulting in a new IET with the same number of 
atoms \cite{Rauzy,Veech}. 
This induction process is a dynamical system over a finite-dimensional
space of IETs, related to the continued fractions algorithm, which 
affords a good description of the parameter space of IETs \cite{Yoccoz}.

An important connection with Diophantine arithmetic is provided by 
the Boshernitzan and Carrol theorem \cite{BoshernitzanCarroll};
it states that in any IET defined over a quadratic number field, 
inducing on any of the atoms results in only finitely many distinct 
IETs, up to scaling.
For a two-interval exchange, this finiteness result reduces
to Lagrange's theorem on the eventual periodicity of the
continued fractions expansions of quadratic surds.
Furthermore, if a (uniquely ergodic) IET is renormalizable, then 
the scaling constant involved in renormalization is a unit in a 
distinguished ring of algebraic integers 
\cite{PoggiaspallaLowensteinVivaldi}.

%%%%%%%%% FIGURE
\begin{figure}[t]
\hfil\epsfig{file=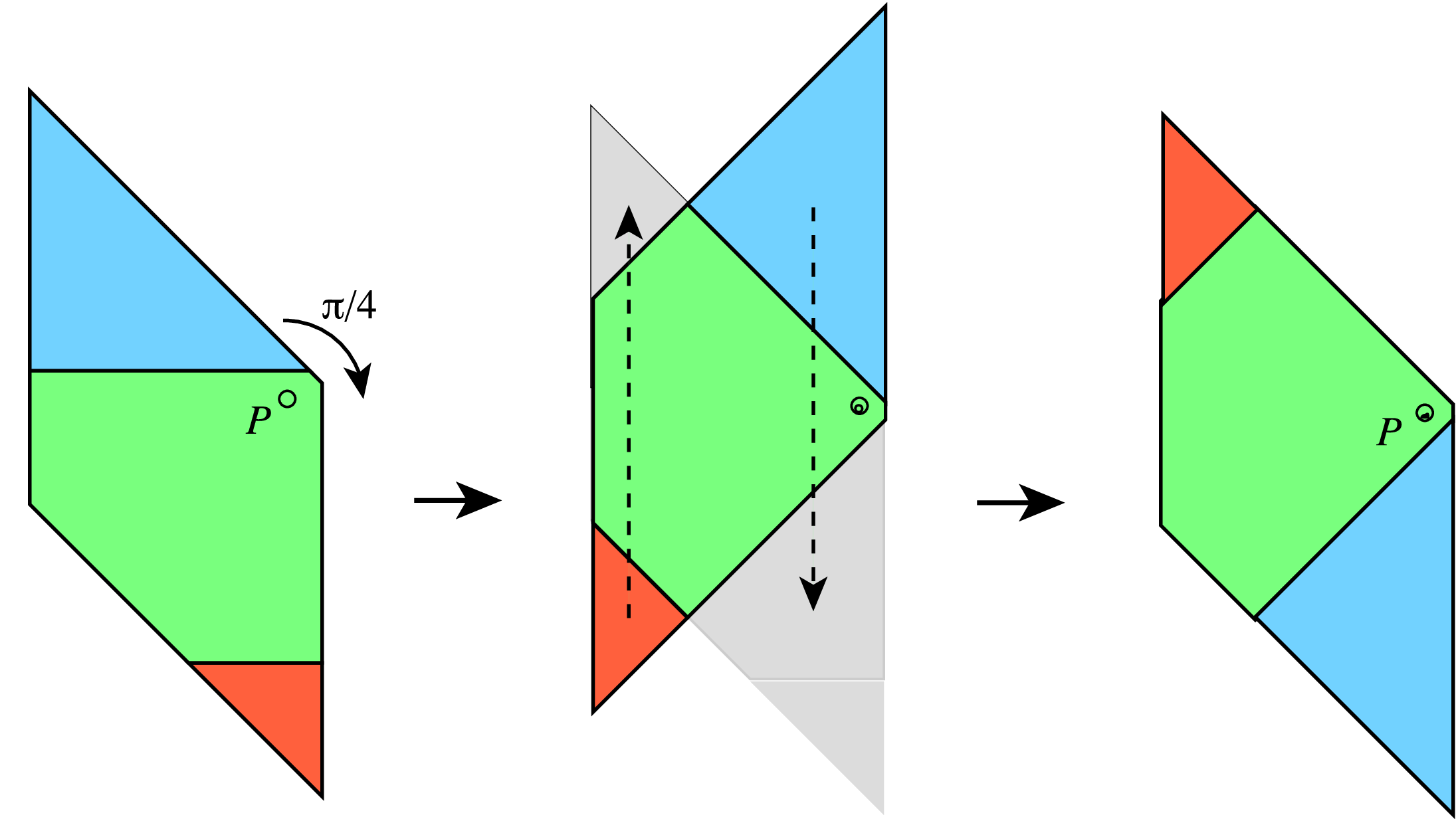,width=11 cm}\hfil
\caption{\label{fig:rhombusmap} \small Our model: a rhombus $\mathrm{R}$ with vertex 
angle $\pi/4$ is rotated clockwise by $\pi/4$ about a parameter-dependent fixed point 
$P$ located on the short diagonal. This leaves two triangular pieces outside the 
boundary of $\mathrm{R}$.  
These are translated vertically back into $\mathrm{R}$ to complete the piecewise isometry. } 
\end{figure}
%%%%%%%%%

In two dimensions general results are scarce 
\cite{Poggiaspalla:03,Poggiaspalla:06}. 
Until recently detailed results on renormalization were limited 
to special cases, defined over quadratic number fields 
\cite{LowensteinHatjispyrosVivaldi,AdlerKitchensTresser,KouptsovLowensteinVivaldi,%
AkiyamaBrunottePethoSteiner,Schwartz:09}. 
These results point consistently towards the existence of a 
two-dimensional analogue of the B\&C theorem. 
A more intricate form of renormalization has also been found in 
a handful of cubic cases \cite{GoetzPoggiaspalla,LowensteinKouptsovVivaldi}. 
In all cases, the renormalization constants are units in the ring 
of integers of the field of definition of the PWI.

More recently, renormalization has been studied in parametrised families. 
Hooper considered a two-parameter family of rectangle exchange 
transformations, and used techniques connected to the renormalization of Truchet 
tilings to establish results on the measure of the periodic and aperiodic sets 
of the map \cite{Hooper}.
In a substantial monograph \cite{Schwartz:14}, Schwartz determined
the renormalization group of a one-parameter family of polygon-exchange 
transformations, where the exchange is achieved by translations only. 
In this model, inducing on suitable domains leads to a conjugacy of the map 
to its inverse, accompanied by a change of parameter given by a 
piecewise-M\"obius map (a variant of Gauss' map).
Metric and topological properties of the limit set are also established. 

In the present work we consider a family of PWIs of a rhombus $\mathrm{R}$ 
with rotation angle $-\pi/4$ and a fixed point $P$ which is allowed to 
vary along the short diagonal, controlled by a real parameter $s$ (see figure \ref{fig:rhombusmap}).  
The choice of rotation determines the quadratic number field $\mathbb{Q}(\sqrt{2})$.
Some special parameter values in this field ($s=0,1+\sqrt{2},2(1+\sqrt{2}))$ have 
received much attention 
\cite{AdlerKitchensTresser,KouptsovLowensteinVivaldi,AkiyamaBrunottePethoSteiner};
they correspond to $P$ being the centre and a vertex of the rhombus.
These PWIs show exact self-similarity, and, as a result, their dynamics is well-understood.

For all maps in our family, the edges of the atoms have normal vectors in 
$\mathbb{Q}(\sqrt{2})^2$, while the Cartesian coordinates of the vertices 
belong to $\mathbb{Q}(\sqrt{2})+\mathbb{Q}(\sqrt{2})s$, which is a 
two-dimensional vector space over $\Q(\sqrt{2})$. The same is true of the 
first-return map induced on any of the atoms, and recursively, of any 
higher-level induced return maps. 
This arithmetical environment will have a profound effect on the renormalization. 

\begin{figure}[h]
\hfil\epsfig{file=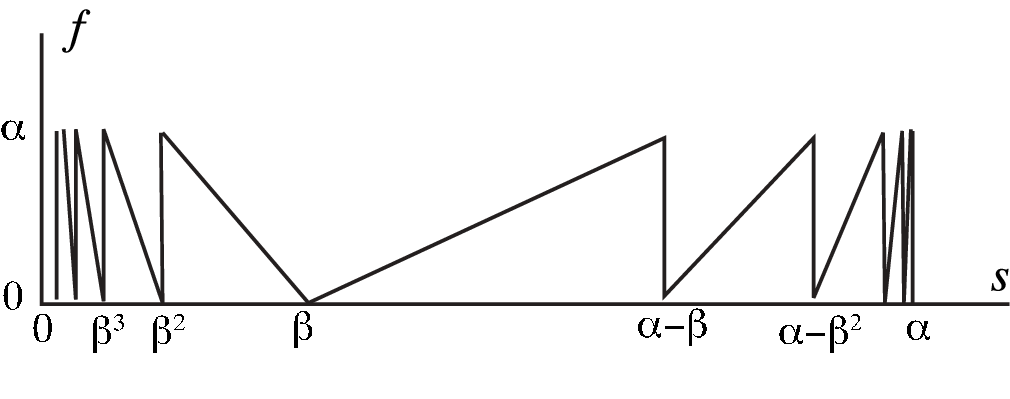,width=8cm}\hfil
\caption{\label{fig:f}\small
The piecewise-affine renormalization function $f$ which controls the 
renormalization of the piecewise isometry $\rho=\rho(s)$ defined
in section \ref{section:Model}.
Here $\alpha=\sqrt{2}$, $\beta=\alpha-1$, and the function $f$ 
maps the parametric interval $I=[0,\alpha)$ into itself. 
}
\end{figure}

We have two main results, theorems 1 and 2, stated in section \ref{section:MainResults}. 
In theorem 1 we prove that the parametrised rhombus map, restricted to the parameter 
interval $s\in I=[0,\sqrt{2}]$, induces in one of its triangular atoms a 
\textit{renormalizable} five-atom PWI, i.e., one which, after repeated inductions, 
recurs at smaller length scales with a piecewise-affine change of parameter $s\mapsto r(s)$. 
We show that $r$ is the second iterate of a modified L\"{u}roth map $f$ 
(a piecewise affine version of the Gauss map for continued fractions), 
which is shown in figure \ref{fig:f}. We also show that all scaling 
constants are units in the ring $\mathbb{Z}[\sqrt{2}]$.

Exact self-similarity is achieved if the induction process eventually 
reproduces a value of $s$ which has already been encountered, i.e., if 
$s$ is an eventually periodic point of $r$. In theorem 2 we prove that these
parameter values are precisely the elements of $\Q(\sqrt{2})\cap I$.
Note that, unlike the classic case of continued fractions, here 
eventual periodicity is associated with a single quadratic field.
This arithmetical characterisation of renormalizability provides 
additional evidence for the existence of an analogue of the B\&C 
theorem for polygon-exchange transformations.
%%%%%%%%%%%%%%%%%%%%%%%%%%%%%%%%%%%%%%%%%%%%%%%% FIGURE 
%%%%%%%%%%%%%%%%%%%%%%%%%%%%%%%%%%%%%%%%%%%% 

The discontinuities of the function $r$ accumulate at the infinitely many zeros 
of $f$ (see figure \ref{fig:r}).
As a result, the return-map dynamics is highly non-uniform, with ever-increasing 
return times as $s$ approaches such discontinuities.
The number of qualitatively distinct renormalization scenarios can be reduced to ten.  
The simplest of these involves a single induction, and it applies to the case in 
which both $s$ and $f(s)$ are in the middle of the interval $I$ .  

\begin{figure}[h]
\hfil\epsfig{file=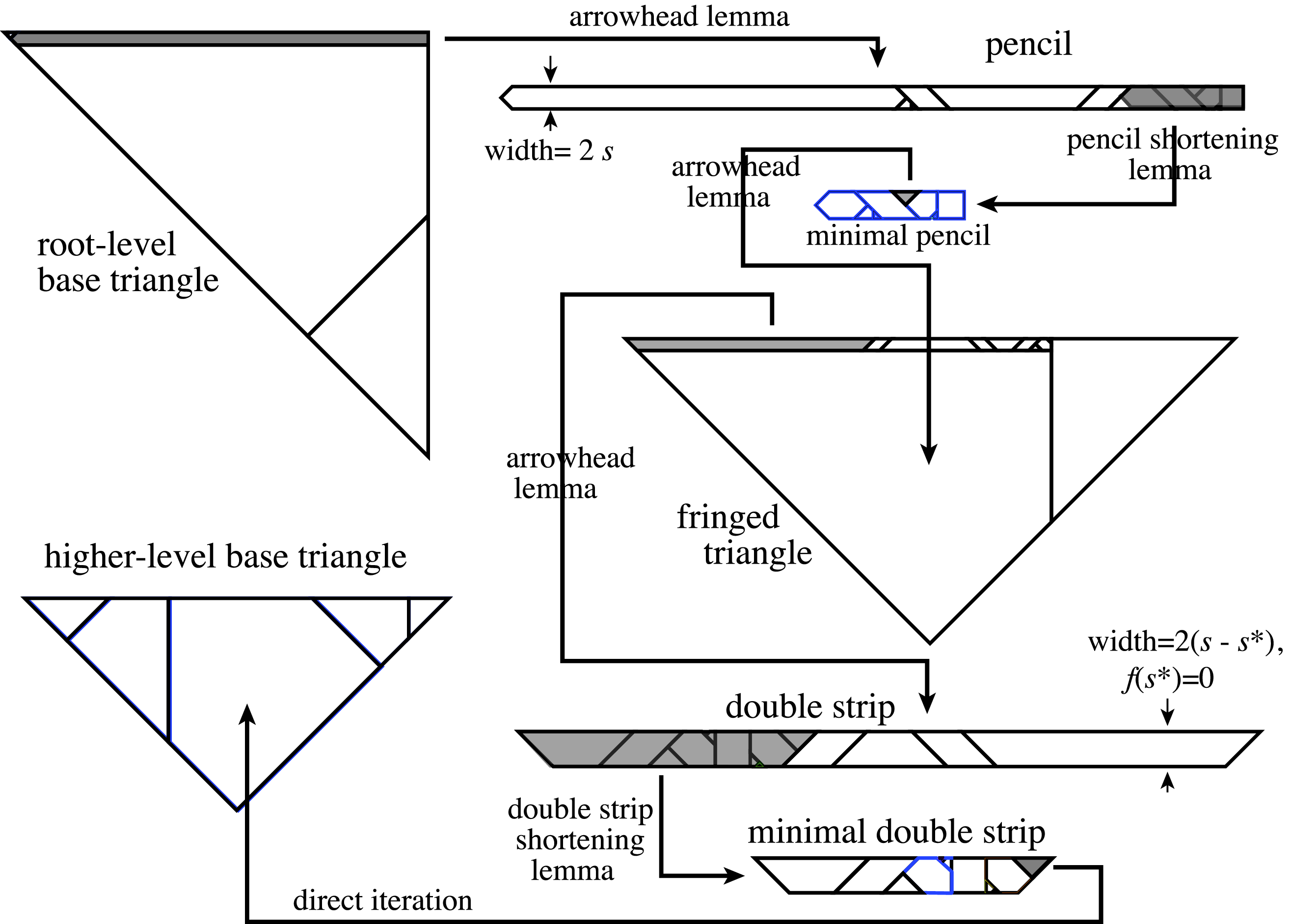,width=5 in}\hfil
\caption{\label{fig:RubeGoldberg} \small Sketch of the renormalization scenario for $s$ 
close to 0 and $f(s)$ close to $\sqrt{2}$.  
In the figure, the relative sizes of the various objects are greatly distorted in order 
to reveal their structure.}
\end{figure}

At the opposite extreme are, among others, those parts of $I$ for which both $s$ 
and $\sqrt{2}-f(s)$ are small.  
As a preview of renormalization dynamics, let us briefly sketch this scenario,
represented schematically in figure \ref{fig:RubeGoldberg}.
In this scheme, the number of inductions remains fixed at six, and it involves
the same types of induced PWIs (the \textit{pencil}, the \textit{fringed triangle}, 
and the \textit{double strip}). For either $s$ or $f(s)$ approaching zero or 
$\sqrt{2}$, the numbers of atoms increases without bounds, but in a tightly 
controlled manner (the complexity increases logarithmically).
The return times are are also unbounded.
%%%%%%%%% FIGURE
%%%%%%%%%

The map $f$ serves as a phase function for the evolution of the pencil with decreasing 
$s$, with two new atoms emerging at the boundary whenever $f(s)$ passes through unity.  
Analogously, bifurcations of the double strip occur at the zeros of $r(s)$.  
The spatial scaling accompanying the renormalization is governed by the 
successive narrowing of the widths of the pencil and the double strip. 

Of the six induction steps, two are `shortenings', exploiting the repetitive, 
quasi-one-dimensional structures of pencils and double strips.  
Steps 1, 3, and 4 appear to be more complicated, but in fact are all grounded 
in one simple dynamical sub-system, the \textit{arrowhead}.  
Once we have formulated that underlying dynamics in the Arrowhead Lemma of 
section \ref{section:Arrowheads}, all aforementioned induction steps can
be split into two manageable parts: the first, with short, fixed-length 
return orbits which can be constructed by explicit iteration, and the 
second, in which the Arrowhead Lemma accounts for all of the 
$s$-dependent bifurcations.  

\subsection{Structure of the paper}
%The rest of the paper is organized as follows: \medskip

\noindent \textit{Sections 2--4: Preliminaries, and formulation of the main results}\smallskip

In section \ref{section:Definitions} we provide definitions and notation 
which will be used throughout, followed by the specification of our model 
in section \ref{section:Model}. 
In section \ref{section:MainResults} we define the renormalization
function $f$ and $r$ and state our two main theorems.  
The first of these, theorem 1, describes the renormalization of our central 
geometric structure, the two-parameter base triangle with its five-atom piecewise isometry. 
The second, theorem 2, introduces the symbolic dynamics of the renormalization 
and establishes the connection between exact self-similarity and membership in 
the field $\mathbb{Q}(\sqrt{2})$. \smallskip

\noindent \textit{Sections 5,6: Symbolic dynamics, L\"{u}roth expansion, proof of theorem 2}\smallskip

After introducing a symbolic representation of parameter space in section  
\ref{section:SymbolicDynamics}, we derive the L\"{u}roth expansion 
and prove theorem 2 using a contraction argument on a lattice.
\medskip

\noindent \textit{Section 7:  Renormalization overview}\smallskip

Here we provide roadmaps for navigating the rest of the paper, organized  
around a renormalization graph, analogous to the Rauzy graphs of 
interval exchange theory \cite{Rauzy}. 
Oriented edges of the graph correspond to various return-map inductions, 
and closed circuits describe the ten distinct pathways by which one succeeds 
in renormalizing the initial base triangle on various subsets of the parameter 
interval $I$.\medskip

\noindent \textit{Sections 8,9:  Machinery for proving theorem 1}\smallskip

In section \ref{section:PrototypeDressedDomains} we define the 
\textit{parametric dressed domains} (tiled polygons with associated piecewise 
isometries) which correspond to the vertices of the renormalization graph: 
the pencil, the fringed triangle, and the double strip. In addition, we prove 
two Shortening Lemmas which allow one to reduce pencils and double strips to their minimum lengths.  
Section \ref{section:Arrowheads} is devoted to the arrowhead and its dynamics, 
culminating in the Arrowhead Lemma.\medskip

\noindent \textit{Section 10: Proof of theorem 1}\smallskip

Through a sequence of propositions, we treat all of the various induction 
steps needed to prove theorem 1 over the entire interval $I$.  
Some of the induction steps require only explicit iteration of PWIs 
over short orbits of fixed length, while others require, in 
addition, application of one of the Shortening Lemmas or the Arrowhead 
Lemma.\medskip  

\noindent \textit{Section 11: Temporal scaling and fractal structure}\smallskip

As a by-product of section 10, we obtain the incidence matrices which determine 
the renormalization's return times and asymptotic temporal scaling.  
We obtain explicit formulae for these matrices in section \ref{section:IncidenceMatrices}, 
where we also explore briefly the implications for the $s$-dependence of the 
Hausdorff dimension of the exceptional set. \medskip

Many steps in the proofs have required computer assistance, invariably to 
verify statements concerning finite orbits of PWIs.
This requires geometrical transformations of polygons (translations, rotations, 
reflections, etc.), polygon inclusion and disjointness tests, plus 
a fair amount of book-keeping. All arithmetic is performed exactly in 
$\mathbb{Q}(\sqrt{2})$ and  $\mathbb{Q}(\sqrt{2})+\mathbb{Q}(\sqrt{2})s$.
The complexity of the computation is manageable, because all of the heavy 
lifting is taken care by the Arrowhead and Shortening lemmas, which are proved 
analytically. A complete description and listing of our procedures, in the form 
of \textrm{Mathematica}\textsuperscript\textregistered  functions, together with all of the calculations 
participating in the proof of theorem 1, may be found in the Electronic 
Supplement \cite{ESupplement}. 

\bigskip\noindent
{\sc Acknowledgements:} \/
JHL and FV would like to thank, respectively, the School of
Mathematical Sciences at Queen Mary, University of London,
and the Department of Physics of New York University, for their hospitality.

%----------------------------------------------------------------%%%SECTION
\section{Definitions and notation}\label{section:Definitions}

Throughout this article we adopt the notation
$$
\alpha=\sqrt{2},\quad
\beta=\sqrt{2}-1,\quad
\omega=\sqrt{2}+1.
$$
The arithmetical environment is the algebraic number field 
$\mathbb{Q}(\alpha)$ with its ring of integers
$\mathbb{Z}[\alpha]$, given by
\begin{equation}\label{eq:FieldRing}
\mathbb{Q}(\alpha)=\{x+y\alpha\,:\,x,y\in\mathbb{Q}\},
\qquad
\mathbb{Z}(\alpha)=\{m+n\alpha\,:\,m,n\in\mathbb{Z}\}.
\end{equation}
The number $\omega$, which will be shown to determine the scaling, 
is the fundamental unit in $\mathbb{Z}[\alpha]$ (see \cite[chapter 6]{Cohn}).
Note that $\beta=\omega^{-1}$ is also a unit.

Our system depends on a parameter $s$, and to represent parameter 
dependence we consider the set
\begin{equation}\label{eq:ModuleS}
{\mathbb{S}}=\Q(\alpha) + \Q(\alpha) s
\end{equation}
Here, and below, $s$ is regarded as an indeterminate; hence the set 
${\mathbb{S}}$ is a two-dimensional vector space over $\Q(\alpha)$ 
(a $\Q(\alpha)$-module) whose elements are degree-one polynomials in $s$.

\subsection{Planar objects\label{section:PlanarObjects}}

A \textit{tile} is an open convex polygon whose edges have outward normal vectors 
taken from the set $\{\mathbf{u}_m=(\cos{m\pi/4},\sin{m\pi/4}): m=0,\ldots,7\}$.  
The equations of the edges of an $n$-sided tile are thus of the form
$$ 
\mathbf{u}_{m_i}\cdot (x,y) = b_i, \quad i=1,\ldots,n, 
$$
where the `octagonal coordinates' $b_i$ belong to $\mathbb{S}$.
The parameter $s$ allows for the continuous deformation of tiles. 

We represent an $n$-sided tile $X$ with edge orientations $m_1,\ldots,m_n$ and 
octagonal coordinates $b_1,\ldots,b_n$ with the bracket notation
\begin{equation}\label{eq:Tile}
X=[(m_1,\ldots,m_n), (b_1,\ldots,b_n)].
\end{equation}
Tile names will always be capital italic letters. 

A \textit{domain} is a union of open polygons whose edges are specified by octagonal 
orientations $\mathbf{u}_m$ and coordinates in the module $\mathbb{S}$. 
The polygons need not be convex. Domain names will always be capital Roman letters.

A \textit{tiling} ${\bf X}$ is a set of tiles,
$$
{\bf X}= \{X_1,\ldots,X_n\}.
$$
A tiling ${\bf X}$ is always associated with a domain ${\rm X}$, the {\it span} of ${\bf X}$ , 
$$
{\rm X}= {\rm span}({\bf X})\stackrel{\rm def}{=} {\rm Int}(\overline{\cup_{i=1}^n X_i}).
$$
The union of two tilings is a tiling.

\subsection{Isometry group}
We employ a group $\mathfrak{G}$ of transformations of planar objects, 
to specify their locations and orientations, and to describe their dynamical evolution. 
The group comprises the rotations and reflections of the symmetry group of the 
regular octagon (the dihedral group $D_8$) together with translations in $\mathbb{S}^2$. 
We define $\mathfrak{G}_+$ to be the subgroup of orientation-preserving transformations in $\mathfrak{G}$, 
i.e., those with Jacobian determinant equal to $+1$.

We adopt the following notation:
\begin{itemize}
\item [] $\mathtt{U}_n$: reflection about the lines generated by the
extended set of vectors
$$
\mathbf{u}_k=(\cos(k\pi/4),\sin(k\pi/4)),\quad k=0,\frac{1}{2},1,\frac{3}{2},\ldots,
$$
where we allow for half-integer indices, to be taken modulo 8;
\item [] $\mathtt{R}_n$: rotation by the angle $n\pi/4,\; n=0,\ldots,7$;
\item [] $\mathtt{T}_{\bf d}$: translation by ${\bf d}\in \mathbb{S}^2$.
\end{itemize}
Thanks to the product formulae
$$
\mathtt{T}_{\bf d} \mathtt{T}_{\bf e} = \mathtt{T}_{{\bf d}+{\bf e}},\qquad
\mathtt{R}_m \mathtt{R}_n = \mathtt{R}_{m+n},\qquad
\mathtt{U}_m \mathtt{U}_n = \mathtt{R}_{2(m-n)},
$$
and commutation relations
\begin{equation}\label{eq:CommutationRelations}
\mathtt{R}_n \mathtt{T}_{\bf d} = \mathtt{T}_{\mathtt{R}_n {\bf d}} \mathtt{R}_n,\qquad
\mathtt{R}_n \mathtt{U}_m =\mathtt{U}_{m+n} \mathtt{R}_n,\qquad
\mathtt{U}_n \mathtt{T}_{\bf d} = \mathtt{T}_{\mathtt{U}_n {\bf d}} \mathtt{U}_n,
\end{equation}
we can write an arbitrary element $\mathtt{G}$ of ${\mathfrak G}$ in the canonical form
$$
\mathtt{G}_{n,m,{\bf d}}= \mathtt{T}_{\bf d} \mathtt{R}_m \mathtt{U}^n_0,
$$
with
$$
n\in\{0,1\},\quad
m\in\{0,1,\ldots,7\},\quad
\mathbf{d}\in \mathbb{S}^2.
$$

In general, we will write $\mathcal{X}\sim\mathcal{Y}$ to indicate that 
$\mathcal{X}=\mathtt{G}(\mathcal{Y})$ for some $\mathtt{G}\in\mathfrak{G}$.  
As $\mathfrak{G}$ is a group, this is an equivalence relation.  
An equivalence class consists of planar objects which are congruent 
(same shape and size) up to a reflection.
Appendix A contains a catalogue of standard representatives of equivalence 
classes, used extensively in this work.

\subsection{Dressed domains and subdomains}
We define a \textit{dressed domain} to be a triple
\begin{equation}\label{eq:DressedDomain}
\mathcal{X}=({\mathrm{X}, \mathbf{X}, \rho}),
\end{equation}
where $\mathrm{X}$ is a domain, $\mathbf{X}$ is a tiling with span $\mathrm{X}$, 
and $\rho$ is a mapping which acts on each element of $\mathbf{X}$ as an 
isometry in $\mathfrak{G}_+$.  
We will describe $\rho$ as a \textit{piecewise isometry} or \textit {domain map} 
acting on $\mathrm{X}$, with $\mathbf{X}$ comprising the set of its \textit{atoms}.  
Dressed domains will always be denoted by capital script letters.  
Under the action of an isometry $\mathtt{G} \in \mathfrak{G}$, a dressed domain transforms as
$$
\mathtt{G}(\mathcal{X})= \mathtt{G}({\mathrm{X}, \mathbf{X}, \rho})= (\mathtt{G}(\mathrm{X}),\{\mathtt{G}(X_1),\mathtt{G}(X_2),\ldots\}, \mathtt{G}\circ\rho\circ \mathtt{G}^{-1}).
$$
To emphasize the association of a mapping $\rho$ with a particular dressed domain $\mathcal X$, 
we will use the notation $\rho_{\stackrel{}{\mathcal{X}}}$. 

Let $\mathcal{X}= (\mathrm{X}, \mathbf{X}, \rho_{\stackrel{}{\mathcal{X}}})$ be a dressed domain, 
and let $\mathrm{Y}$ be a sub-domain of $\mathrm{X}$.  
We denote by $\rho_{\stackrel{}{\mathcal{Y}}}$ the first-return map on $\mathrm{Y}$ induced by
$\rho_{\stackrel{}{\mathcal{X}}}$.  
We call the resulting dressed domain 
$\mathcal{Y}= (\mathrm{Y}, \mathbf{Y}, \rho_{\stackrel{}{\mathcal{Y}}})$ 
a \textit{dressed subdomain} of $\mathcal{X}$, writing
\begin{equation}\label{eq:XtoY}
\mathcal{X} \rightarrow \mathcal{Y}.
\end{equation}

A \textit{prototype} is a canonical representative of an equivalence class of dressed domains.
If $\widehat{\mathcal{X}}$ is a prototype and $\mathcal{Y}\sim\widehat{\mathcal{X}}$, then
the \textit{parity} $\pi(\mathcal{Y})$ of $\mathcal{Y}$ is the Jacobian determinant of the 
isometry in $\mathfrak{G}$ relating $\widehat{\mathcal{X}}$ to $\mathcal{Y}$.

The dressed subdomain relation (\ref{eq:XtoY}) enjoys the important property of 
\textit{scale invariance}, namely invariance under an homothety. Indeed if
$\mathtt{S}_\omega$ denotes scaling by a factor $\omega$, then in the data 
(\ref{eq:Tile}) specifying a tile, the orientations $m_k$ remain unchanged,
while the octagonal coordinates $b_k$ scale by $\omega$. 
Moreover, the identity
$$
\mathtt{S}_\omega \mathtt{T}_\mathbf{d}\mathtt{R}_n= \mathtt{T}_{\omega \mathbf{d}}\mathtt{R}_n\mathtt{S}_\omega.
$$
shows that the piecewise isometries $\rho$ scale in the same way.
We conclude that the relation (\ref{eq:XtoY}) is preserved if the dressed 
domain parameters are scaled by the same factor for both members.

We shall be dealing with renormalizability of dressed domains depending on 
a parameter $s$ ---the \textit{parametric dressed domains}.
The parameter $s$, ranging over an interval $I$, controls the `shape' of the domain. 
For reasons that will become clear below, it is useful to re-parametrise the system
with a pair $(l,h)$ where $l$ is a `size' parameter ranging over the positive real numbers,
and $h=sl$. So we shall write $\mathcal{X}=\mathcal{X}(l,h)$.
Note that a parametric dressed domain need not have a fixed number of atoms.
Indeed many of the parametric dressed domains introduced in section 
\ref{section:PrototypeDressedDomains} feature an infinite sequence of bifurcations,
each producing a change in the number and shapes of its atoms.

\subsection{Renormalizability of dressed domains}\label{section:Renormalizability}
A dressed domain $\mathcal {X}=({\mathrm{X}, \mathbf{X}, \rho_{\stackrel{}{\mathcal X}}})$ is 
\textit {strictly renormalizable} if $i)$ there exists a dressed subdomain $\mathcal{Y}$ of 
$\mathcal{X}$ and a dressed subdomain $\mathcal{Y}^*$ of $\mathcal{Y}$, which differs from 
$\mathcal{Y}$ by a  contracting scale transformation (homothety) composed with an isometry from $\mathfrak{G}$; \/  
$ii)$ the domain $\mathrm{X}$ has the \textit{recursive tiling property}, namely
it is completely tiled (ignoring sets of zero measure) by the return orbits of the atoms of 
$\mathcal{Y}^*$, together with the periodic orbits of a finite set of tiles. 

This is the simplest version of renormalizability.
Its implications for a planar piecewise isometry are well-known 
(see, for example, \cite{Poggiaspalla:06}, \cite[chapter 2]{Lowenstein}).  
Thus one can iterate the process at will, and with each iteration more and more 
periodic domains of finer and finer scales are revealed, leading to a full 
measure of periodic tiles in the limit.  
Simultaneously, the return orbits of the rescaled copies of $\mathcal{Y}$, provide 
finer and finer coverings of the \textit{exceptional set} complementary to all periodic tiles.  
While the latter has vanishing measure, its dimension is not trivial.  
Standard arguments \cite{Falconer,Lowenstein} show that the Hausdorff dimension of the exceptional set is given 
by $d_H= -\log \tau / \log \omega$, where $\omega$ and $\tau$ are, respectively the 
asymptotic spatial and temporal scale factors associated with the renormalization. 
The asymptotic spatial scaling is known, since each renormalization step is 
accompanied by multiplication by the same $\omega$. 
The temporal scaling is more subtle, requiring construction and diagonalization of the 
stepwise \textit{incidence matrix} $\mathrm{M}$, whose $i,j$th component gives ---in the 
above notation--- the number of times that the return orbit of atom $Y^*_j$  visits atom $Y_i$.   
The scale factor $\tau$ governing the asymptotic increase in length of the return orbits 
is given by the largest eigenvalue of $\mathrm{M}$.

A parametric dressed domain $\mathcal{Y}(l,h),\; l\in\R_+,\; s=h/l \in I$, is said to
be \textit{renormalizable} if there exist a piecewise smooth renormalization map 
$r:I\rightarrow I$, and an auxiliary scaling function $\kappa:\R_+\rightarrow (0,1)$ such 
that for every choice of $l$ and $h$, the dressed domain $\mathcal{Y}(l,h)$ has a dressed
subdomain congruent to $\mathcal{Y}(l',h')$ with $(l',h')=(\kappa(s)l,r(s)\kappa(s)l)$, 
and moreover the recursive tiling property is satisfied.  
In the present work, the renormalization map is piecewise-affine (as opposed to the 
piecewise-M\"obius map of \cite{Schwartz:14} and Gauss' map) with derivative equal 
to $1/\kappa$. Furthermore, all values assumed by $\kappa$ are units in the ring
$\mathbb{Z}[\sqrt{2}]$.
Note that $r$ and $\kappa$ depend only on $s$, a requirement of scale invariance. 
A parametric dressed domain which, for all valid parameter values, has a renormalizable 
parametric dressed sub-domain, with recursive tiling, will also be regarded as renormalizable.

If a parametric dressed domain $\mathcal{Y}(l,h)$ is renormalizable, we can consider 
those parameter values for which $\mathcal{Y}$ is strictly renormalizable. Because of 
scale invariance, if $\mathcal{Y}(l,h(l,s))$ is strictly renormalizable for $s=s_0$ 
and some $l$, then it is so for any $l$. 
It then follows that the $s$-values of strict renormalizability
are precisely the eventually periodic points of the function $r$.
A virtue of our model is an arithmetical characterization of these parameter
values: they are precisely the elements of the quadratic number field $\mathbb{Q}(\alpha)$.

The above definition of renormalizability is tailored to our model and it is
conceivable that in more general situations the recursive tiling property 
may require participation of more than one renormalizable parametric dressed 
sub-domain.  

%%%%%SUBSECTION
\subsection{Computations}\label{section:DirectIteration}
All computations reported in this work are exact, employing integer and polynomial 
arithmetic with Mathematica\textsuperscript\textregistered. For fixed parameter value, the computations take place in the
algebraic number field $\Q(\alpha)$ ---see (\ref{eq:FieldRing})--- whereas the parametric 
dependence requires computations in the module $\mathbb{S}$ defined in (\ref{eq:ModuleS}). 

All relevant objects are represented by data structures of elements of these two arithmetic sets.
In particular, we shall be concerned with finite orbits of polygonal domains under the domain
map of a dressed domain, which is an isometry in $\mathfrak{G}$.
To perform these computations we employ the procedures of our CAP Toolbox,
available in the Electronic Supplement \cite{ESupplement}.

In such processes, one must determine membership of points to polygons and intersections 
of polygons, which requires the evaluation of inequalities. Since the latter are expressed 
by affine functions of $s$, it suffices to check the inequalities at the endpoints of the 
assumed $s$-interval. All these boundary values belong in the field $\Q(\alpha)$, and the 
inequalities are evaluated by estimating $\alpha=\sqrt{2}$ via a pair of sufficiently close 
convergents in its continued fraction expansion. In this way we are able to establish 
statements valid over intervals of parameters.

Typically we are given a one-parameter family of piecewise isometries $\rho(s)$ of a dressed 
domain $\mathcal{R}(s)$, which we use to move each tile in the domain from an initial
position to a pre-assigned destination, checking at each step that no tile arrives at the 
wrong destination. Each iteration involves testing a number of half-plane inequalities 
to determine which atom $R_i(s)$ of $\rho(s)$ contains a particular tile, followed by 
application of the relevant isometry $\rho_i(s)$ to map the tile forward. 
When constructing a finite orbit (typically a return orbit), we keep track of the 
atoms visited, obtaining at the end the symbolic paths and incidence matrices of the orbits.  
The recursive tiling property defined in section \ref{section:Renormalizability} is established
by adding up the areas of the tiles of all the orbits, and comparing it with the total area 
of the parent domain.

Henceforth we will refer to this technique as \textit{direct iteration}.

%
%
%-------------------------------------------------------------------------------
\section{The model}\label{section:Model}

\begin{figure}[t]
\hfil\epsfig{file=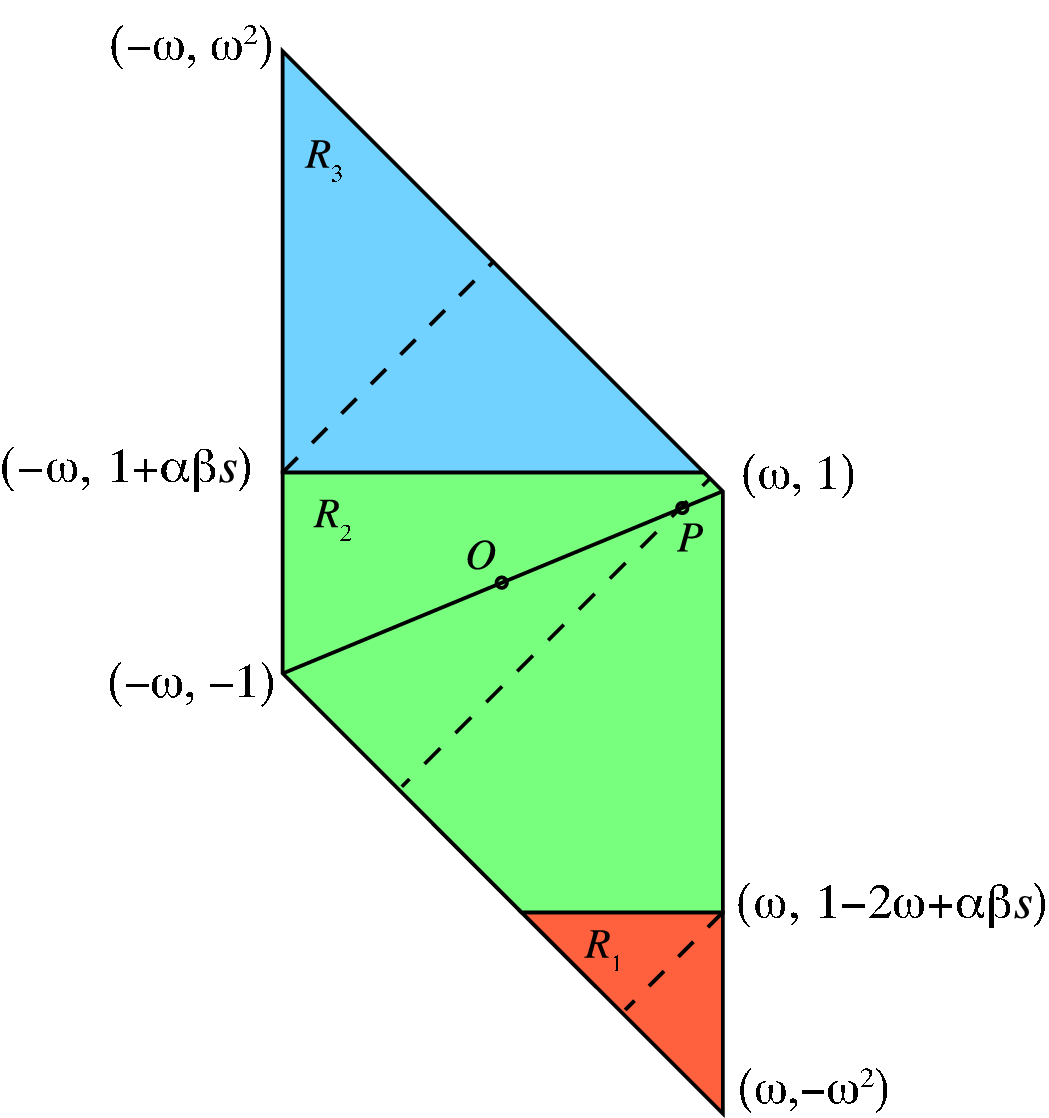,width=3 in}\hfil
\caption{\label{fig:Atoms} \small Partition of the rhombus $\mathrm{R}$ into atoms 
$R_i,\; i=1,2,3$.
The piecewise isometry $\rho(s)$ may be represented as a composition of two involutions: 
simultaneous reflection of the atoms about their respective symmetry axes (dashed lines), 
followed by reflection about the short diagonal of the rhombus (solid line).   
The intersection of these symmetry lines is a fixed point $P$, located $s$ 
to the left and $\beta s$ below the vertex $(\omega,1)$. }  
\end{figure}

We consider a one-parameter family $\rho(s)$ of piecewise isometries 
on a fixed rhombus $\mathrm{R}$ of side $2\alpha\omega$ and vertex angle $\pi/4$
(figure \ref{fig:Atoms}), specified by the half-plane conditions 
$\mathrm{u}_i\cdot (x,y)< \omega,\; i=0,1,4,5$, i.e.,
\begin{equation}\label{eq:R}
\mathrm{R}= [(0,1,4,5),(\omega,\omega,\omega,\omega)].
\end{equation}
The map $\rho(s)$ acts as an orientation-preserving isometry $\rho_i(s)\in \mathfrak{G}_+$ 
on each of its atoms: 
$R_i,\; i=1,2,3$,
\begin{eqnarray}
R_1(s)&=&[(0,2,5),(\omega, -1-2\alpha+\alpha\beta s,\omega)],\nonumber\\
R_2(s)&=&[(0,1,2,4,5,6),(\omega,\omega, 1+\alpha\beta s,\omega,\omega,1+2\alpha-\alpha\beta s)],\label{eq:Ri}\\
R_3(s)&=&[(1,4,6),(\omega,\omega,-1-\alpha\beta s)].\nonumber
\end{eqnarray}
Specifically, each atomic isometry is a clockwise rotation by $\pi/4$ followed by an $s$-dependent vertical translation, 
\begin{eqnarray}
\rho_1(s)&=&\mathtt{T}_{(0,2\beta-2\beta s)} \mathtt{R}_7,\nonumber \\
\rho_2(s)&=&\mathtt{T}_{(0,-2-2\beta s)} \mathtt{R}_7, \label{eq:rhoi}\\
\rho_3(s)&=&\mathtt{T}_{(0,-2\omega-2\beta s)} \mathtt{R}_7.\nonumber 
\end{eqnarray}
For each value of $s$ in the interval $[0,2\omega]$, the map $\rho(s)$ has a fixed-point
$$
P(s)= (\omega-s,1-\beta s)
$$
on the short diagonal of the rhombus.
The renormalizability for the cases $s=0,\omega,2\omega$ is known
\cite{AdlerKitchensTresser,KouptsovLowensteinVivaldi,AkiyamaBrunottePethoSteiner}.
%%%%%%%%% FIGURE 
%%%%%%%%%

The piecewise isometry $\rho(s)$ is {\it reversible}, namely it can be
written as the composition of two orientation-reversing involutions,
$$
\rho(s)= \mathtt{G}\circ \mathtt{H}(s),\qquad \mathtt{G}^2 = \mathtt{H}(s)^2 ={\bf 1},
$$ 
where $\mathtt{H}(s)$ is the simultaneous reflection of the three atoms about their respective symmetry axes, 
and $\mathtt{G}$ is the reflection of the rhombus $\mathrm{R}$ about its short diagonal.  
Note that the fixed point $P$ is symmetric, namely it lies at the intersection of fixed lines of $\mathtt{G}$ and $\mathtt{H}$.  
Moreover, $\mathtt{H}=\mathtt{G}\circ \rho$, and either $\mathtt{G}$ or $\mathtt{H}$ serves as a time-reversal operator for the map $\rho$:
$$
\mathtt{G}\circ \rho(s)\circ \mathtt{G} =  \mathtt{H}\circ \rho(s) \circ \mathtt{H} =  \rho(s)^{-1}.
$$

Another useful symmetry relation follows from the invariance of the rhombus under $\mathtt{R}_4$
(rotation by $\pi$), which takes $(x,y)$ into $(-x,-y)$.  One readily verifies that
$$
\mathtt{R}_4\circ \rho(s)\circ  \mathtt{R}_4= \rho(2\omega -s).
$$
In studying the renormalizability of the family over $[0,2\omega]$, we are thus permitted to 
restrict ourselves to $s\in [0,\omega]$. For reasons which will soon become clear, we will 
mainly be focusing on the shorter interval, 
\begin{equation}\label{eq:I}
I=[0,\alpha].
\end{equation}

Equipped with the $s$-dependent piecewise isometry $\rho(s)$, the rhombus $\mathrm{R}$ becomes the 
span of a parametric dressed domain $\mathcal{R}(s) = (\mathrm{R}, \mathbf{R}(s), \rho(s))$. 
To demonstrate the renormalizability of $\mathcal{R}(s)$, as defined in section 2, is the principal 
goal of this investigation.  To do this, we concentrate on the atom $R_1(s)$, showing that it is a 
dressed subdomain of $\mathcal{R}(s)$ and moreover is an example of a two-parameter family of 
\textit{base triangles}.  
The renormalizability of base triangles will then occupy our efforts for the remainder of the article.
%

%%%%%SUBSECTION
\subsection{The base triangle}
We define a two-parameter parametric dressed domain, the \textit{base triangle}.
For parameters $l\in \R_+$ and $h\in [0,\alpha l]$, we define a prototype
$\widehat{\cB}(l,h)=(\widehat{\mathrm{B}},\widehat{\mathbf{B}},\rho)$ to represent 
its equivalence class with respect to $\mathfrak{G}$.  
The dressed domain induced on the atom $R_1(s)$ of $\mathcal{R}$ will be shown below 
to be congruent to $\widehat{\cB}(1,s)$.

\begin{figure}[h]
\hfil\epsfig{file=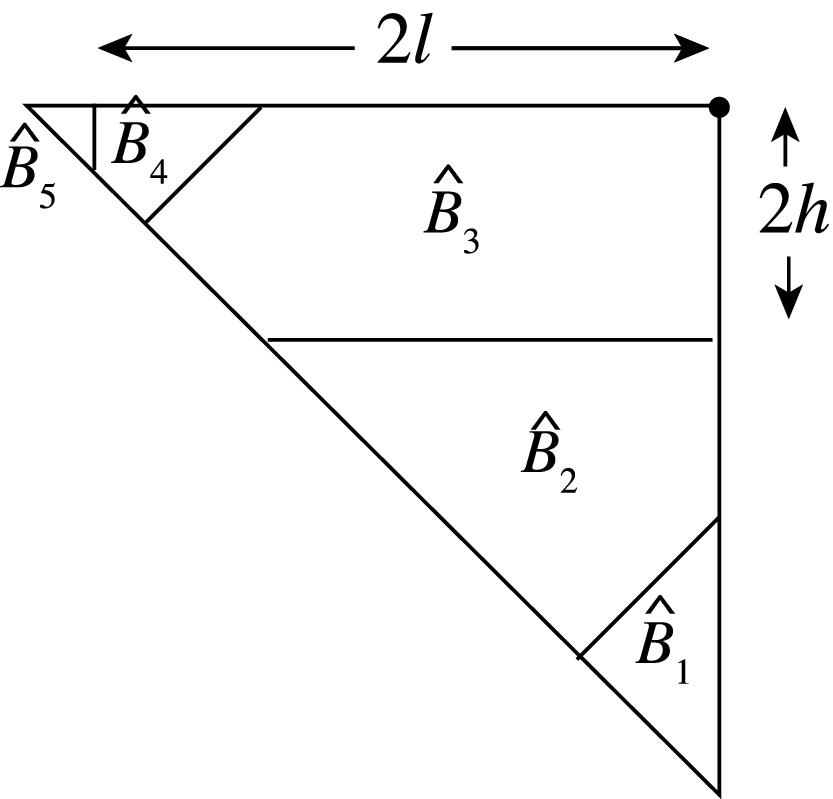,width=2 in}\hfil
\caption{\label{fig:protoB} \small The prototype base triangle $\widehat{\cB}(l,h)$ for $0<h/l<\alpha$. }
\end{figure}

The tiling $\widehat{\mathbf{B}}$ is illustrated in figure \ref{fig:protoB}.
The defining data are presented in Table \ref{tbl:Bdef}.
For simplicity, we choose a frame of reference with the right-angle vertex of
$\widehat{\mathrm{B}}$ at the origin and the remaining vertices at points of the 
negative $x$ and $y$ axes.

%%%%%%%%%%%TABLE
\begin{table}[h]
\caption{\label{tbl:Bdef}\small Tiling table of the prototype base triangle $\widehat{\cB}(l,h)$, for $0<h/l<\alpha$. }
$$
\begin{array}{|c|c|c|c|c|c|c|}\hline
\multicolumn{3}{|c|}{\mbox{Source Polygon}}&\multicolumn{2}{|c|}{\mbox{Initial Placement}}&\multicolumn{2}{|c|}{\mbox{Destination}}\\ \hline
\mbox{Tile}&Q_{\#}&\mbox{Parameters}&{\tt R}_{\#}&\mbox{Translation}&\mbox{Tile}&{\tt R}_{\#}\\ \hline\hline
\widehat{\mathrm{B}}&1&\alpha l+\beta h&7&(0,0)&\mbox{---}&\mbox{---}\\ \hline \hline
\widehat{B}_1 & 1&\alpha \beta l-\beta h & 2 & (-\alpha l - h,-\alpha \beta l+ \beta h)&\mathtt{U}_1 \widehat{B}_1&2\\ \hline
\widehat{B}_2 & 7&2 l - \alpha h & 0 & (-2 l+\alpha h,-2 h)&\mathtt{U}_1\widehat{B}_2&3\\ \hline
\widehat{B}_3 & 6&\alpha l-h,h   &  0 & (-h,-h)&\mathtt{U}_1\widehat{B}_3&2\\ \hline
\widehat{B}_4 & 7&\alpha h & 5   & (-2 l+\alpha h,0)&\mathtt{U}_1\widehat{B}_4&1\\ \hline
\widehat{B}_5 & 1&\beta h & 7    & (-2 l,0)&\mathtt{U}_1\widehat{B}_5&0\\ \hline
\end{array}
$$
\end{table}
%%%%%%%%%%%%%%%%
In the table, the atoms of $\widehat{\mathbf B}$ and their span are listed by giving their respective 
orientations and translation vectors relative to a representative tile in the catalogue of 
Appendix A. For example, we learn from table \ref{tbl:Bdef} that 
$$
\widehat{B}_1=\mathtt{T}_{\mathbf{d}}\,\mathtt{R}_2\, Q_1(\alpha \beta l-\beta h),
$$
where $\mathbf{d}=(-\alpha l-h,-\alpha \beta l+\beta h)$ is the location of the offset tile's anchor point (local origin; 
see Appendix A).  
The isometry $\rho_1$ associated with atom $\widehat{B}_1$ is uniquely specified by the information listed in the last two columns.  Because $\mathbf{d}$ lies on the symmetry line of the tile, it is taken into $\mathtt{U}_1(\mathbf{d})$ by $\rho_1$. More generally, if $n_1$ is the rotation index of the last column of the table, we calculate
$$
\rho_1= T_{\mathtt{U}_1\mathbf{d}-\mathtt{R}_{n_1}\mathbf{d}}\, \mathtt{R}_{n_1}.
$$

%%%%%%%%% FIGURE 
%%%%%%%%%%%%%%

The definition of the base triangle can be extended to the boundary of the parametric domain.
This domain has two atoms for $h=0$ (and any $l$) and three for $h=\alpha l$; the atoms 
are still described by table \ref{tbl:Bdef} with the stipulation that all zero-parameter 
entries are to be deleted.

The following result establishes the dynamics of the base triangle, which is the basis of
the renormalization process.
\medskip

\begin{proposition}\label{prop:BaseTriangle}
For all $s\in I$, let $\mathcal{R}(s)=(\mathrm{R}(s),(R_i(s)),(\rho_i(s)))$ be the parametric dressed 
rhombus defined in equations (\ref{eq:R})--(\ref{eq:rhoi}).  
Then the atom $R_1(s)$, equipped with the return map induced by $\rho(s)$, is a dressed subdomain 
$\mathcal{B}$ congruent to the prototype base triangle $\widehat{\mathcal{B}}(1,s)$. 
The rhombus $\mathrm{R}(s)$ is tiled by the return orbits of the atoms of $\mathcal{B}$ and 
also of the periodic tiles
$$
\begin{array}{l}
 \Pi_1=\mathtt{T}_{(-1-s,1-\beta s)}\, Q_{10}(s,\alpha-s),\\
 \Pi_2=\mathtt{T}_{(-1-s,-\beta-\beta s)}\, Q_5(\alpha-s),\\
 \Pi_3=\mathtt{T}_{(\omega-s,1-\beta s)}\, Q_5(s),
 \end{array}
$$
apart from a set of zero measure. The incidence matrix for the return orbits of the atoms is:  
$$
\left(
\begin{array}{ccccc}
1&1&1&1&1\\
8&5&2&4&6\\
5&2&2&2&2
\end{array}
\right).
$$
\end{proposition}

\noindent {\sc Proof.} The proof is a straightforward application of the 
\textit{direct iteration method} described in section \ref{section:DirectIteration}.
The initial and final tiles of the return orbits can be gleaned from table \ref{tbl:Bdef}, 
and we know that the periodic orbits should begin and end on $\Pi_k,\, k=1,2,3$.  
In the course of constructing the return orbits, we keep track of the atoms visited, obtaining 
at the end the symbolic paths and incidence matrices of the orbits. By adding up the areas of 
the tiles of all the orbits, and comparing 
with the total area of the parent domain, we prove the completeness of the tiling. 

The details of the computer-assisted calculation may be found in the 
Electronic Supplement\cite{ESupplement}.  \hfill $\Box$

\begin{figure}[h]
\hfil\epsfig{file=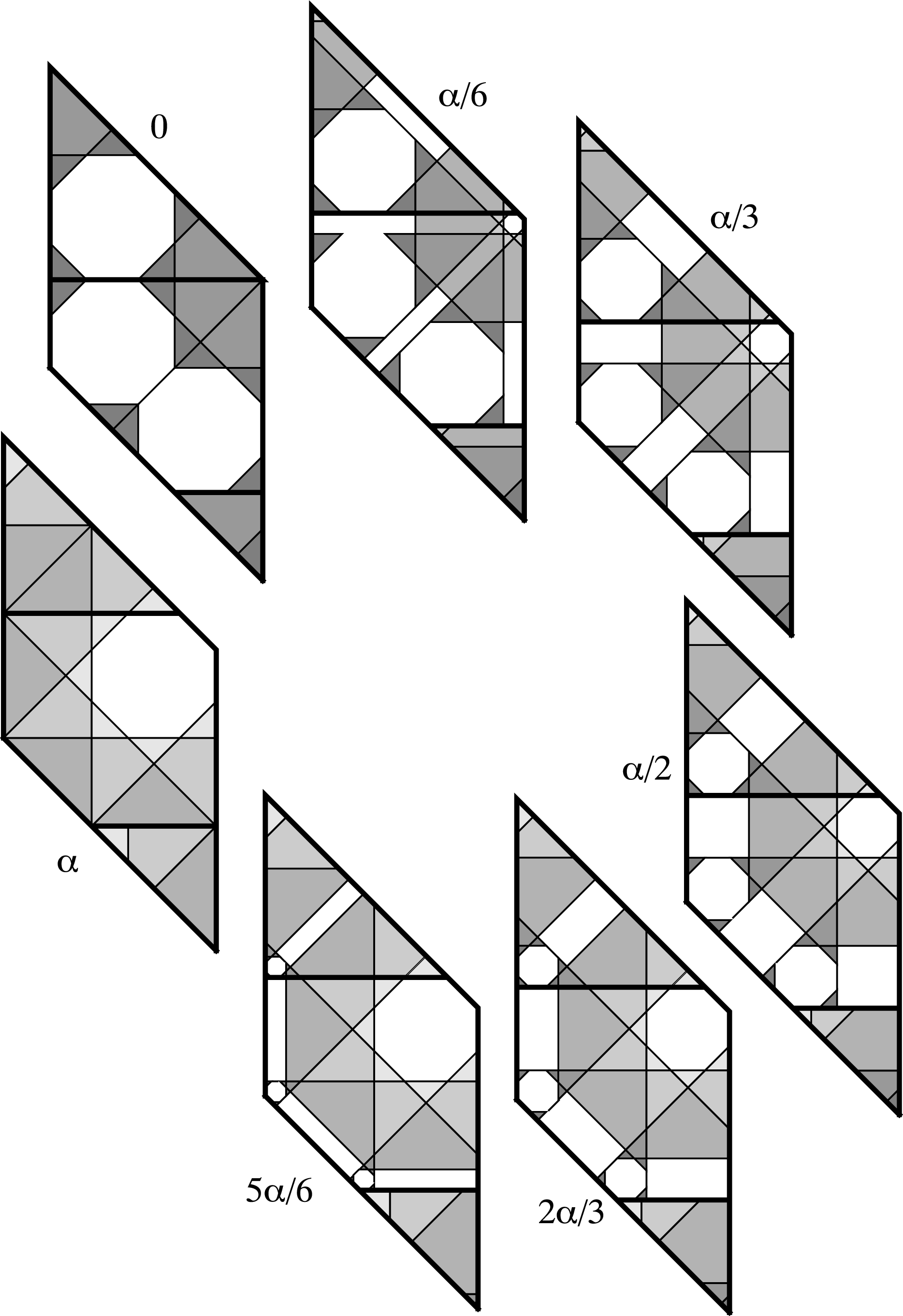,width=12cm}\hfil
\caption{\label{fig:BtilingR} \small Tiling of the rhombus $\mathrm{R}$ by return orbits of the atoms 
of $R_1$ (shades of gray) and the periodic tiles $\Pi_1,\Pi_2,\Pi_3$ (white), for equally 
spaced values of $s=h/l$ in $I=[0,\alpha]$. 
The three atoms of $\mathcal{R}$ have been drawn with thicker boundaries.}  
\end{figure}

The tiling is illustrated for several values of $s$ in figure \ref{fig:BtilingR}. 
The reader may find it instructive to follow each of the orbits around the rhombus, 
applying `by eye' the product of local and global involutions at each step. 

%%%%%%%%% FIGURE 
%%%%%%%%%
\clearpage

%
%%%%SECTION
\section{Main results}\label{section:MainResults}
\begin{figure}[h]
\hfil\epsfig{file=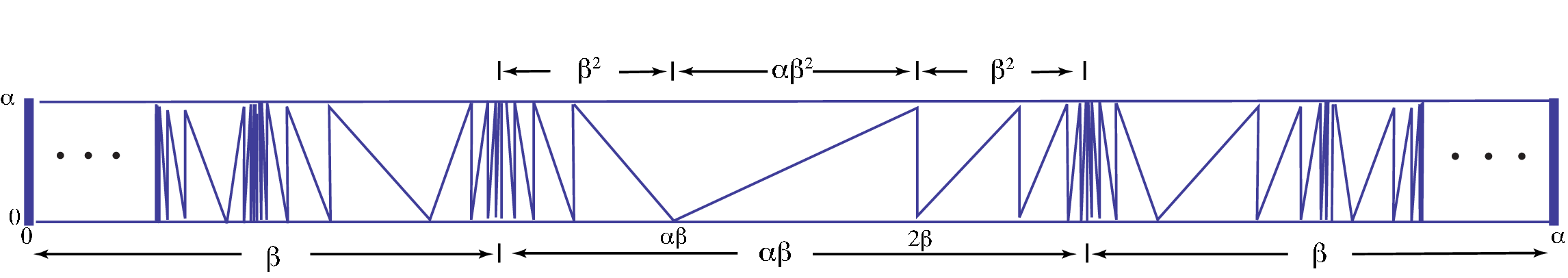,width=5in}\hfil
\caption{\label{fig:r}\small
The piecewise-affine renormalization function $r=f^2$}.  
\end{figure}

We now state the main results of our investigation.
We begin by defining the \textit{renormalization functions} 
$f$ and $r=f^2$ (see figures \ref{fig:f} and \ref{fig:r}), for all $s\in I=[0,\alpha]$:
$$
f(0)=f(\alpha)=0,
$$
\begin{equation}\label{eq:f}
f(s)= \omega^{|i|+1}\left\{\begin{array}{ll} (\Delta_i-s) 
   & s\in I_i,\; i<0\\ (s-\beta)& s\in I_0\\ (s-\Delta_i) &s\in I_i,\;  i>0, \end{array}\right.
\end{equation}
where
\begin{equation}\label{eq:IDelta}
I_i=\left\{\begin{array}{ll} (\Delta_{i-1},\Delta_i]& i<0\\
(\beta,1)& i=0\\   
\left[\Delta_i,\Delta_{i+1})\right. & i> 0 \end{array}\right.\qquad
\Delta_i=\left\{\begin{array}{ll}\beta^{|i|}& i<0\\ \beta,&i=0\\ \alpha-\beta^{|i|} & i> 0\end{array}\right.
\end{equation}
and
$$
r(0)=r(\alpha)=0,
$$
\beq\label{eq:r(s)}
r(s)=\omega^{|i|+|j|+2}\left\{\begin{array}{ll}
(\Delta_{i,j}-s) &s\in I_{i,j},\quad  j<0 \mbox{ or } (j=0 \mbox{ and } i<0),\\ 
(s-\Delta_{i,j}) &s\in I_{i,j},\quad  j> 0 \mbox{ or } (j=0 \mbox{ and } i\geq0), \end{array}\right.
\eeq
where
$$
I_{i,j}=\left\{\begin{array}{ll} 
(\Delta_{i,j-1},\Delta_{i,j}]& j<0 \\
(\Delta_{i,-1},\Delta_{i,0}) & j=0,\,i<0\\ 
(\Delta_{i,0},\Delta_{i,1})  & j=0,\,i\geqslant 0\\
 \left[\Delta_{i,j},\Delta_{i,j+1})\right. &j>0,
 \end{array}\right.
$$
\medskip
$$
\Delta_{i,j}=\left\{\begin{array}{ll}
\beta^{|i|+1}+\beta^{|i|+|j|+1}& i\leqslant 0,\, j\leqslant 0\\ 
\beta^{|i|}-\beta^{|i|+|j|+1}  & i<0,\, j>0\\
\alpha-\Delta_{-i,-j}          & i>0 \mbox{ or } (i=0 \mbox{ and } j>0).
 \end{array}\right.
$$

%%%%%%%%%%%%%%%%%%%%%%%%%%%%%%%%%%%%%%%%%%%%%%%% FIGURE 
%%%%%%%%%%%%%%%%%%%%%%%%%%%%%%%%%%%%%%%%%%%%

The following two theorems constitute our main results. \pagebreak
\begin{theorem}\label{thm:Main} {\sl
Let $\cB$ be a parametric base triangle congruent to the prototype $\widehat{\cB}(l,s l)$. 
Then $\cB$ is renormalizable for all positive real $l$ and $s\in I$.  Specifically: 
\medskip

\noindent {$i)$} For $s\in\{0, \alpha\}$, $\cB$ has a dressed subdomain $\cB^*$ congruent 
to $\widehat{\cB}(\kappa(s) l,0)$, with
$$
\kappa(0)= \beta, \qquad 
\kappa(\alpha)= 1 .
$$
The parity $\pi(\cB^*)$ is $-\pi(\cB)$ for $s=0$ and $\pi(\cB)$ for $s=\alpha$. 
The return-map orbits of the atoms of $\cB^*$, together with those of a finite number of periodic 
tiles,  tile the spanning triangle of $\cB$, up to a set of measure zero.\bigskip

\noindent {$ii)$} For all $i\in\Z\setminus\{0\}$ and $s=\Delta_i$, as defined in (\ref{eq:IDelta}), the
domain $\cB$ has two dressed subdomains, $\cB^*$ and $\cB^\dagger$ congruent, respectively, 
to $\widehat{\cB}(\kappa^*(\Delta_i)l,0)$ and $\widehat{\cB}(\kappa^\dagger(\Delta_i)l,0)$, with
$$
\kappa^*(\Delta_i)= \beta^{|i|}/\alpha , \qquad 
\kappa^\dagger(\Delta_i)=  \beta^{|i|+2}/\alpha. 
$$
The parities $\pi(\cB^*)$ and $\pi(\cB^\dagger)$  are both $(-1)^{|i|+1}\pi(\cB)$.
The return-map orbits of the atoms of $\cB^*$ and $\cB^\dagger$, together with those of a finite 
number of periodic tiles,  tile the spanning triangle of $\cB$, up to a set of measure zero.
\bigskip

\noindent {$iii)$} For all $i,j\in\Z$ and $s\in I_{i,j}$, the dressed domain $\cB$ has a dressed subdomain 
$\cB^*$ congruent to $\widehat{\cB}(l^*, s^* l^*),\;l^*=\kappa(s)l,\; s^*=r(s)$,  with $r(s)$ given by (\ref{eq:r(s)}) and
$$
\kappa(s)=\beta^{|i|+|j|+2},\qquad
\pi(\cB^*)= (-1)^{|i|+|j|}\pi(\cB).
$$
The return-map orbits of the atoms of $\cB^*$, together with those of a finite number of periodic tiles,  
tile the spanning triangle of $\cB$, up to a set of measure zero.  
The tilings vary continuously with respect to $s\in I_{i,j}$, with the return paths 
(hence the incidence matrix) constant over the interior of the interval.}
\end{theorem}
\medskip

There is a tight connection between strict renormalizability and the field $\Q(\sqrt{2})$, 
due to the following result.

\begin{theorem}\label{thm:f}
{\sl The renormalization function $f$ is conjugate to a left shift acting in a space of one-sided 
symbol sequences with alphabet $\Z\cup\{-\infty,+\infty\}$. 
A point $s\in I$ is eventually periodic under $f$ if and only if $s\in\Q(\sqrt{2})$.  
Hence the set of values of $s$ for which a base triangle congruent to $\widehat{\cB}(l,s l)$ is 
strictly renormalizable is $\Q(\sqrt{2})\cap I$. }
\end{theorem}

Theorem \ref{thm:Main} will be proved from the analysis of the return-map dynamics. 
Before that rather lengthy analysis, we will prove theorem \ref{thm:f} and study some 
interesting properties of the function $f$ and its symbolic dynamics.
%
%---------------------------------------------------------------------------------- 
%SECTION
\section{Symbolic dynamics}\label{section:SymbolicDynamics}
In this section we introduce a symbolic dynamics for $f$ which will give us a 
useful expansion for $s\in I$ ---equation (\ref{eq:Luroth}).
This is a variant of the so-called L\"uroth expansions \cite{Luroth,Galambos,BarrionuevoEtAl}. 

The renormalization map $f$, see (\ref{eq:f}) and figure \ref{fig:f}, is piecewise-affine, and
its restriction to the interval $I_i$ has slope $\pm\omega^{|i|+1}$. 
Since $\omega>1$, this map is expanding, and since the length 
of $I_i^\pm$ is equal to $\alpha\beta^{|i|+1}=\alpha/\omega^{|i|+1}$,
we see that $f(I_i)=I$ (with the origin missing if $i=0$). 
It follows that $f$ preserves the Lebesgue measure.
%The Lyapounov exponent is 
%$$
%\beta\log \omega+2\sum_{n=1}^\infty \beta^{n+1}\log\omega^{n+1}=2\log\omega.
%$$

Next we define the function $i:I\to \Z\cup\{-\infty,+\infty\}$ which assigns to each $s$ 
the interval $I_i$ to which it belongs. Explicitly,
\beq
i(s)=\left\{\begin{array}{ll}
       \lfloor\log_\beta(s)\rfloor & \qquad 0<s\leqslant\beta\\
       \lfloor\log_\beta(\alpha-s)\rfloor & \qquad \beta < s<\alpha\\
        -\infty &\qquad s=0\\
        +\infty &\qquad s=\alpha.
\end{array}\right. \label{eq:i}\\
\eeq
In terms of this function, $f$ can be rewritten as
\begin{equation}\label{eq:f2}
f(s)  =\left\{\begin{array}{ll}
        \sigma(i(s))\,\omega^{|i(s)|+1}\,\bigl(s-\Delta_{i(s)})
         &\quad s\not=0, \alpha\\
        0&\quad s=0, \alpha,
        \end{array}\right.
\end{equation}
where
\beq\label{eq:sigma}
\sigma(i)=\left\{\begin{array}{ll}-1 & i<0\\ 1 & i\geqslant 0.\end{array}\right.
\eeq

  Let us now consider the orbit of $f$ with initial condition $s_0=s\in I$.
For all $k\geqslant 0$ we let $s_{k+1}=f(s_k)$. With the notation 
\beq
i_k=i(s_k), \qquad
\sigma_k=\sigma(i_k),\qquad
 \pi_k=\alpha(1+\sigma_k)/2,
 \eeq
equation (\ref{eq:f2}) becomes 
\begin{equation} \label{eq:s}
s_{k+1}=
 \left\{\begin{array}{ll}
         \sigma_k\omega^{|i_k|+1}\,(s_k-\pi_k+\sigma_k\beta^{|i_k|})&\quad s_k\not=0,\alpha\\
         0                                                  &\quad s_k=0,\alpha.
        \end{array}\right.
\end{equation}

In this way, every $s\in I$ is associated with a unique sequence
\begin{equation}\label{eq:SymbolicDynamics}
(i_0,i_1,i_2, \ldots\,)
\end{equation}
with $i_k\in\mathbb{Z}\cup\{-\infty,+\infty\}$. The only constraints which we impose on 
these sequences is that the sub-strings $(\mp\infty,r),\;r\neq -\infty$ are forbidden,
and that the symbol $+\infty$ can only appear as the leading symbol: $(+\infty,-\infty,\ldots)$.
In the space of such sequences, a left shift is conjugate to the map $f$ on $I$.

To establish that each allowed sequence (\ref{eq:SymbolicDynamics}) corresponds to a unique $s$, 
we begin with the sequences $(\mp\infty,-\infty,-\infty,\ldots)$, which represents exclusively 
the value $0$ (negative sign) and $\alpha$ (positive sign).  
All other sequences either have no symbols $-\infty$, or else have an infinite tail of 
$-\infty$ symbols preceded by a finite sequence in which $\mp\infty$  does not appear.  
In either case, we can assume $s_k\neq 0,\alpha$ and invert (\ref{eq:s}) to 
obtain\begin{equation}\label{eq:finv}
s_{k}=\pi_k-\sigma_k\beta^{|i_k|}+\sigma_k\beta^{|i_k|+1}s_{k+1}.
\end{equation}
Iterating, we find:
\begin{eqnarray*}
s_{0}&=&\pi_0-\sigma_0\beta^{|i_0|}+\sigma_0\beta^{|i_0|+1}s_{1}\\
&=&\pi_0-\sigma_0(\omega-\pi_1)\beta^{|i_0|+1}-\sigma_0\sigma_1\beta^{|i_0|+|i_1|+1}
           +\sigma_0\sigma_1\beta^{|i_0|+|i_1|+2}s_2.
\end{eqnarray*}
An easy induction gives
\begin{equation}\label{eq:Luroth}
s=s_0=\pi_0-\Lambda
\hskip 40pt
\Lambda=\sum_{n\geqslant 1}a_n\beta^{n+b_n},
\end{equation}
where 
\begin{equation}\label{eq:abc}
a_n=(\omega-\pi_n)\prod_{k<n}\sigma_k
\hskip 40pt
b_n=\sum_{k<n}|i_k|.
\end{equation}

The sequence $(b_n)$ is non-negative and non-decreasing, and we have 
$b_n\equiv 0$ only for the fixed point $s=\sqrt{2}/2$. 
The sequence $(a_n)$ depends only on the $\sigma_k$s; indeed, 
\begin{equation}
\omega-\pi_n=\left\{\begin{array}{ll}
        \omega&\quad \sigma_n=-1\\
        1 & \quad\sigma_n=+1.
        \end{array}\right.
\end{equation}
Thus $a_n\in\{\pm1,\pm\omega\}$.

The sum in (\ref{eq:Luroth}) is absolutely convergent, and it provides 
an expansion for any $s\in[0,\alpha)$. On the other hand, having excluded
$+\infty$ from all but one sequence, distinct sequences correspond to
distinct values of $s\in[0,\alpha]$.
This completes the proof of the claimed bi-unique correspondence.
 
The expansion (\ref{eq:Luroth}) is finite if the orbit of $s_0$ reaches the origin, 
and infinite otherwise. In the former case $\Lambda$ is a finite sum of 
elements of $\mathbb{Z}[\alpha]$, and hence $s\in\mathbb{Z}[\alpha]$.
If the sequence (\ref{eq:SymbolicDynamics}) is eventually periodic
with limiting period $N$, then the sequence $a_n$ is eventually periodic 
with the same transient and period $N$ or $2N$, while the sequence $b_n$ 
eventually becomes the sum of an affine function plus a periodic function 
with period dividing $N$. Then the sum $\Lambda$ decomposes into
the sum of finitely many geometric progressions, and so
$\Lambda$, and hence $s$, belong to $\mathbb{Q}(\alpha)$.

In the next section we shall demonstrate that the converse is also true, namely 
that any $s\in\mathbb{Q}(\alpha)\cap [0,\alpha]$ has an eventually periodic 
symbol sequence of the type (\ref{eq:SymbolicDynamics}).

%----------------------------------------------------------------------------------
\section{Lattice dynamics}\label{section:LatticeDynamics}
Let the ring $\mathbb{Z}[\alpha]$ and the interval $I$ be given 
by (\ref{eq:FieldRing}) and (\ref{eq:I}), respectively.
For $d\geqslant 1$, we define
\begin{equation}\label{eq:M}
\mathcal{M}_d=\frac{1}{d}\mathbb{Z}[\alpha]\cap I
\end{equation}
which is the restriction to $I$ of the module $d^{-1}\mathbb{Z}[\alpha]$.
Because $f(s)$ is obtained from $s$ by ring operations in
$\mathbb{Z}[\alpha]$, and $f(I)=I$, we have that $f(\mathcal{M}_d)\subset\mathcal{M}_d$. 
We have established the first part of the following lemma:
%%%%%%LEMMA
\begin{lemma}\label{lemma:ff}
{\sl For each $d$, we have}
$$
f(\mathcal{M}_d)\subset \mathcal{M}_d,\qquad
f^{-1}(\mathcal{M}_d\setminus \{\alpha\} )\subset \mathcal{M}_d.
$$
\end{lemma}
%%%%%%%%%%%

\noindent {\sc Proof.} It remains to show the invariance with respect to $f^{-1}$. 
In the statement of the lemma, the element $\alpha\in \mathcal{M}_1$ had 
to be removed since it is not in the domain of the function $f^{-1}$.   

We have $0\in \mathcal{M}_1$, and, by construction we have
$$
f^{-1}(\{0\})=\{\ldots,\Delta_{-2},\Delta_{1},\Delta_0,\Delta_1,\Delta_2,\ldots\}\cup \{0,\alpha\}
    \subset \mathcal{M}_1.
$$
Let now $s\not=0,\alpha$ and let $s^\prime=f(s)\in\mathcal{M}_d$.
Using (\ref{eq:f}) we find
$$
s=s^\prime\sigma(i(s))\beta^{|i(s)|+1}+\Delta_{i(s)}.
$$
For any choice of the values of $i$,
we see that $s$ is an affine function of $s^\prime$ with coefficients
in $\mathbb{Z}[\alpha]$. Since $\mathbb{Z}[\alpha]/d$ is a module
over $\mathbb{Z}[\alpha]$, it follows that $s\in\mathcal{M}_d$ (see (\ref{eq:M})).
\hfill $\Box$

\bigskip

Since, by the lemma, the inverse function cannot increase the denominator, 
then the forward function cannot decrease it. This can be rephrased as follows.
For any $\xi\in\mathbb{Q}(\alpha)$, let $d^*=d^*(\xi)$ be the smallest
natural number $d$ such that $\xi\in\mathcal{M}_{d}$. 
Then $d^*$ is a constant of the motion for the map $f$.

For $s\in\mathcal{M}_d$ we let $\zeta=s d$.
Then $\zeta\in\mathbb{Z}[\alpha] \cap dI$, and alongside the interval map 
$s\mapsto f(s)$, we have the ring map 
\begin{equation}\label{eq:fd}
f_d:\mathbb{Z}[\alpha]\to\mathbb{Z}[\alpha]\qquad
\zeta \mapsto \sigma(i(s)) \omega^{|i(s)|+1} 
     \left(\zeta-d \Delta_{i(s)}\right)
\quad
s=\frac{1}{d}\,\zeta
\end{equation}
(with $f_d(0)=f_d(d\alpha)=0$) which represents the restriction of $f$ to 
$\mathcal{M}_d$ after clearing denominators. 
\bigskip

\noindent {\sc Conclusion of the proof of theorem 2.}
We introduce the natural bijection
$$
\phi:\mathbb{Z}[\alpha]\to\mathbb{Z}^2\qquad
m+n\alpha\mapsto (m,n)
$$
which conjugates $f_d$ to a lattice map on $\Z^2$, for which still 
use the same symbol.

For any $d> 0$, we define the infinite strip 
$$
\Sigma_d=\left\{(m,n)\in\mathbb{Z}^2\,:\,
  -m\leqslant n\alpha \leqslant -m+d\alpha\right\},
$$
which is invariant under $f_d$ (because $I$ is invariant under $f$).

We claim that all orbits of the map $f_d$ are eventually periodic.
Since $\Sigma_d\subset \mathbb{Z}^2$ this means that all orbits
of $f_d$ are bounded.
Multiplication by $\omega$ in $\mathbb{Z}[\alpha]$ induces
a linear map $M$ on $\mathbb{Z}^2$, with eigenvalues $-\omega,\beta$.
The lines $\alpha y=\pm x$ are the corresponding eigendirections.
Let $z\in\Sigma_d$, and let $(z_-,z_+)$ be the 
components of $z$ with respect to a basis of eigenvectors.

Since the expanding eigendirection is transversal to $\Sigma_d$,
there is a constant $c_d$ such that $|z_+|<c_d$. So it suffices
to show that the component $z_-$ remains bounded.
For all $z\in\mathbb{Z}^2$ and $i\geqslant 0$ we have
\begin{equation}\label{eq:FirstBound}
|M^{j+1}(z)_-|\leqslant \beta |z_-|<|z_-|.
\end{equation}
Furthermore, from (\ref{eq:IDelta}) we have 
that $\omega^{|i|+1}\Delta_{i}$ is a monomial or binomial
in $\omega$ of degree at most ${|i|+1}$ with coefficients $1,-1,\alpha$.
Defining
$$
C=d\,(|\phi(1)_-| + |\phi(-1)_-| + |\phi(\alpha)_-|)
$$
from (\ref{eq:FirstBound}) we have
\begin{equation}\label{eq:SecondBound}
|M^{|i|+1}(\phi(d\Delta_i))_{-}|\leqslant |\phi(d\Delta_i)_{-}|\leqslant C.
\end{equation}
Finally, from (\ref{eq:fd}), (\ref{eq:FirstBound}), and 
(\ref{eq:SecondBound}), we obtain
$$
|f_d(z)_-|\leqslant \beta|z_-|+C
$$
and since $\beta<1$, for all sufficiently large
$|z_-|$ the map $f_d$ is a contraction mapping. Thus its
orbits are bounded hence eventually periodic.

This completes the proof of theorem 2.\hfill  $\Box$

%%%%%%%%%%%% SECTION
\section{Overview of the renormalization dynamics}\label{section:Renormalization}
As a prelude to the proof of theorem 1, we now turn our attention to the dynamical 
underpinnings of the renormalizability of the parametric base triangle $\cB(l,s)$. 
This analysis is based on the construction of the return-map tree through successive 
inductions on sub-domains, a process which is far more complex than one might guess from the 
simple functional form of the function $r$.
To account for the renormalizability of the entire parameter interval, 
ten distinct \textit{renormalization scenarios} need to be considered, 
each characterized by the participation of distinctive parametric dressed domains.  
We have given such special domains names suggestive of their geometric structure: 
the {\it pencil} $\cP$, the {\it fringed triangle} $\cT$, and the {\it double strip} $\cD$.  
The induction sequence for each of the ten scenarios and the corresponding
parameter intervals are specified in Table \ref{tbl:scenarios}.  
In the labelling of the scenarios, Roman numerals I through IV are used to indicate 
the major classifications, with asterisks and binary subscripts $\mu,\nu=\pm1$ 
indicating finer distinctions (to be clarified in the next section). 
%%%%%TABLE
\begin{table}[h]
\caption{\label{tbl:scenarios}\small Renormalization scenarios, each one corresponding to 
a simple closed loop on the renormalization graph of figure \ref{fig:RenormGraph}.  
(See also figures \ref{fig:sAxis} and \ref{fig:4Scenarios}.)}
$$
\begin{array}{|c|c|c|}
\hline
\multicolumn{2}{|c|}{\mbox{Renormalization scenario}} & \mbox{Ranges of indices}\\
\hline
\mbox{I}&{\cB}\rightarrow {\cB}&\begin{array}{c}
(\pm \infty,0),\,(0,\pm \infty)\\
(-1,j),\;(1,-j),\quad j=0,1,2,3\\
(0,j),\quad |j|\leqslant 2.
\end{array}
\\
\hline
\mbox{II}&{\cB}\rightarrow\cP\rightarrow\cP^*\rightarrow  {\cB}&\begin{array}{c}
(i,j), (-i,-j),\;\;i\geqslant 2,\;\; -3\leqslant j\leqslant 2,\\
(1,j),\; (-1,-j),\quad j=1,2,\\
(-i,-\infty),\, (i,+\infty),\quad i\geqslant 1.
\end{array}
\\
\hline
\mbox{III}_{\mu\nu}&{\cB}\rightarrow\cT_\mu\rightarrow\cD_\mu\rightarrow \cD^*_\nu\rightarrow  {\cB}&
\begin{array}{l}
\mu=-1 :\quad  (0,\pm j), \quad j\geqslant 3,\\
\mu=+1 :\quad (1,-j), \; (-1,j), \quad j\geqslant 4,\\
\nu=(-1)^j.
\end{array}
\\
\hline
\mbox{IV}_{\mu\nu}&{\cB}\rightarrow\cP\rightarrow\cT_\mu \rightarrow\cD_\mu\rightarrow \cD^*_\nu \rightarrow  {\cB}&
\begin{array}{l}
\mu=-1 :\quad (i,j),\;  (-i,-j),\; i\geqslant 1,\;\; j\geqslant 3,\\
\mu=+1:\quad (-i,j),\;  (i,-j),\; i\geqslant 2,\;\; j\geqslant 4,\\
\nu=(-1)^j.
\end{array}
\\
\hline
\end{array}
$$
\end{table} 

\begin{figure}[h]
\hfil\epsfig{file=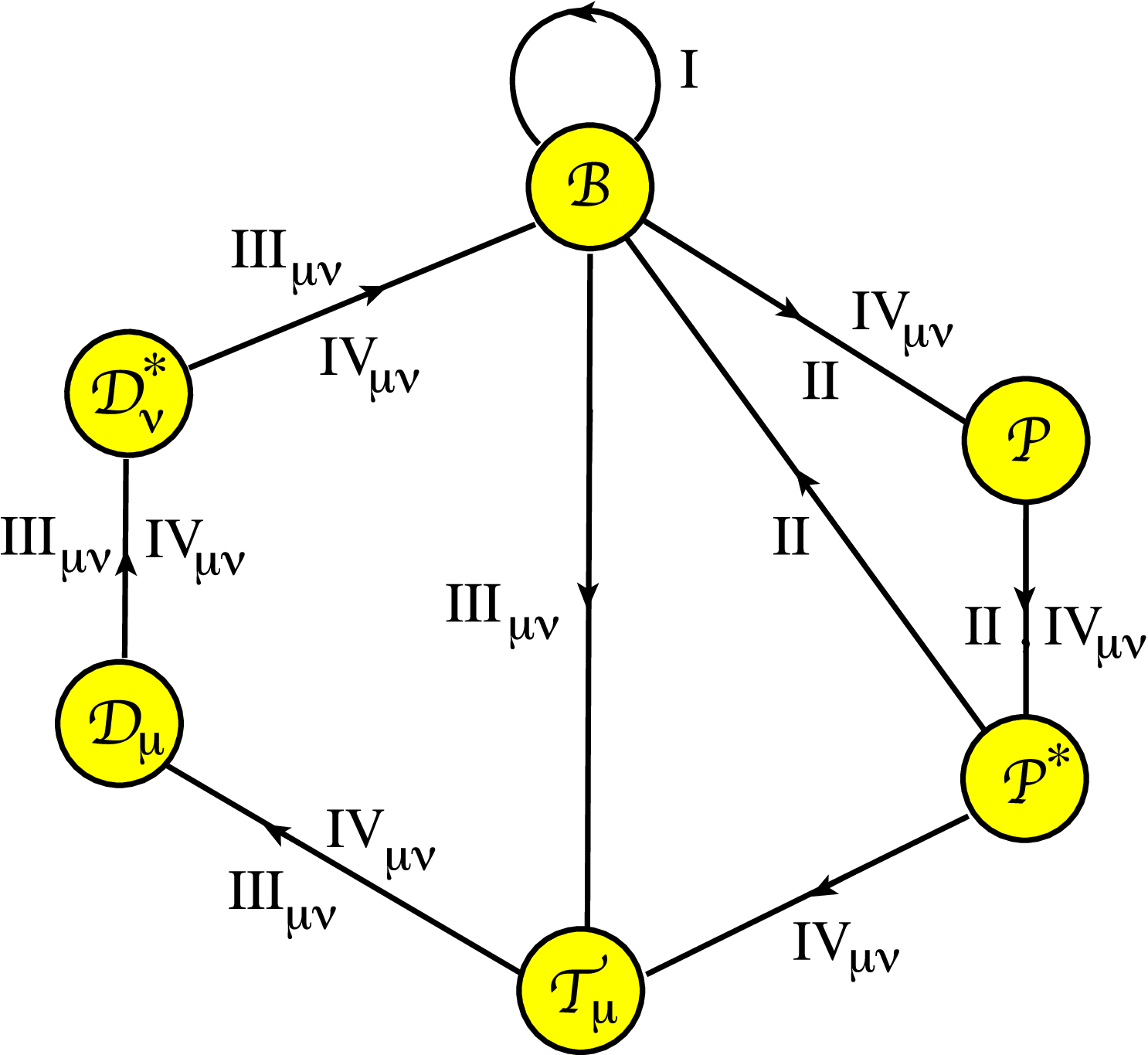,width=9cm}\hfil
\caption{\label{fig:RenormGraph} \small Renormalization graph, whose closed loops 
are the renormalization scenarios of Table \ref{tbl:scenarios}.}
\end{figure}

The rest of this section is devoted to graphical representations of the ten scenarios,
the most important of which is the \textit{renormalization graph} of figure \ref{fig:RenormGraph}. 
Each vertex of the graph corresponds to an equivalence class of parametric dressed domains. 
Thus the vertex $\cB$ represents a base triangle congruent to the prototype $\widehat{\cB}(l,h)$, 
with $l>0$ and $s=h/l\in I$. 
The precise interpretation of the remaining vertices will emerge from the prototype 
definitions and lemmas of section 8, together with the specification of parameter ranges in 
the propositions of section 10.
An oriented edge of the graph, $\mathcal{X} \rightarrow \mathcal{Y}$,
signifies that ${\cal Y}$ is a dressed sub-domain of ${\cal X}$,
which is dressed by $\mathcal{X}$ via induction.  
Each edge is labelled by subscripted Roman numerals, indicating 
the relevant parameter constraints listed in Table \ref{tbl:scenarios}. 
Loops in the graph correspond to renormalization scenarios and  
since $\mu,\nu=\pm1$, there are ten different scenarios in all.

In labelling the vertices of the graph, we have used asterisks to differentiate members of the same family.  
For example, $\cP$ and $\cP^*$ are both pencils, the latter being minimal in a sense to be made clear 
in section \ref{section:Pencil}.  
If $\cP$ is already minimal, then $\cP$ coincides with $\cP^*$ and the edge simply represents the identity.
%%%%%%%%% FIGURE 
%%%%%%%%%

\begin{figure}[h]
\hfil\epsfig{file=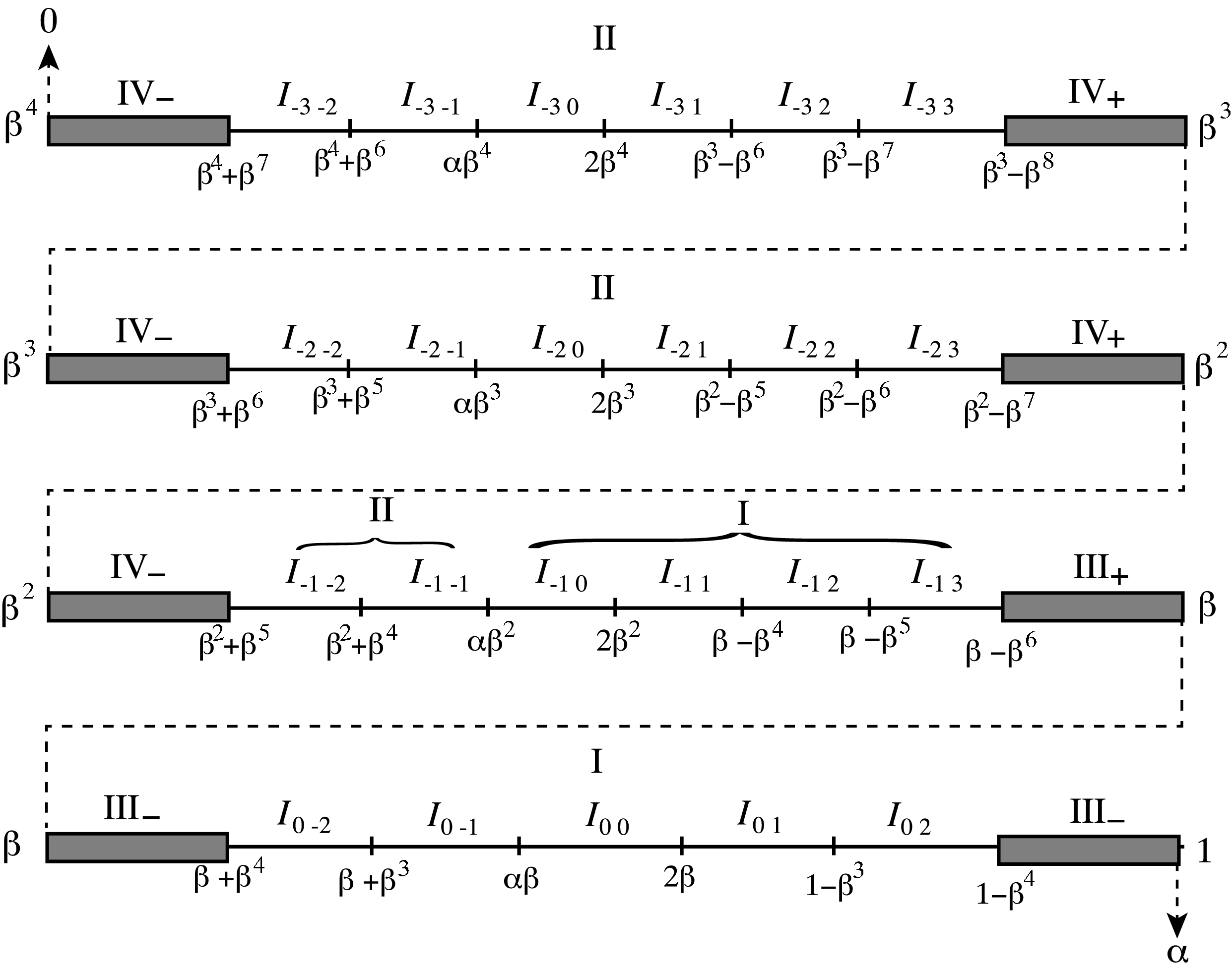,width=5.5 in}\hfil
\caption{\label{fig:sAxis}\small Portion of the $s$-axis, folded, with distorted scale. 
The large-$|j|$ asymptotic regime (scenarios III and IV) is represented by thick segments 
surrounding the accumulation points $\beta^k,\; k=1,2,3,4$. 
The large-$|i|$ regime (scenarios II and IV) consists of the entire interval $[0,\alpha\beta^2]$.}
\end{figure}

\begin{figure}[h]
\hfil\epsfig{file=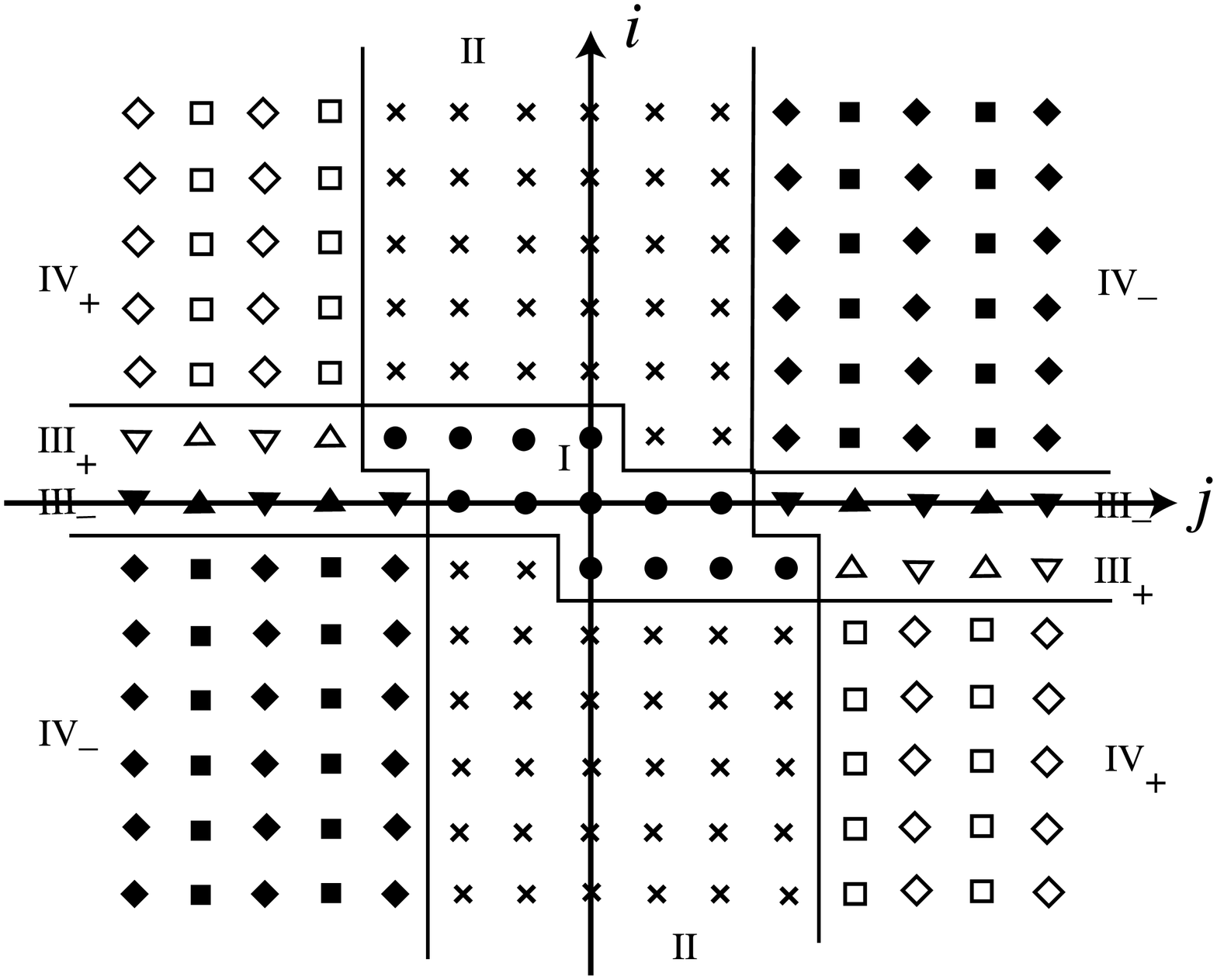,width=5.5 in}\hfil
\caption{\label{fig:4Scenarios}\small The ten renormalization scenarios for the indices 
$(i,j)\in \Z^2$. We write $\pm$ for $\pm1$.
The labels ${\rm III}_\pm$ and ${\rm IV}_\pm$ stand for ${\rm III}_{\pm -}\cup{\rm III}_{\pm +}$ and ${\rm IV}_{\pm -}\cup {\rm IV}_{\pm +}$, respectively.
Key: ${\rm I}=\bullet,\;{\rm II}= \times, {\rm III}_{++}=\vartriangle,\;{\rm III}_{+-}=\triangledown,\;{\rm III}_{-+}=\blacktriangle,\;{\rm III}_{--}=\blacktriangledown,
{\rm IV}_{++}=\square,\; {\rm IV}_{+-}=\diamondsuit,\;{\rm IV}_{-+}=\blacksquare,\;{\rm IV}_{--}=\blacklozenge.$
 }
\end{figure}

Figure \ref{fig:sAxis} emphasizes the organisation of the ten scenarios on the $s$-axis. 
Specifically, we display the sub-intervals $I_{i,j}$ and their assignment to the 
scenarios I--IV for $\beta^4\leqslant s \leqslant 1$. 
The same information is illustrated in figure \ref{fig:4Scenarios} on the $(i,j)$-lattice.
Here the labels I--IV denote subsets of $\mathbb{Z}^2$, with suitable added points at infinity.
We see that scenario I is restricted to pairs $(i,j)$ where both indices are small, plus 
four points at infinity corresponding to $s=\alpha,0$ and $s=\beta,1$, respectively. 
Scenario II corresponds to small values of $j$, with unbounded $i$
(corresponding to $s$ approaching $0,\alpha$), plus a sequence of points at infinity
corresponding to the accumulation points $\beta^{|i|+1}$ and $\alpha-\beta^{|i|+1}$,
for $|i|\geqslant 1$.
Scenario III features small values of $i$ and unbounded $j$
(corresponding to $s$ tending to an accumulation point). 
Finally, scenario IV covers the doubly asymptotic cases.
With reference to table \ref{tbl:scenarios}, we shall use the short-hand notation:
\begin{equation}\label{eq:ShortHandScenarios}
\mathrm{III}=\bigcup_{\mu,\nu} \mathrm{III}_{\mu\nu}
\hskip 30pt
\mathrm{IV}=\bigcup_{\mu,\nu} \mathrm{IV}_{\mu\nu}.
\end{equation}

Our classification scheme leaves open the possibility of more than one realization of each scenario.
Thus, even though mirror sub-intervals $I_{i,j}$ and $I_{-i,-j}$ always belong to the 
same scenario (as evident in the symmetry of Table \ref{tbl:scenarios}), 
the difference in their return paths may be sufficient to have distinct temporal 
scaling properties.   
This essentially doubles the number of incidence matrices which we need to calculate.

The dynamical architecture of renormalization in the present model bears a strong, 
if imperfect, resemblance to that of Rauzy \cite{Rauzy,Veech} to construct renormalizable 
interval exchange transformations (IETs), i.e., maps on an interval $I$ which 
permute $n$ sub-intervals which form a partition of $I$.
At the heart of the Rauzy construction is an irreducible graph (\textit{Rauzy graph}), 
each vertex of which is a `parametric IET' corresponding to a permutation of $n$ 
symbols and parametrized by a positive $n$-vector $\mathbf{s}$ of interval lengths.  
The IET's of the graph are known as a \textit{Rauzy class}. Each vertex has two 
outgoing edges, corresponding to two possible induced return maps, and is associated 
with a matrix transformation in parameter space, $\mathbf{s}\mapsto \mathbf{r}_i\cdot \mathbf{s},\; i=0,1$.  
To search for strictly renormalizable IET's, one considers the closed loops of the Rauzy graph,  
multiplying the $\mathbf{r}_i$ matrices around a loop and seeking a positive eigenvector of the product matrix.

In the present work, our strategy for proving renormalizability is clearly analogous to 
Rauzy's, but of course there are important differences.  In particular, 
our `Rauzy class' of base triangles, pencils, fringed triangles, and double 
strips is a more variegated collection of parametric dressed domains, with 
bifurcating parameter dependence and no uniform rules of induction to compare with Rauzy's.  
Nevertheless, the general strategy (also applied in the context of polygon-exchange 
transformations by Schwartz \cite{Schwartz:14}) is the same.

%%%%%%%%% FIGURE
%%%%%%%%%
%%%%%%%%% FIGURE
%%%%%%%%%
\clearpage
%
%%%%%SECTION
\section{Parametric dressed domains with strips}\label{section:PrototypeDressedDomains}
The reader has already been alerted to the fact that certain classes of parametric dressed domains 
(pencils, fringed triangles, and double strips) play central roles in our renormalization story.  
We now define these objects.  
A common feature of all of them is the presence of a special, quasi-one-dimensional sub-tiling 
which we call the \textit{strip}. 
%%%SUBSECTION
\subsection{The strip}
\begin{figure}[h]
\hfil\epsfig{file=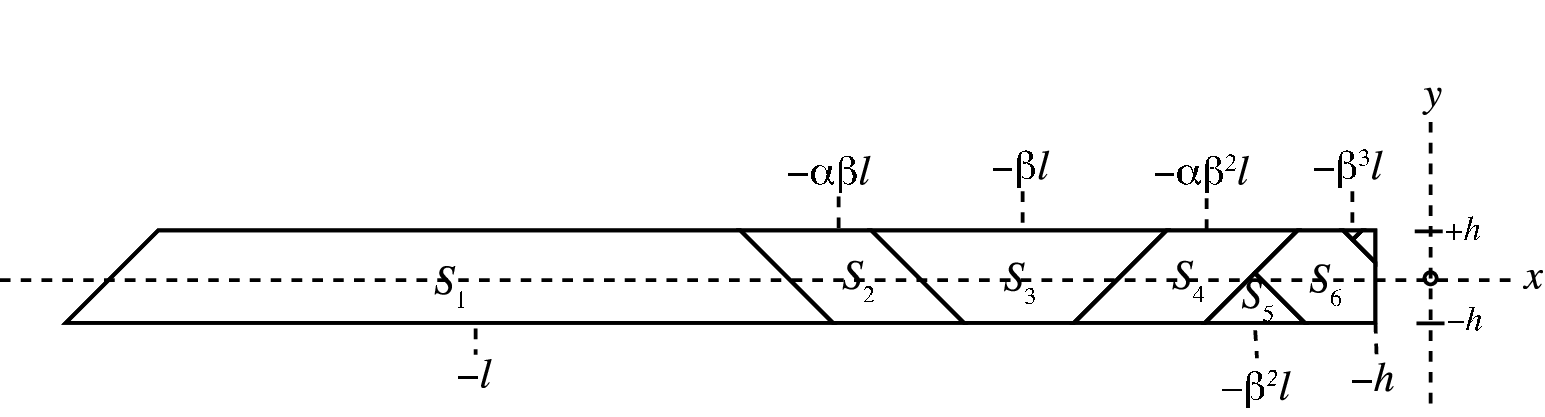,width=5.5 in}\hfil
\caption{\label{fig:trapstrip} \small A strip with 8 tiles.} 
\end{figure}

The prototype strip $\widehat{\bf S}(l,h)$ is a tiling with a variable number
\begin{equation}\label{eq:2J}
2 J = 2 \lfloor \log_\beta(h/l)\rfloor
\end{equation} 
of tiles, all of which are reflection-symmetric. If we let $h$ tend monotonically to zero, the strip undergoes a bifurcation every time $h/l$ assumes a value $\beta^k,\;k=2,3,\ldots$.  The number of tiles increases by two, with the additional tiles being born at one of the vertices at $x= -h$.
The precise structure of  $\widehat{\bf S}(l,h)$ is specified in table 
\ref{tbl:protostrip}, and illustrated, for $J=4$, in figure \ref{fig:trapstrip}.
%%%%%%%%% FIGURE 
%%%%%%%%%
%%%%%%%%%%%TABLE
\begin{table}[h]
\caption{\label{tbl:protostrip}\small Tiling table of the prototype strip $\widehat{\bf{S}}(l,h)$. 
The origin of coordinates lies on the mid-line of the strip, at distance $h$ to the right of its vertical edge.}
$$
\begin{array}{|c|c|c|c|c|}\hline
\multicolumn{3}{|c|}{\mbox{Source Polygon}}&\multicolumn{2}{|c|}{\mbox{Placement}}\\ \hline
\mbox{Tile}&Q_{\#}&\mbox{Parameters}&{\tt R}_{\#}&\mbox{Translation}\\ \hline\hline
\widehat{\rm S}&12&\alpha(l-h),2 h&0&(-\alpha l+\beta h,-h)\\ \hline\hline
\widehat{S}_1& 4&\beta l-\beta h,\beta l -\omega h& 0& (-l,0)\\ \hline
\widehat{S}_2& 3& h& 7 & (-\alpha \beta l,0)\\ \hline
\begin{array}{c}\widehat{S}_{2 j-1}\\ _{2\leqslant j \leqslant J-2}\end{array}& 4 &\beta^j l-\beta h,\beta^j l-\omega h & \begin{array}{cc}4&j \mbox{ even}\\0&j \mbox{ odd}\end{array}&(-\beta^{j-1}l, (-1)^{j-1}(\beta^j l -\alpha h))\\ \hline
\begin{array}{c}\widehat{S}_{2j}\\ _{2\leqslant j \leqslant J-2}\end{array}&3&h& \begin{array}{cc}0&j \mbox{ even} \\7& j\mbox{ odd}\end{array}&(-\alpha\beta^j l,0)\\ \hline
\widehat{S}_{2J-3}&1& \beta^{J-1} l -\beta h&\begin{array}{cc}0&J \mbox{ even}\\4&J \mbox{ odd}\end{array}&(-\beta^{J-2}l, (-1)^J(\beta^{J-1}l -\alpha h))\\ \hline
\widehat{S}_{2J-2}&8&-\alpha\beta^{J-1} l+\alpha \beta h,2 h&\begin{array}{cc}1&J \mbox{ even}\\6&J\mbox{ odd}\end{array}&(-\alpha\beta^{J-1} l-\beta h,(-1)^{J-1} h)\\ \hline
\widehat{S}_{2J-1}&1& \beta^{J} l -\beta h&\begin{array}{cc}4&J \mbox{ even}\\0&J\mbox{ odd}\end{array}&(-\beta^{J-1}l, (-1)^{J-1}(\beta^{J}l -\alpha h))\\ \hline
\widehat{S}_{2J}&7&\alpha\beta^{J-1} l-\alpha h&\begin{array}{cc}3&J \mbox{ even}\\6&J \mbox{ odd}\end{array}&(-h,(-1)^{J-1}(\alpha\beta^{J-1} l-\omega h)))\\ \hline
\end{array}
$$
\end{table}

\clearpage

%%%SUBSECTION
\subsection{The pencil}\label{section:Pencil}
\begin{figure}[t]
\hfil\epsfig{file=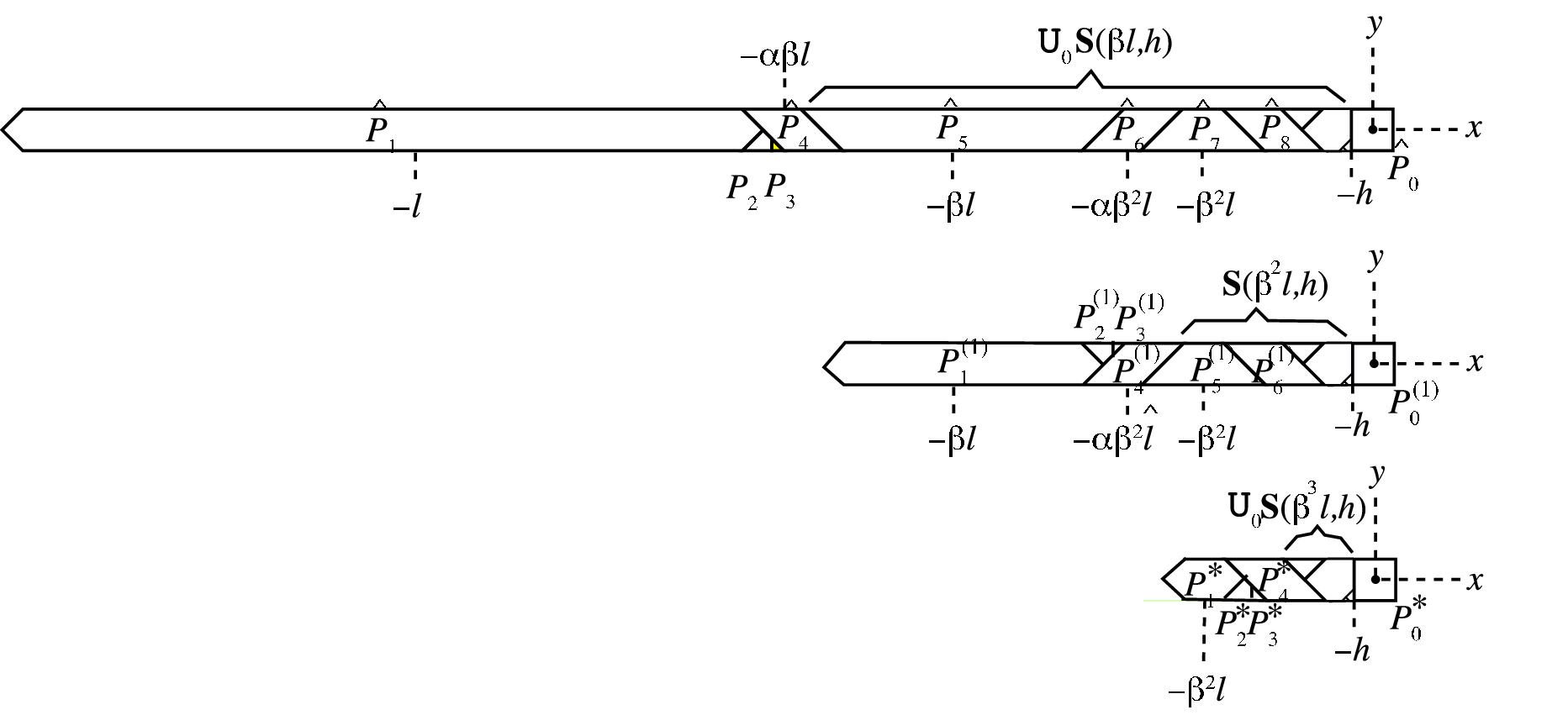,width=5.5 in}\hfil
\caption{\label{fig:A-pencil} \small  The pencil prototype $\widehat{\cP}(l,h)$ with $2L+1$ atoms, 
for $L=6$, together with the two sub-pencils ($L=5$ and $L=4$) obtained by application of the 
Pencil Shortening Lemma.  
Each shortening step (i) contracts the parameter $l$ by a factor $\beta$, 
with no change in the parameter $h$, (ii) reduces the number of atoms by two, and (iii) 
reverses the parity. The nine-atom pencil at the bottom cannot be shortened further 
without violating the definition of a pencil.}
\end{figure}

The prototype pencil $\widehat{\cP}(l,h)$ is a parametric dressed pentagon with a variable number 
$$
2L+1 = 4+\lfloor \log_\beta(\frac{h}{l}) \rfloor
$$ 
of atoms.  Its tiling $\widehat{\mathbf{P}}$ is the union of five individual tiles and a strip, namely
$$
\widehat{\bf P}=P_0 \cup P_1 \cup P_2\cup P_3 \cup P_4
   \cup {\tt U}_0 \widehat{\bf S}(\beta l, h),
$$  
Here we assume a coordinate system aligned with the axis of the pencil and with origin 
at the centre of the square tile $\widehat{P}_0$. 
Explicitly, arranging the tiles, apart from $\widehat{P}_0$, in left-to-right order of their 
anchor points, we have
\begin{eqnarray*}
\widehat{P}_0&=&Q_2(h)\\
\widehat{P}_1&=&{\tt T}_{(-l,0)}\, Q_9(\beta l/\alpha-h,h),\\
\widehat{P}_2&=&{\tt T}_{(-\alpha\beta l-\omega h,-h)}\, {\tt R}_1 \,Q_7(\alpha h),\\
\widehat{P}_3&=&{\tt T}_{(-\alpha\beta l- h,-h)}\, {\tt R}_3 \,Q_1(\beta h),\\
\widehat{P}_4&=&{\tt T}_{(-\alpha\beta l,0)}\, {\tt R}_7 \,Q_3(h),\\
\widehat{P}_{k+4}&=&{\tt U}_0 \,S_k(\beta l,h), \;k=1,\ldots,2L-4.
\end{eqnarray*}
The span of the tiling is
$$
\widehat{\rm P}={\rm span}(\widehat{\bf P})={\tt T}_{(-\alpha l + \alpha h,0)}\, Q_6(\alpha l-\omega h, h).
$$
The structure of the pencil is illustrated in fig. \ref{fig:A-pencil}.
%%%%%%%%% FIGURE 
%%%%%%%%%

The domain map $\rho_{\stackrel{}{\widehat{\cP}}}$ for the pencil is defined in terms of a composition of 
involutions, ${\tt U}_0\circ {\tt H}$, where ${\tt H}$ is a simultaneous reflection of each 
tile about an assigned axis. The trapezia, triangles, kites, and the 
hexagon $\widehat{P}_{2L-2}$ have a unique reflection symmetry.
For the rhombi, we choose the short diagonal, as in the case of the strip.
This leaves the square $\widehat{P}_0$ and hexagon $\widehat{P}_1$, for which we assign axes 
parallel to $\mathbf{u}_7$ and $\mathbf{u}_2$, respectively. 
\medskip 

Studying the renormalization of pencils with arbitrarily many atoms is made manageable 
by the Pencil Shortening Lemma, which relates any pencil to a minimal one with only nine atoms.  
Note that in specifying the associated incidence matrix, we use as matrix indices the 
canonical atom labels shown in figure \ref{fig:trapstrip}.  
For the sake of transparency, we will adopt this convention for all of our incidence 
matrices throughout this article.

\noindent {\bf Pencil Shortening Lemma.}  
{\sl Let $\cP(l,h)=(\mathrm{P}(l,h),{\bf P}(l,h),\rho(l,h))$ be a pencil 
congruent to $\widehat{\cP}(l,h)$, with $2 L+1$ tiles (cf.~equation (\ref{eq:DressedDomain})).
The return maps induced by $\rho(l,h)$ on the tiles
${\rm P}(\beta^k l,h),\; k=1,\ldots,L-4,$ promote the latter
to pencils of parity $(-1)^k$ congruent to $\widehat{P}(\beta^k l,h)$.
For a given $L$, the minimal pencil induced by this shortening process has nine atoms,
parity $(-1)^{L}$, and an incidence matrix (with respect to $\cP(l,h)$) given, 
for $i=0,\ldots,2L+1$, by (we label the rows and columns of $\mathrm{M}$ by the
canonical tile names $\cP$ and $\cP^*$, respectively)}
\begin{eqnarray}
{\rm M}_{P_i,P^*_0}&=&\delta_{i,0},\nonumber\\
{\rm M}_{P_i,P^*_1}&=&\delta_{i,2L-5},\nonumber\\
{\rm M}_{P_i,P^*_2}&=&\delta_{i,2}+ 2 \sum_{j=4}^{2L-5}\delta_{i,j},\label{eq:shortpencilM}\\
{\rm M}_{P_i,P^*_3}&=&2(\delta_{i,2L-6}+\delta_{i,2L-5}),\nonumber\\
{\rm M}_{P_i,P^*_j}&=&\delta_{i,j+2L-8}, \; 4\leqslant j\leqslant 8.\nonumber
\end{eqnarray}
{\sl The polygon $\mathrm{P}(l,h)$ is tiled, up to a set of zero measure, by the return
orbits of the tiles of the minimal pencil, as well as a finite number of periodic tiles.}
\medskip

\noindent {\sc Proof.}
Without loss of generality we assume $\cP(l,h)=\widehat{\cP}(l,h)$.  
We wish to show that the piecewise isometries $\rho^{(k)}$ of the shortened
pencils $\cP^{(k)}={\tt U}_0^k \cP(\beta^k l,h)$ are induced return maps
of $\rho^{(0)}= \rho(l,h)$. It suffices to prove it for $k=1$, since the 
step can be repeated until the pencil is minimal.
The proof is by direct iteration of  $\rho^{(0)}$ on the tiles of $\cP^{(1)}$.  
Only a small number of tiles have non-trivial return orbits.  
To see this, we refer to figure \ref{fig:A-pencil}.  
All of the tiles in the strip ${\bf S}(\beta^2 l,h)$ are 
mapped the same by $\rho^{(0)}$ and by $\rho^{(1)}$.  
The same is true of $P^{(1)}_0$ and $P^{(1)}_4$, and even  $P^{(1)}_1$.  
The remaining tiles, $P^{(1)}_2$ and $P^{(1)}_3$, have short return orbits which we 
calculate explicitly: we find that they pass through, in order,  $P_5,P_4,P_4,P_5$ and $P_5,P_4,P_2,P_4,P_5$, respectively. 

From the structure of the return orbits, we can write down immediately the incidence 
matrix for the recursive step from $L$ to $L-1$. We label the rows and columns of 
the incidence matrix by the canonical tile names of $\mathrm{P}$:
\begin{eqnarray}
{\rm M}_{P^{(0)}_i,P^{(1)}_0}&=&\delta_{i,0},\nonumber\\
{\rm M}_{P^{(0)}_i,P^{(1)}_1}&=&\delta_{i,5},\nonumber\\
{\rm M}_{P^{(0)}_i,P^{(1)}_2}&=&\delta_{i,2}+ 2( \delta_{i,4}+\delta_{i,5})\label{eq:incidenceLL-1}\\
{\rm M}_{P^{(0)}_i,P^{(1)}_3}&=&2(\delta_{i,4}+\delta_{i,5}),\nonumber\\
{\rm M}_{P^{(0)}_i,P^{(1)}_j}&=&\delta_{i,j},\quad 6 \leqslant i\leqslant 2L,\quad  4\leqslant j\leqslant 2L-2,\nonumber
\end{eqnarray}
where in the first four equations the index $i$ varies over its full range: $0\leqslant i\leqslant 2L$.

For the full shortening process of $L-4$ steps, ending with a minimal pencil, the proof is by mathematical induction on $L$.  
The starting point is the case $L=5$, where the one-step incidence matrix is given by (\ref{eq:incidenceLL-1}), 
which coincides with (\ref{eq:shortpencilM}).  Given formula (\ref{eq:shortpencilM}) for a given $L$, we get the 
incidence matrix for $L+1$ by multiplication on the right by the recursion matrix defined by (\ref{eq:incidenceLL-1}). 
One readily verifies that this reproduces the general formulae with $L$ incremented by one.

To prove the completeness of the tiling, it is again sufficient to restrict ourselves to a single shortening step.  
The periodic cells are readily identified cells: the hexagonal period-1 atom $P_1$, the triangular 
period-3 atom $P_3$, and an octagonal period-1 tile inscribed in the rhombic atom $P_4$.  
We explicitly verify that the total area of all return orbits is equal to that of the original pencil.   
That the minimal pencil has 9 atoms follows from the definition of a pencil, while the parity 
of $(-1)^L$ is a consequence of the fact that each of the $L-4$ shortening steps is 
accompanied by a reflection ${\tt U}_0$.  \hfil $\Box$
%
%
%%%%SUBSECTION
\subsection{Fringed triangle}
There are two prototype fringed triangles $\widehat{\mathcal{T}}_\pm(l,h)$, each containing
a strip congruent to $\widehat{\bf S}(l,h)$ with a variable number $2J$ of atoms.  
The total numbers of atoms are  $2J+2$ for $\widehat{\cT}_-(l,h)$ and $2J+7$ for $\widehat{\cT}_-(l,h)$.  
 %%%%%%%%% FIGURE 
\begin{figure}[h]
\hfil\epsfig{file=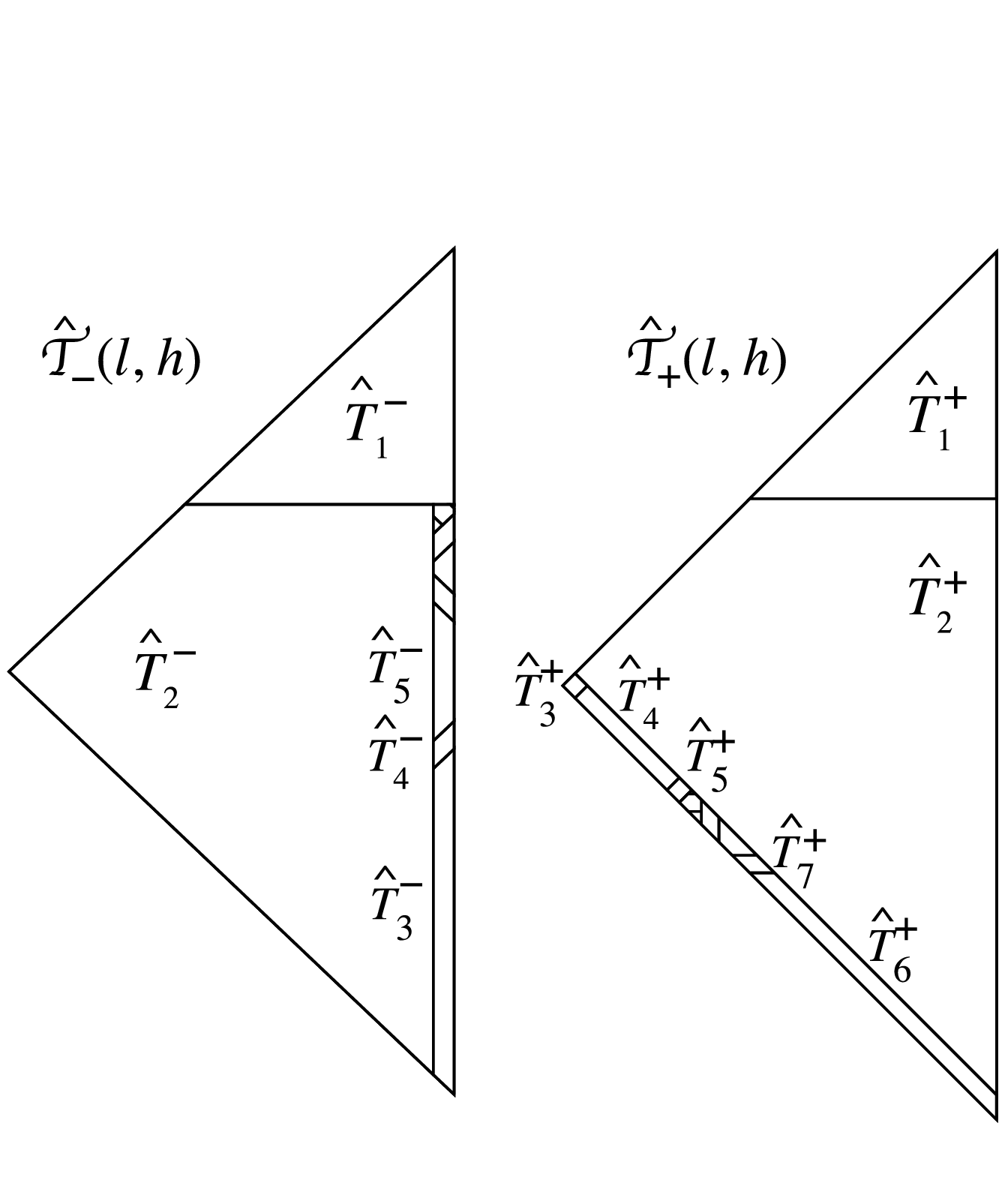,width=4 in}\hfil
\caption{\label{fig:FringedTriangle} \small  Tiling of the prototypes $\widehat{\cT}_{\pm}(l,h)$.}
\end{figure}
%%%%%%%%%

We begin with $\widehat{\cT}_{-}(l,h)$.
Its tiling $\widehat{\bf T}_{-}$ is a union of two individual tiles with a strip, 
namely (see figure \ref{fig:FringedTriangle}) 
$$
\widehat{\bf T}_{-}=\widehat{T}^-_{1}\cup\widehat{T}^-_{2}\cup {\tt R}_2 \,{\bf S}(l, h).
$$  
Here we assume a coordinate system whose origin coincides with that of the strip, 
with the mid-line of the strip lined up along the negative $y$-axis.
Explicitly, 
\begin{eqnarray*}
\widehat{T}^-_{1}&=&{\tt T}_{(h,-h)}\, {\tt R}_5 \,Q_1(\beta l+\beta h),\\
\widehat{T}^-_{2}&=&{\tt T}_{(-h,-\alpha l+\omega h)}\, {\tt R}_3 \,Q_7(\alpha l-\alpha \omega h),\\
\widehat{T}^-_{k+2}&=&{\tt R}_2 \,S_k( l,h), \;k=1,\ldots,2J.
\end{eqnarray*}
The span of the tiling is
$$
\widehat{\rm T}_{-}={\rm span}(\widehat{\bf T}_{-})= {\tt T}_{(-l+\alpha h,-\beta l)}\,{\tt R}_2\,Q_1(l-\beta h).
$$

\bigskip
Next we turn to the prototype fringed triangle $\widehat{\cT}_{+}(l,h)$.
Its tiling is the union of seven individual tiles and a strip, namely
(see figure \ref{fig:FringedTriangle}) 
$$
\widehat{\bf T}_{+}=\bigcup_{k=1}^7\widehat{T}^+_{k}\cup {\tt R}_3 \,\widehat{\bf S}(l, h).
$$  
Here we assume a coordinate system whose origin coincides with that of the strip, with the mid-line 
of the strip lined up along the negative $y$-axis.
Explicitly, 
\begin{eqnarray*}
\widehat{T}^+_{1}&=&{\tt T}_{((3+\alpha) l-h,(1+2\alpha) l+h))}\, {\tt R}_5 \,Q_1(\omega l-\beta h),\\
\widehat{T}^+_{2}&=&{\tt T}_{((3+\alpha) l-h,-(3+\alpha) l+\omega h)}\, {\tt R}_3 \,Q_7(\alpha\omega^2 l-\alpha h),\\
\widehat{T}^+_{3}&=&{\tt T}_{(-\alpha l,\alpha l)}\, {\tt R}_ 7\,Q_2(h),\\
\widehat{T}^+_{4}&=&{\tt T}_{(-l/\alpha,l/\alpha)}\, {\tt R}_7 \,Q_{10}(h,l-h),\\
\widehat{T}^+_{5}&=&{\tt R}_ 7\,Q_2(h),\\
\widehat{T}^+_{6}&=&{\tt T}_{((3+\alpha) l-h,-l-h)}\, {\tt R}_7 \,Q_4(\omega l-\beta h,\omega l-\omega h),\\
\widehat{T}^+_{7}&=&{\tt T}_{(l,-l)}\, {\tt R}_ 7\,Q_3(h),\\
\widehat{T}^+_{k+7}&=&{\tt R}_3 \,S_k(l,h), \;k=1,\ldots,2J.
\end{eqnarray*}
The span of the tiling is
$$
\widehat{\rm T}_{+}={\rm span}(\widehat{\bf T}_{+})
    ={\tt T}_{(-\alpha l+\alpha h,\alpha l)}\,{\tt R}_2\,Q_1(\omega^2 l+\beta h).
$$

The domain maps $\rho_{\stackrel{}{\cT_{\pm}}}$ of the fringed triangles are defined in terms of 
a composition of involutions, namely a simultaneous reflection of each atom about an 
assigned axis, followed by a reflection about the triangle's symmetry axis.  
As in the case of the pencil, the assigned axis of each rhombic atom is its short diagonal.  
For the atoms $\widehat{T}^+_{3}$, $\widehat{T}^+_{4}$, and $\widehat{T}^+_{5}$, the assigned axes 
are parallel to $\mathbf{u}_2, \mathbf{u}_3$, and $\mathbf{u}_2$, respectively.

\clearpage
%
%%%%SUBSECTION
\subsection{Double strip}\label{section:DoubleStrip}
\begin{figure}[h]
\hfil\epsfig{file=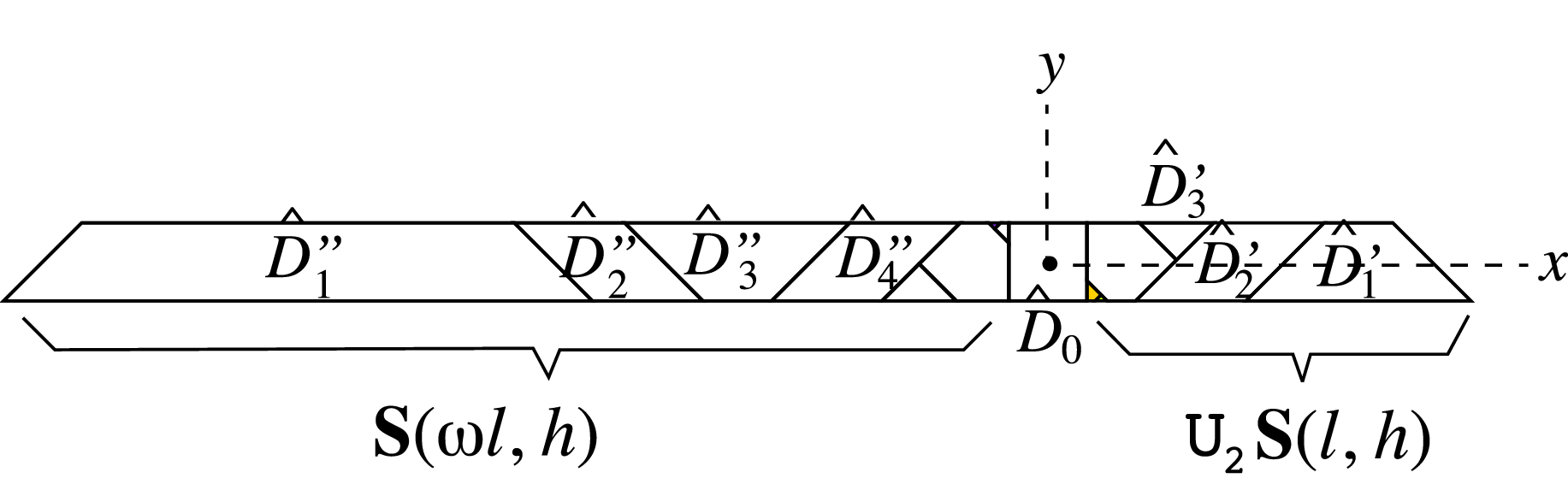,width=4 in}\hfil
\caption{\label{fig:A-dblstrip} \small  Prototype double strip.}
\end{figure}

The prototype \textit{double strip} $\widehat{\cD}_\nu(l,h),\;\nu=\pm 1,$ is a dressed domain 
constructed out of a square and two strips, one on the left with positive parity and $2K$ atoms, 
the other on the right with negative parity and $2K-2$ atoms.  
Since a well-defined strip has at least four tiles, a double strip requires at least 11 atoms 
(i.e., $K\geqslant 3$).
For both signs $\nu$, we define the prototype to have the tiling
$$
 \widehat{\bf D}(l,h)=  {\bf S}(\omega l,h)\, \cup \, \widehat{D}_0 \,\cup {\tt U}_2\,{\bf S}(l,h), 
  \hskip 30pt 
 \widehat{D}_0=Q_2(h).
$$
Note the appearance of the reflection operator ${\tt U}_2$ to correctly place and orient one of 
the component strips.  
A prototype double strip is illustrated in figure \ref{fig:A-dblstrip}. 
Here and in what follows we adopt a canonical labelling of the tiles of any double strip $\cD$, in order along the midline, 
$$
D^{''}_1, D^{''}_2, \ldots, D^{''}_{2K}, D_0, D^{'}_{2K-2},\ldots,D^{'}_1,
$$
 %%%%%%%%% FIGURE 
%%%%%%%%%

The distinction between $\widehat{\cD}_+(l,h)$ and $\widehat{\cD}_-(l,h)$ enters when we specify the piecewise 
isometry $\rho_{\stackrel{}{\widehat{\cD}_\nu}}$.  
As before, we define the map as a composition of two involutions, the reflection 
of each atom about an assigned symmetry axis, followed by reflection about the vertical 
symmetry axis of the double strip as a whole.  Once again the axes of the rhombi are 
their short diagonals.   The square $ \widehat{D}_0$, on the other hand, is assigned 
the diagonal $\mathbf{u}_1$ for $\nu=+$ and $\mathbf{u}_3$ for $\nu=-$.
A considerable simplification of the renormalization structure results from the following 
`shortening' lemma:

\medskip

\noindent {\bf Double-Strip Shortening Lemma.}  
{\sl Let $\cD(l,h)=(\mathrm{D}(l,h),{\bf D}(l,h),\rho(l,h))$ be a double strip 
congruent to $\widehat{\cD}_{\nu}(l,h)$, with $4 K-1$ tiles. 
The return maps induced by $\rho(l,h)$ on the tiles ${\tt U}_2^k{\rm D}(\beta^k l,h),\; k=1,\ldots,K-3,$ 
promote the latter to double strips of parity $(-1)^k$ congruent to $\widehat{D}_{\nu^{(k)}}(\beta^k l,h)$, 
with $\nu^{(k)}=(-1)^k \nu$. For a given $K$, the minimal double strip induced by this shortening 
process has 11 atoms, parity $(-1)^{K-3}$, index $\nu_{K-3}=(-1)^{K-3}\nu$, and an incidence 
matrix (with respect to $\cD(l,h)$) given, for $K$ odd, by}
\begin{eqnarray}
{\rm M}^{\rm odd}_{D, D_0}&=&\delta_{D, D_0} + {\rm N}_D^{\rm odd},\nonumber\\ 
{\rm M}^{\rm odd}_{D, D^{''}_k}&=&\delta_{D, D^{''}_{2K+k-6}} + {\rm N}_D^{\rm odd},\quad 
k=1,\ldots,6,\label{eq:shortModd}\\ 
{\rm M}^{\rm odd}_{D, D^{'}_k}&=&\delta_{D, D^{'}_{2K+k-6}} + {\rm N}_D^{\rm odd},\quad
k=1\ldots,4,\nonumber
\end{eqnarray}
\beq\label{eq:shortNodd}
{\rm N}_D^{\rm odd}= \sum_{k=1}^{\frac{K-3}{2}}\left(2 \cdot 4^{k-1}\,\delta_{D, D^{''}_{2K-4k-5}} 
    + 4^{k-1}\, \delta_{D, D^{'}_{2K-4k-5}}\right),
\eeq
{\sl and, for $K$ even, by}
\begin{eqnarray}
{\rm M}^{\rm even}_{D, D_0}&=&\delta_{D, D_0} + {\rm N}_D^{\rm even},\nonumber\\ 
{\rm M}^{\rm even}_{D, D^{''}_k}&=&\delta_{D, D^{'}_{2K+k-8}} 
     + {\rm N}_D^{\rm even},\quad k=1,\ldots,6,\label{eq:shortMeven}\\ 
{\rm M}^{\rm even}_{D, D^{'}_k}&=&\delta_{D, D^{''}_{2K+k-4}} 
     + {\rm N}_D^{\rm even},\quad k=1,\ldots,4,\nonumber
\end{eqnarray}
\beq\label{eq:shortNeven}
{\rm N}_D^{\rm even}=  \sum_{k=1}^{\frac{K-2}{2}}\left(2 \cdot 4^{k-1}\,\delta_{D, D^{'}_{2K-4k-7}} + 4^{k-1}\, \delta_{D, D^{''}_{2K-4k-3}}\right),
\eeq
{\sl where the subscript $D$ denotes an arbitrary element of $\mathbf{D}(l,h)$.}

{\sl The polygon ${\rm D}(l,h)$ is tiled, up to a set of zero measure, by the return orbits of the 
tiles of the minimal double strip, as well as a finite number of periodic cells.}

\medskip

\noindent {\sc Proof.}
To show that the piecewise isometries $\rho^{(k)}$ of the shortened double strips 
$\cD^{(k)}={\tt U}_2^k \cD_{\nu_k}(\beta^k l,h)$ are induced return maps of $\rho^{(0)}= \rho(l,h)$, 
it suffices to prove it for $k=1$, since this step can be repeated until the double strip is minimal.  
Here we utilize the decomposition of each $\rho^{(k)}$ into a product of  involutions, 
$\rho^{(k)}= {\tt U}^{(k)} {\tt H}^{(k)}, k=0,1,$ where $ {\tt H}^{(k)}$ reflects each 
tile about its own specified symmetry axis, while ${\tt U}^{(k)} $ is a global reflection 
about the symmetry axis of  $\cD^{(k)}$ as a whole. 

A key observation is that the tiles of $\cD^{(1)}$ coincide with the $4K-5$ rightmost tiles 
of $\cD^{(0)}$, and the span of these tiles, ${\rm D}^{(1)}$, is mapped by a single 
application of $\rho^{(0)}$ onto $D^{(0)'' }_1$, the leftmost (and largest) tile of $\cD^{(0)}$.  
Under the global reflection ${\tt U}^{(0)}$, these two trapezoids are reflected about their 
respective symmetry axes and interchanged.  
One important consequence is the identity (for points of ${\rm D}^{(1)}$),
\beq\label{eq:U0H0U0}
{\tt U}^{(1)}= {\tt U}^{(0)}{\tt H}^{(0)}{\tt U}^{(0)}.
\eeq

Two iterations of  $\rho^{(0)}= {\tt U}^{(0)} {\tt H}^{(0)}$ map a point of ${\rm D}^{(1)}$ 
back into that polygon for the first time, hence constitute the first-return map induced by $\rho^{(0)}$. 
We must still show that $\rho^{(0)2}=\rho^{(1)}$ on ${\rm D}^{(1)}$. But this follows from
$$
\rho^{(0)2}={\tt U}^{(0)}{\tt H}^{(0)}{\tt U}^{(0)}{\tt H}^{(0)} = {\tt U}^{(1)}{\tt H}^{(1)},
$$
where we have used (\ref{eq:U0H0U0}) and the fact that ${\tt H}^{(0)}$ and ${\tt H}^{(1)}$ coincide 
on ${\rm D}^{(1)}$. The opposite signs of $\nu^{(0)}$ and $\nu^{(1)}$ are crucial here to maintain 
a consistent symmetry axis for the square tile. That the parity of the double strip changes 
with each shortening step is an obvious concomitant of the 
action of the reflection operator ${\tt U}_2$. 

To see the completeness of the tiling, it is again sufficient to restrict ourselves to the single step, 
from $\cD^{(0)}$ to  $\cD^{(1)}$.  We can focus on those tiles of the former which are not covered by 
the return orbits of the tiles of the latter.  These are precisely $D^{(0)}_j,\; j=2,3,4$. 
From the decomposition $\rho^{(0)}= {\tt U}^{(0)} {\tt H}^{(0)}$, it follows that $D^{(0)}_3$, 
a trapezoid whose symmetry axis coincides with the global symmetry axis, is a period-1 cell, while the 
symmetrically deployed rhombi $D^{(0)}_2$ and $D^{(0)}_4$ form a 2-cycle. 
Thus all points of ${\rm D}^{(0)}$ are covered, up to boundary points, by the return orbits of $\cD^{(1)}$ 
and the periodic cells just discussed.

Finally we turn to the incidence matrices. From our discussion of the two-step return orbits, we can 
immediately write down the incidence matrix for the shortening process from a double strip labeled by 
$K$ to the shortened one labeled by $K-1$.  
Here we label the columns of the incidence matrix by the canonical tile names of $\cD^{(K-1)}$, 
while the row index $D$ stands for an arbitrary tile label of $\cD^{(K)}$.
\begin{eqnarray}
{\rm M}_{D,D_0}&=& \delta_{D,D_0} + \delta_{D,D_1''},\nonumber\\
{\rm M}_{D,D_k''}&=& \delta_{D,D_k'} + \delta_{D,D_1''},\quad k=1,\ldots,2K-2\label{eq:incidenceKK-1}\\
{\rm M}_{D,D_k'}&=& \delta_{D,D_{4+k}''}+ \delta_{D,D_1''} ,\quad k=1,\ldots,2K-4.\nonumber
\end{eqnarray}

For the full shortening process of $K-3$ steps, ending with a minimal double strip, the proof is by 
mathematical induction on $K$.  
The starting point is the case $K=4$, where the one-step incidence matrix is given by (\ref{eq:incidenceKK-1}).  
Given formulae (\ref{eq:shortModd}) and (\ref{eq:shortNodd}), or (\ref{eq:shortMeven}) and (\ref{eq:shortNeven}), 
for a given $K$, we get the incidence matrix for $K+1$ by multiplication on the right by the recursion 
matrix defined by (\ref{eq:incidenceKK-1}). One readily verifies that this reproduces the general formulae 
with $K$ incremented by one.
\hfill $\Box$

%_{\stackrel{}{\mathcal{A}}}

%%%%%SECTION
\section{Arrowheads}\label{section:Arrowheads}

In the preceding section we have obtained a detailed description of the dressed domains
participating in the renormalization. Together they account for all of the vertices
of the renormalization graph. We are now left with the task of establishing the edges.   
Two of the latter have already been discussed:  the inductions $\cP\rightarrow\cP^*$
and $\cD_\mu\rightarrow\cD^*_\nu$  are implemented by the pencil and double-strip shortening
lemmas of the preceding section.  
Three of the others, namely $\cB\rightarrow\cB$, $\cP^*\rightarrow\cB$, and
$\cD^*_\nu\rightarrow \cB$, will be established in the next section by direct iteration of
the parent piecewise isometry. 
As we shall see, the remaining links all involve return-map partitions which produce
strips in the child dressed domain, a process which has at its heart the dynamics of a
parametric, \textit{partially} dressed domain, the {\it arrowhead}. 
In the present section we study arrowhead dynamics, establishing an important lemma 
which will be applied numerous times in the proofs of section \ref{section:ProofOfTheorem1}.

What distinguishes the arrowhead from the parametric dressed domains of the previous 
section is that its piecewise isometry is left undefined on one of its three tiles.  
Thus it cannot be viewed as a self-standing dynamical system.  
As a dressed sub-domain, however, the arrowhead is fully functional, with the missing 
isometry supplied, via induction, by the PWI of its parent. 
The flexibility of this arrangement will allow us, in our proof of various renormalization 
scenarios, to bring to bear the strip-building machinery of the arrowhead in a variety 
of different contexts.  
 
%%%%SUBSECTION
\subsection{Prototype}
For $h\in (0,l)$, we define the arrowhead prototype as
$$
\widehat{\cA}(l,h)= (\mathrm{A},{\bf A},\rho_{\stackrel{}{\mathcal{A}}}), \qquad \mathrm{A}=\mathrm{span}(\mathbf{A}),
$$
where, for $0<h/l<\beta, \; \mathbf{A}=(A_1,A_2,A_3)$, with
\begin{eqnarray*}
A_1&=&{\tt T}_{(-\alpha(l-h),0)}\,{\tt R}_2\, Q_1(\omega(l-h)),\\
A_2&=&{\tt T}_{(-\beta(l-h),l-h)} \,{\tt R}_5 \,Q_{13}(\alpha(l-h),2h),\\
A_3&=&{\tt T}_{(-l-\beta h,-l+\omega h)}\, {\tt R}_4\,Q_1(l-\omega h),
\end{eqnarray*}
and, for $\beta\leqslant h/l<1, \; \mathbf{A}=(A_1,A_2)$, with
\begin{eqnarray*}
A_1&=&{\tt T}_{(-\alpha(l-h),0)}\,{\tt R}_2\, Q_1(\omega(l-h)),\\
A_2&=&{\tt T}_{(-\alpha(l-h),0)}{\tt R}_4\, Q_1(l-h).
\end{eqnarray*}

\begin{figure}[h]
\hfil\epsfig{file=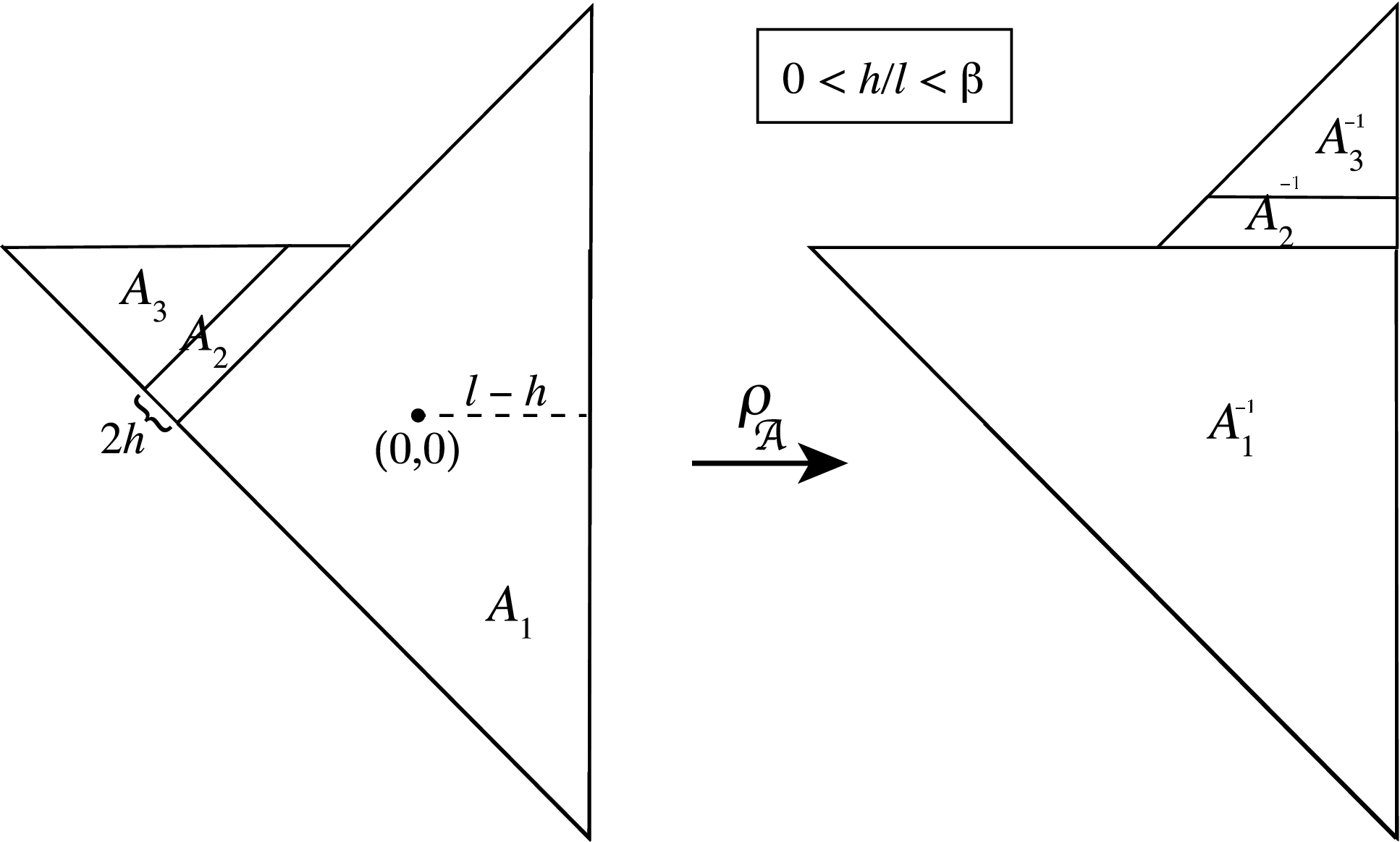,width=4.5 in}\hfil
\caption{\label{fig:A-2} \small The prototype arrowhead $\cA(l,h)$ and its domain map $\rho_{\stackrel{}{\cA}}$, for $0<h/l<\beta$.   }
\end{figure}

\begin{figure}[h]
\hfil\epsfig{file=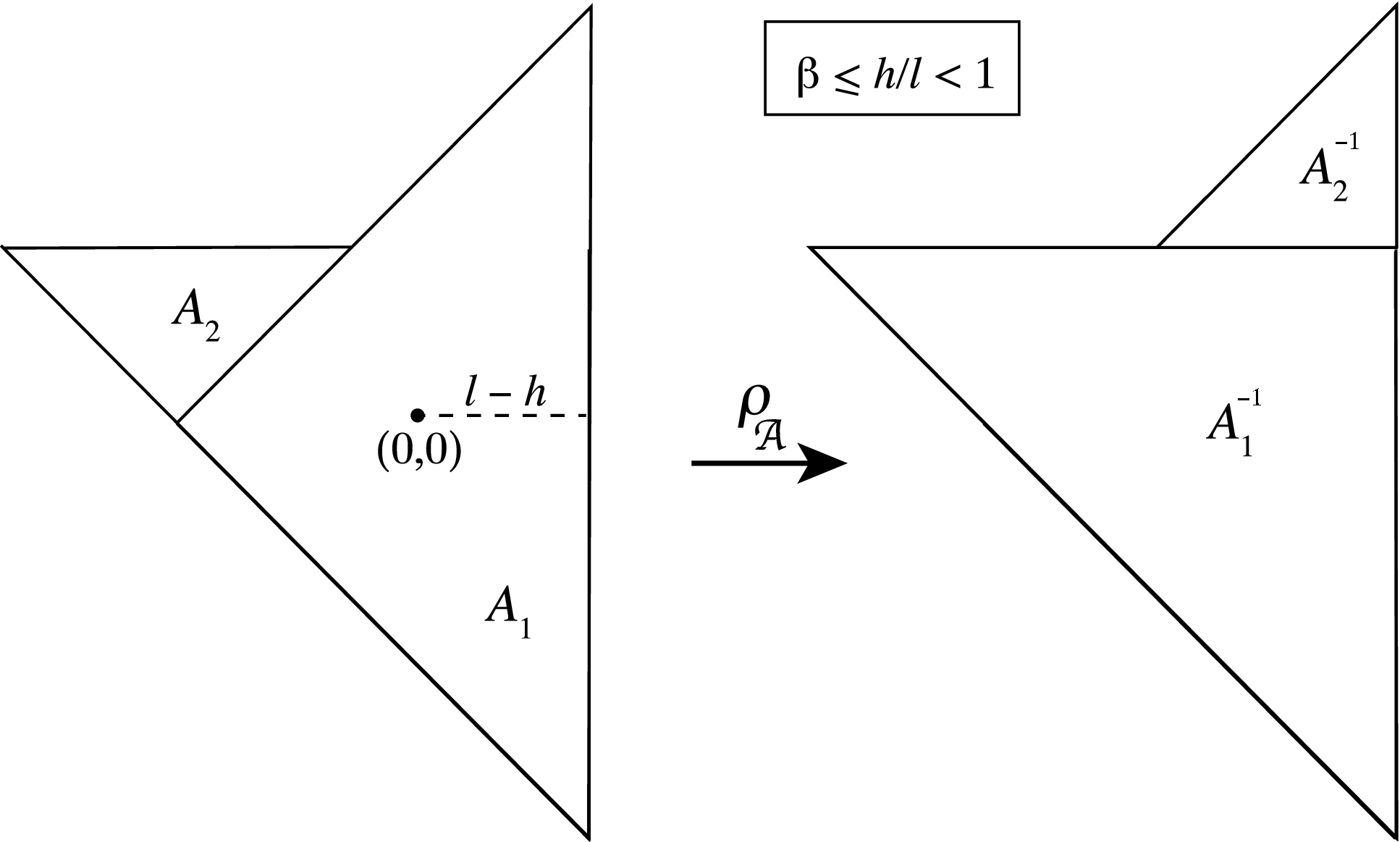,width=4.5 in}\hfil
\caption{\label{fig:A-2b} \small The prototype arrowhead $\cA(l,h)$ and its domain map $\rho_{\stackrel{}{\cA}}$, for $\beta\leqslant h/l<1$. }
\end{figure}

The non-convex polygon $\mathrm{A}$ is equal to the union
of the isosceles right triangle $A_1$ with its reflection about 
the axis $\mathbf{u}_{\frac{5}{2}}$:
$$
{\rm A}= A_1\cup {\tt U}_{\frac{5}{2}}A_1.
$$
Note that the origin of coordinates (anchor point for the arrowhead) has been taken 
to be the in-centre of ${\rm A}$, i.e., the centre of an inscribed circle of radius $l-h$.  
The piecewise isometry $\rho_{\stackrel{}{\mathcal{A}}}$ acts on the tiles of ${\bf A}$ as
$$
\rho_1\stackrel{\rm def}{=}\rho|_{A_1}= {\tt R}_5,\qquad
\rho_3\stackrel{\rm def}{=}\rho|_{A_3}={\tt T}_{(2l,2l)}\, {\tt R}_1,
$$
with the isometry $\rho_2$ on $A_2$ left to be defined by induction in cases where the 
arrowhead is a dressed sub-domain.
Since in all of our applications, the induced map takes $A_2$ outside the arrowhead, we
shall refer to the latter as the {\it exit tile}. The inverse map is given by
$$
\rho^{-1}_{\stackrel{}{\mathcal{A}}}={\tt U}_{\frac{5}{2}}\circ \rho_{\stackrel{}{\mathcal{A}}}\circ {\tt U}_{\frac{5}{2}}.
$$
The atoms $\{A_1^{-1},A_2^{-1},A_3^{-1}\}$ of $\rho^{-1}_{\stackrel{}{\mathcal{A}}}$ are just the
reflected images of those of $\rho_{\stackrel{}{\mathcal{A}}}$, with 
$\rho^{-1}_{\stackrel{}{\mathcal{A}}}$ undefined intrinsically on {\it entrance tile} $A^{-1}_2$.
The parametrization of the arrowhead $\cA(l,h)$ and the action of $\rho_{\stackrel{}{\mathcal{A}}}$ 
is illustrated in figures \ref{fig:A-2} and \ref{fig:A-2b}.
%%%%%%%%% FIGURE 
%%%%%%%%%
%%%%%%%%% FIGURE 
%%%%%%%%%

\subsection{Transfer map and the Arrowhead Lemma}
\begin{figure}[h]
\hfil\epsfig{file=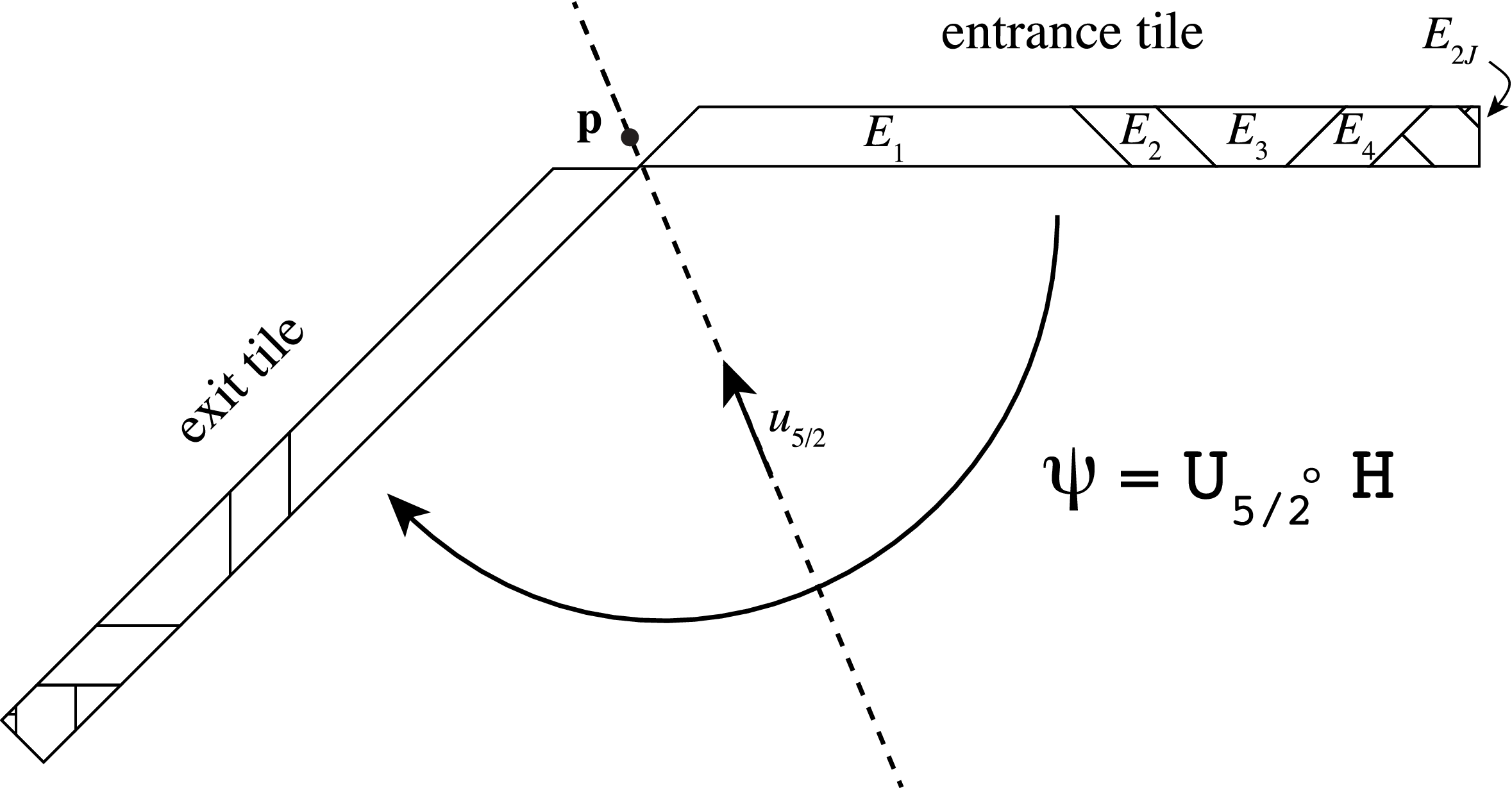,width=4.5 in}\hfil
\caption{\label{fig:AHpsi}\small Illustration of statements $i)$ and $ii)$ of the
Arrowhead Lemma, in an example where the transfer map partition of the entrance
tile is a strip with 8 tiles.  
The point $\mathbf{p}$, at the intersection of the midlines of the entrance and 
exit tiles and the arrowhead symmetry axis, plays an important role in the proof 
of the Arrowhead Lemma.}
\end{figure}

\begin{figure}[h]
\hfil\epsfig{file=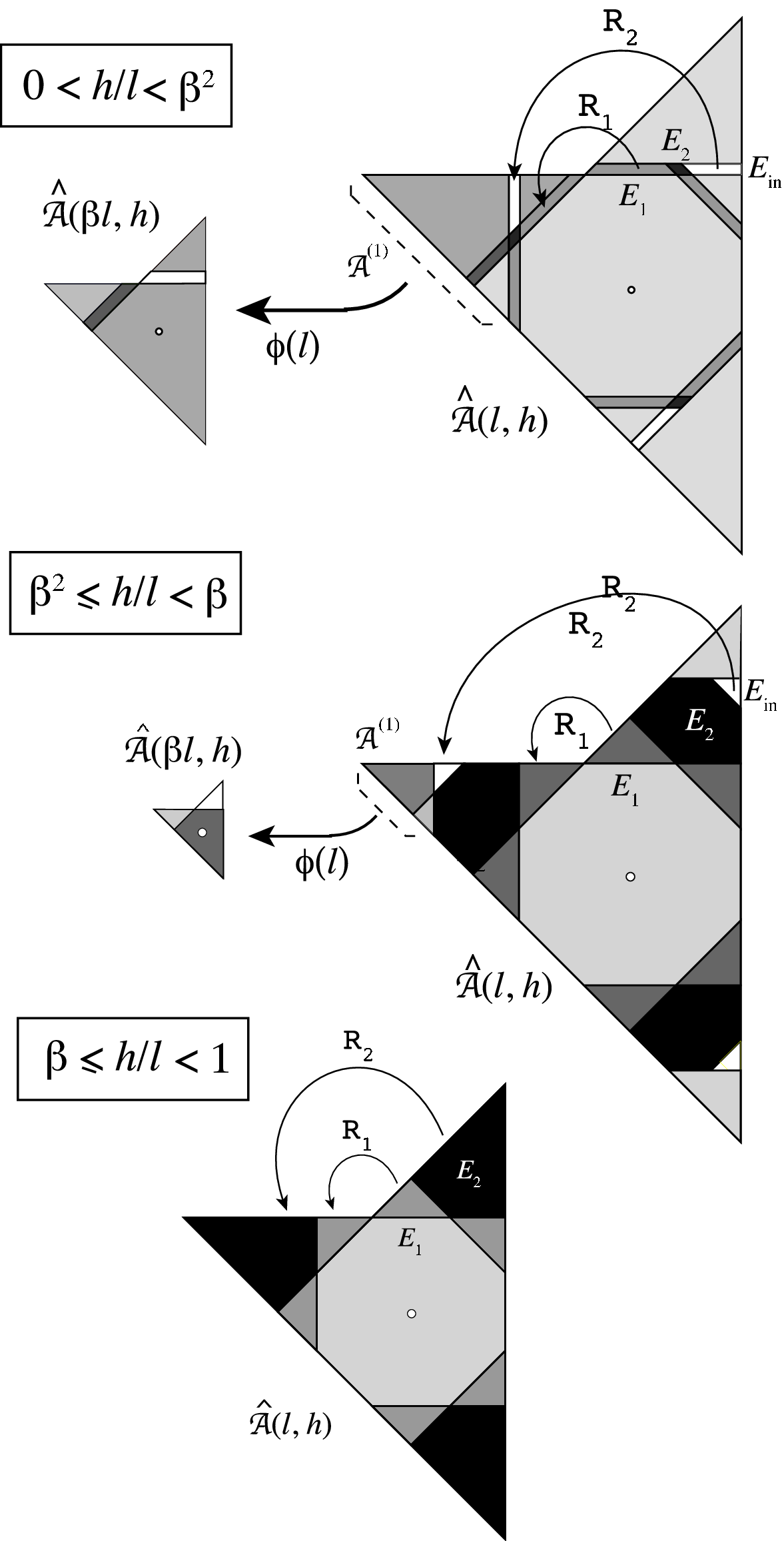,width=4 in}\hfil
\caption{\label{fig:A-3} \small Single-step transfer paths for the three parameter 
ranges of the Auxiliary Lemma.}
\end{figure}

For the  arrowhead $\widehat{\cA}(l,h)$, we define the
{\it transfer map} $\psi$ to be the piecewise isometry induced by 
$\rho_{\stackrel{}{\mathcal{A}}}$ mapping the entrance tile $A^{-1}_2$ onto the exit tile $A_2$. 
The Arrowhead Lemma below shows that this map is well-defined as a composition of two involutions. 
In particular, there is a partition of $A^{-1}_2$ into $2 J(l,h)=\lfloor\log_\beta(h/l)\rfloor$ 
tiles, each of which gets mapped isometrically into $A_2$ by iterations of 
$\rho_{\stackrel{}{\mathcal{A}}}$. 
The area-preserving property of the domain map ensures that the transfer orbits 
are finite.
Figure \ref{fig:AHpsi} illustrates the principal features of $\psi$ in a case 
where $J(l,h)>1$.  
In the special case $J(l,h)=1$, which arises for $h/l\in [\beta,1)$, the transfer 
orbits are displayed in figure \ref{fig:A-3}.
 %%%%%%%%% FIGURE 
%%%%%%%%%%%%%%

%
\medskip

%%%%%%%%%%%ARROWHEAD LEMMA
\noindent\textbf{Arrowhead Lemma.}
{\sl
Let $\cA=\widehat{\cA}(l,h)$, with $l>0$ and $h\in(0,l)$. The following holds:

$i)$ For $h/l\in(0,\beta)$, the tiling $\mathbf{E}=\{E_1,\ldots,E_{2J(l,h)}\}$ of the entrance tile of $\cA$ by
the transfer map $\psi$ coincides with the strip $\mathtt{T}_{(l,l)}\widehat{\bf S}(l,h)$.  For $h/l\in [\beta,1)$, the tiling $\mathbf{E}=\{E_1,\ldots,E_2\}$ is given in table \ref{tbl:AHaux}.

$ii)$ The transfer map acts as a composition of two involutions: a simultaneous reflection of
the tiles $E_i$ of the entrance strip about their respective symmetry axes, followed by a
reflection about the symmetry axis of ${\cA}$. For rhombi, the relevant symmetry axis is the short diagonal.

$iii)$ The incidence matrix column ${\tt N}(E_j)$, listing the number of times the transfer 
orbit of the entrance strip atom $E_j$ visits tiles $A_1$ and $A_3$, is given by
\begin{eqnarray}
{\tt N}(E_{2k-1})&=&\left(\begin{array}{c}
-\frac{1}{2}-(-1)^k+\frac{3}{2} 3^k\\
-\frac{1}{2}+(-1)^k+\frac{1}{2} 3^k\end{array}\right),\nonumber\\
\label{eq:NEk}\\
{\tt N}(E_{2k})&=&\left(\begin{array}{c}
-\frac{1}{2}-\frac{1}{4}(-1)^k+\frac{3}{4} 3^k\\
-\frac{1}{2}+\frac{1}{4}(-1)^k+\frac{1}{4} 3^k\end{array}\right)\nonumber
\end{eqnarray}
with $k=1,\ldots,J(h,l)$.

$iv)$ The orbits of $\psi$, including their destination tiles, together with the periodic 
orbits of the octagonal tiles 
$\Pi^{(k)} = \mathtt{T}_{\mathbf{\gamma}^{(k)}}\,Q_5(\beta^k l-h),\; k=0,1,\ldots,$ $\lfloor\log_\omega(l/h)\rfloor$, 
where
$$
\mathbf{\gamma}^{(k)}=\left\{\begin{array}{ll}
(-\alpha l(1-\beta^k),0) & k\;\mbox{even}\\
(-\alpha l, \alpha\beta^k l) & k\;\mbox{odd},
\end{array}\right.
$$
completely tile ${\rm A}$, up to sets of measure zero.  
The respective paths of the periodic orbits are $\sigma^k(1)$, with the substitution $\sigma$ given by, 
$$
\sigma(1)=(3,1^2),\qquad
\sigma(3)=(1^3),
$$
and so their periods are $3^k$.
}

%The choice of symmetry of the tiles $E_i$ follows the conventions of 
%section \ref{section:PrototypeDressedDomains}.
%The proof of this lemma, which will be applied repeatedly in 
%section \ref{section:ProofOfTheorem1},
%is a recursive one, based on the introduction of a sub-arrowhead 
%$\cA^{(1)}\sim \widehat{\cA}(\beta l, h)$, as well as a transfer map 
%from $A_2^{-1} \cup {\rm A}^{(1)}$ to $A_2\cup {\rm A}^{(1)}$. 
%This is accomplished in the following auxiliary lemma.

Our strategy for proving the Arrowhead Lemma is a recursive one, calculating at each step 
the transfer orbits of a single pair of tiles of $\mathbf{E}$ and mapping the rest 
isometrically into the entrance tile of a sub-arrowhead whose first parameter $l$ has been 
contracted by $\beta$, with $h$ unchanged.  
The top panel of figure \ref{fig:A-3} illustrates this single-step transfer map for 
$0<h/l<\beta^2$.  
The reader can follow by eye the orbits of $E_1$ and $E_2$, from their initial positions 
in $A_2^{-1}$ to their final destinations in $A_2$, along paths $(1,1,1,1,1)$ and $(1,1)$ 
respectively.
Meanwhile, the residual part of the entrance tile, $E_{\mathrm{in}}$, is mapped by two 
iterations of $\rho_{\stackrel{}{\cA}}$ into the entrance tile of the sub-arrowhead $\cA^{(1)}$, 
which is congruent to the prototype $\widehat{\cA}(\beta l,h)$ via an orientation-reversing 
isometry $\phi(l)$.

Repeating the process generates additional tiles $E_j$, until we reach the 
penultimate step, where the parameter ratio is in the range $(\beta^2,\beta)$.
The final induced transfer map, with parameter ratio exceeding $\beta$, is 
completely described by the orbits of two tiles, with no residual part
of the entrance strip, and so the recursion terminates.

%%%%%%AUXILIARY LEMMA
\begin{lemma}\label{lemma:AHaux} \textrm{(\textbf{Auxiliary Lemma}).}
Let $\cA=\widehat{\cA}(l,h)$ with $l>0$ and $0<h/l<1$.  Further, let 
$E_1,E_2,E_{\rm in}, A^{(1)}_1,A^{(1)}_2,A^{(1)}_3$ be tiles within 
$\widehat{\rm A}$ specified in the first and second columns of 
Table \ref{tbl:AHaux} for various ranges of $h/l$.  
The domain map $\rho_{\stackrel{}{\cA}}$ induces a joint 
transfer map $\psi'$ from $A_2^{-1} \cup {\rm A}^{(1)}$ to $A_2\cup {\rm A}^{(1)}$, 
for which the listed tiles are atoms, with respective isometries and transfer 
paths listed in the third and fourth columns of the table. 
The orbits of $\psi'$, including the destination tiles in $A_2$, 
together with the periodic octagon $\Pi$ given in the table, 
completely tile ${\rm A}$, up to a set of measure zero.  
The map $\psi'$, restricted to the domain 
${\rm A}^{(1)}={\rm span}(\{A^{(1)}_1,A^{(1)}_2,A^{(1)}_3\})$, 
promotes the latter to the status of an arrowhead, namely
\beq\label{eq:Astar}
\cA^{(1)}=\phi(l)^{-1} \widehat{\cA}(\beta l,h),\quad 
\phi(l)= {\tt U}_1 \mathtt{T}_{(\alpha l,-\alpha\beta l)}.
\eeq 
\end{lemma}
%%%%%%%%% FIGURE 
%%%%%%%%%%%%%%%%
\clearpage

\begin{table}\caption{Data for the transfer map $\psi'$ of Lemma \ref{lemma:AHaux}.} \label{tbl:AHaux}
$$
\begin{array}{|c|c|c|c|} \hline
h/l\;\;\mbox{Range} & \mbox{Tile} & \mbox{Isometry} & \mbox{Path}\\ \hline
{0<h/l < \beta^2} &
{\begin{array}{l}
E_1={\tt T}_{(0,\alpha(l-h))Q_4(\beta l-\beta h,\beta l-\omega h)}\\
E_2={\tt T}_{(\beta l,l)}{\tt R}_7 Q_3(h)\\
E_{\rm in}={\tt T}_{(\beta(l+h),l+h)}Q_{13}(\alpha\beta l-\alpha h,2h)\end{array}}&
{\begin{array}{l}{\tt R}_1\\{\tt R}_2\\{\tt R}_2\end{array}}&
{\begin{array}{l}1^5\\1^2\\1^2\end{array}}\\ \hline
{\beta^2<h/l < \beta} &
{\begin{array}{l}
E_1={\tt T}_{(0,\alpha l-\alpha h)}\,Q_1(\beta(l-h))\\
E_2={\tt T}_{(\beta(l-h),l-h)}\,{\tt R}_1\,Q_8(\alpha\beta(l-h),2 h)\\
E_{\rm in}={\tt T}_{(l-h,2h)}\,{\tt R}_7\,Q_1(\beta l-h)
\end{array}}&
{\begin{array}{l}{\tt R}_1\\{\tt R}_2\\{\tt R}_2\end{array}}&
{\begin{array}{l}1^5\\1^2\\1^2\end{array}}\\ \hline
{\beta\leqslant h/l} &
{\begin{array}{l}
E_1={\tt T}_{(0,\alpha l-\alpha h)}\,Q_1(\beta(l-h))\\
E_2={\tt T}_{(l-h,\omega(l-h))}\,{\tt R}_6\,Q_7(\alpha(l-h))
\end{array}}&
{\begin{array}{l}{\tt R}_1\\{\tt R}_2\end{array}}&
{\begin{array}{l}1^5\\1^2\end{array}}\\ \hline
{0<h/l<1} &{ \begin{array}{l}
A^{(1)}_1={\tt T}_{-\alpha l,\alpha h)} {\tt R}_4Q_1(l-3 h)\\
A^{(1)}_2={\tt T}_{-l-h,\beta l+\beta h)} {\tt R}_54Q_{12}(\alpha\beta l-\alpha h,2h)\\
A^{(1)}_3={\tt T}_{-\alpha l+\alpha h,0)} {\tt R}_2Q_1\beta l-\omega h)\\ 
\Pi^{(0)}=Q_5(l-h)\end{array} }&{ \begin{array}{c}{\tt T}_{(-2l,l)}{\tt R}_3\\ \mbox{---}\\ {\tt R}_7\\{\tt R}_5\end{array}} & {\begin{array}{l}3,1^2 \\ \mbox{---} \\ 1^3 \\1\end{array}}\\ \hline
\end{array}
$$
\end{table}

%%%%%%PROOF OF AUXILIARY LEMMA
\noindent {\sc Proof of the Auxiliary Lemma.}
For each of the listed parameter ranges, the proof is obtained by explicitly applying the domain map 
$\rho_{\stackrel{}{\cA}}$ along the specified paths in column 4, testing for disjointness  at each step.  
Figure \ref{fig:A-3} illustrates the various orbits of $\psi'$, as well as the conjugacy $\phi(l)$.  
Keeping track of the cumulative mapping relative to the initial tile, we verify that the isometry 
listed in column 3 is correct. To check the completeness of the tiling, we verify that the 
total area of all orbit tiles is equal to that of the polygon ${\rm A}$.  
To prove that the induced isometries on $A_1^{(1)}$ and $A_3^{(1)}$ are indeed those of an 
arrowhead of type $\widehat{\cA}(\beta l,h)$, we verify by a straightforward calculation the identities
$$
\phi(l)\rho_1^2\rho_3\phi(l)^{-1}={\tt R}_5,\qquad
 \phi(l)\rho_1^3\phi(l)^{-1}={\tt T}_{(2\beta l,2\beta l)}\,{\tt R}_1.
$$
$\Box$

%%%%%%PROOF OF ARROWHEAD LEMMA 
\noindent {\sc Proof of the Arrowhead Lemma.} 
For $h/l\in [\beta,1)$, statements $i$) -- $iv$) follow from the Auxiliary Lemma (in this parameter range, $\psi$ coincides with $\psi'$).
For $h/l\in [\beta^k,\beta^{k-1}),\; k>1$, we  partition the entrance tile of $\cA$ as a strip 
with $2J(l,h)$ tiles (see (\ref{eq:2J})),
$$
A^{-1}_2= {\rm span}({\bf E}),\qquad {\bf E}={\tt T}_{(l,l)}{\bf S}(l,h).
$$
We need to prove that each of the tiles $E_i$ is an atom of $\psi$, mapped in accordance with 
statement $ii)$ of the lemma.  
For $i=1,2$, the action of $\psi$ coincides with that of $\psi'$ of Lemma \ref{lemma:AHaux}, namely, 
$$
\psi|_{E_1} = {\tt R}_1 ,\qquad\psi|_{E_2} ={\tt R}_2,\qquad  \psi(E_i)\subset A_2, \; i=1,2.
$$  
For $i>2$,  one shows by explicit calculation that $E_i$ is mapped by $\psi'$ onto 
${\tt R}_2 E_i = E^{(1)}_{i-2} = \phi(l)^{-1} E_{i-2} \subset (A^{(1)}_2)^{-1}$.
Since the image tile is in the entrance strip of the arrowhead $\cA^{(1)}$, we can apply 
Lemma \ref{lemma:AHaux} recursively, with the parameter ratio $h/(\beta^k l)$ increasing by a 
factor $\omega=\beta^{-1}$ at each step. For $1<j\leqslant 2J(l,h)$, the recursion terminates after 
$j$ steps, with 
\begin{eqnarray}
\psi|_{E_{2j-1}}&=&\phi(l)^{-1}\cdots \phi(\beta^{j-2} l)^{-1} 
   {\tt R}_1 \phi(\beta^{j-2} l) {\tt R}_2\cdots \phi(l) {\tt R}_2,\label{eq:psiodd}\\
\psi|_{E_{2j}}&=&\phi(l)^{-1}\cdots \phi(\beta^{j-2} l)^{-1} 
   {\tt R}_2 \phi(\beta^{j-2} l) {\tt R}_2\cdots \phi(l) {\tt R}_2.\label{eq:psieven}
\end{eqnarray}
Inserting
$$
\phi(l)^{-1}={\tt T}_{(-\alpha l,\alpha\beta l)} {\tt U}_1,\qquad
\phi(l){\tt R}_2={\tt U}_3 {\tt T}_{(-\alpha\beta l,-\alpha l)}
$$
and simplifying using operations in the group $\mathfrak{G}$ and commutation relations 
(\ref{eq:CommutationRelations}) p.~\pageref{eq:CommutationRelations}, we get
\beq\label{eq:nestedphiinv}
\phi(l)^{-1}\phi(\beta l)^{-1}\cdots\phi(\beta^m l)^{-1}=
\left\{\begin{array}{ll}
{\tt T}_{(-\alpha l (1-\beta^{m+1}),0)} & m\;\mbox{odd}\\ 
{\tt T}_{(-\alpha l,\alpha\beta^{m+1} l)}{\tt U}_1 & m\;\mbox{even}
\end{array}\right.
\eeq
$$
\phi(\beta^m l) {\tt R}_2 \cdots \phi(l) {\tt R}_2=
\left\{\begin{array}{ll}
{\tt T}_{((\beta^{m+1}-1) l,\beta^{m+1}-1) l,))}&m\;\mbox{odd}\\ 
{\tt U}_0 {\tt T}_{(\beta^{m+1}-1) l,(-\beta^{m+1}-1) l)}&m\;\mbox{even.}
\end{array}\right.
$$
Substituting into (\ref{eq:psiodd}) and (\ref{eq:psieven}) and simplifying, we get 
\begin{eqnarray}
\psi|_{E_{2j-1}}&=&{\tt T}_{(\alpha(\beta^{j-1}-1),\alpha(\beta^{j-1}-1))}\, {\tt R}_1\label{eq:psiodd1}\\
\psi|_{E_{2j}}&=&\left\{\begin{array}{ll}
{\tt T}_{((\beta^{j-1}-1)\beta l,(\beta^{j-1}-1) l)} \, {\tt R}_2& j\;\mbox{odd}\\ 
{\tt T}_{(-\omega l +\beta^{j-1} l, -l +\beta^j l)} & j\;\mbox{even.}
\end{array}\right.\label{eq:psieven1}
\end{eqnarray} 

Next we express the right-hand sides of these formulae in terms of products of reflections.
To this end we write ${\tt R}_m^{\mathbf{w}}$ for the rotation through angle $m\pi/4$ about 
the point $\mathbf{w}$, and $\mathtt{U}_n^{\mathbf{w}}$ for the reflection about the line 
through $\mathbf{w}$ parallel to $\mathbf{u}_n$.  
Now we let $\mathbf{p}=(-\beta l,l)$ be the intersection of the symmetry axis of
the arrowhead with the midline of the entrance and exit tile
---see figure \ref{fig:AHpsi}. Further, we let $\mathbf{q}_i$ 
be the intersection of the preferred symmetry axis of $E_i$ (the short diagonal in the case 
of a rhombus) with the mid line of the entry tile, $y=l$. Explicitly,
$$
\mathbf{q}_{2j-1}=((1-\beta^{j-1})l,l),\qquad   \mathbf{q}_{2j} = ((1-\alpha\beta^j)l,l),\qquad j=1,2,\ldots.
$$
Once again making use of the product and commutation relations (\ref{eq:CommutationRelations}), 
we derive the following expressions for the action of $\psi$ on the atoms $E_i$:
\begin{eqnarray}
\psi|_{E_{2j-1}}&=&{\tt R}_5^{\mathbf{p}}\,{\tt R}_4^{\mathbf{q}_{2j-1}}
    ={\tt U}_{5/2}^{\mathbf{p}}{\tt U}_2^{\mathbf{q}_{2j-1}}\label{eq:psiodd2}\\
\psi|_{E_{2j}}&=&\left\{\begin{array}{ll}
{\tt R}_5^{\mathbf{p}}\,{\tt R}_5^{\mathbf{q}_{2j}}
    ={\tt U}_{5/2}^{\mathbf{p}}{\tt U}_{3/2}^{\mathbf{q}_{2j}}&\quad j\;\mbox{odd}\\ 
{\tt R}_5^{\mathbf{p}}\,{\tt R}_3^{\mathbf{q}_{2j}}
    ={\tt U}_{5/2}^{\mathbf{p}}{\tt U}_{5/2}^{\mathbf{q}_{2j}}&\quad j\;\mbox{even.}
\end{array}\right.\label{eq:psieven2}
\end{eqnarray}
Here the third member of each equation has been obtained by applying the identity
\begin{equation}
{\tt U}_m^{(x,y)} {\tt U}_n^{(w,y)}={\tt U}_m^{(x,y)}{\tt U}_0^{(x,y)}{\tt U}_0^{(w,y)}{\tt U}_n^{(w,y)}
={\tt R}_{2m}^{(x,y)}{\tt R}_{8-2n}^{(w,y)}\nonumber.
\end{equation}
Noting that ${\tt U}_{5/2}^{\mathbf{p}}$ is a reflection about the symmetry axis of the arrowhead, 
we see that formulae (\ref{eq:psiodd2}) and (\ref{eq:psieven2}) give us statement $ii)$ of the lemma.

We next turn to $iii)$.  We recall that the transfer orbit of an atom $E_i$ in the entrance 
tile of $\cA$ passes through a succession of nested arrowheads $\cA^{(j)}$, congruent 
to $\widehat{\cA}(\beta^j l,h)$, on its way to the exit tile.  
The transition from level $j$ to level $j+1$ corresponds to two iterations of the 
isometry $\rho^{(j)}_1$.  The path associated with this transition is related to 
that of its predecessor by the substitution $\sigma: 1\mapsto (3,1^2),  3\mapsto (1^3)$.  
Combining all the pieces in accordance with the last column of Table \ref{tbl:scenarios}, we have 
for the full transfer paths,
$$
{\rm path}(E_1)= 1^5,\qquad
{\rm path}(E_2)= 1^2,
$$
\begin{eqnarray*}
{\rm path}(E_{2k-1})&=& 1^2,\sigma(1)^2,\ldots,\sigma^{k-2}(1)^2,\sigma^{k-1}(1)^5,\qquad k=2,3,\ldots
\\
{\rm path}(E_{2k})&=& 1^2,\sigma(1)^2,\ldots,\sigma^{k-2}(1)^2,\sigma^{k-1}(1)^2,\qquad k=2,3,\ldots
\end{eqnarray*}
Denoting by $n_i(\pi)$ the number of times the symbol $i$ appears in the path $\pi$, we have 
$$
\left(\begin{array}{c} n_1(\sigma^j(1)) \\ n_3(\sigma^j(1))\end{array}\right)=
\left(\begin{array}{cc} 2&3\\1&0\end{array}\right)\left( \begin{array}{c}1\\0\end{array}\right)= 
(-1)^j \left( \begin{array}{c}\frac{1}{4}\\-\frac{1}{4}\end{array}\right)+3^j \left( \begin{array}{c}\frac{3}{4}\\ \frac{1}{4}\end{array}\right),
$$
and hence
$$
{\tt N}(E_1)=\left( \begin{array}{c}5\\0\end{array}\right),\qquad {\tt N}(E_2)=\left( \begin{array}{c}2\\0\end{array}\right),
$$
and for $k>1$,
$$
{\tt N}(E_{2k-1})=\sum_{n=0}^{k-2}(-1)^j \left( \begin{array}{c}\frac{1}{2}\\-\frac{1}{2}\end{array}\right)+3^j \left( \begin{array}{c}\frac{3}{2}\\ \frac{1}{2}\end{array}\right)+
(-1)^{k-1} \left( \begin{array}{c}\frac{5}{4}\\-\frac{5}{4}\end{array}\right)+3^j \left( \begin{array}{c}\frac{15}{4}\\ \frac{5}{4}\end{array}\right),
$$
$$
{\tt N}(E_{2k})=\sum_{n=0}^{k-1}(-1)^j \left( \begin{array}{c}\frac{1}{2}\\-\frac{1}{2}\end{array}\right)+3^j \left( \begin{array}{c}\frac{3}{2}\\ \frac{1}{2}\end{array}\right).
$$
Summing up the geometric series, we get the formulae in (3).

Finally, we turn to $iv)$.  We recall once again the nested sequence of arrowheads $\cA^{(k)}$, 
whose successive in-centres are related by the mappings $\phi(\beta^k l)^{-1}$. 
The in-centre of $\cA^{(k)}$ is thus 
$$
 \gamma^{(k)}=\phi(l)^{-1}\phi(\beta l)^{-1}\cdots\phi(\beta^{k-1} l)^{-1}\left(\begin{array}{c} 0 \\ 0\end{array}\right).
$$
The formula in $iv)$ follows from substitution of (\ref{eq:nestedphiinv}).  
The path follows, by recursive application of the substitution $\sigma$ on the lowest-level path,
${\rm path}(\Pi)=(1)$.
This completes the proof of the Arrowhead Lemma.
$\Box$

%
%
%
%%%%%SECTION 10  PROOF OF THEOREM 1
%
\section{Proof of Theorem 1}\label{section:ProofOfTheorem1}
We are now in a position to establish the edges of the renormalization graph, 
thus completing the proof of Theorem 1.  
As a by-product, we will calculate the incidence matrices which together specify the 
temporal scaling behaviour over the entire parameter interval.  
Most of the induction proofs naturally split into two parts, a preliminary part 
in which a fixed collection of return-map orbits are constructed by direct iteration 
of a given piecewise isometry, and a secondary part, containing all the recursive 
branching, which is handled by application of the Arrowhead Lemma or one of the 
Shortening Lemmas.   
%
%

%%%%%SUBSECTION
\subsection{Tiling plans and incidence matrices}\label{section:TilingPlans}
The computer-assisted elements of our proofs consist of direct calculation of finite 
orbits of polygonal domains under the domain map of a given dressed domain.  
In each case, all of the information needed to set up and execute these calculations 
is presented in tabular form as a \textit{tiling plan} for an edge 
$\mathcal{X}\to\mathcal{Y}$ of the renormalization graph.
In Appendix B we display a selected list of tiling plans; a complete record
of the computer-assisted proof is available in the Electronic 
Supplement \cite{ESupplement}.

Each tiling plan is to be validated for either a single value of the parameter $s$, 
or for an interval of values of $s$, using the direct iteration method described
in section \ref{section:DirectIteration}.  
Employing the procedures of our CAP Toolbox (see Electronic Supplement \cite{ESupplement}),
we construct the orbit of each source tile of the tiling plan, checking that it reaches 
its assigned destination without intersecting any of the other destination tiles prior 
to the final step. This guarantees that the orbits are disjoint.
We also check that the isometric mapping between source and destination is as specified in the plan.
As a by-product of the orbit construction, we obtain for each entry various information 
about the orbits, including the number of iterations and the column of the incidence matrix 
giving the number of visits to each of the atoms of the parent dressed domain.  As the final step in 
the proof, we show the completeness of the tiling by verifying that the sum of the areas of the 
tiles of all the orbits is equal to that of the parent domain.

In the present section we will denote by $\mathrm{M}_\lambda(\mathcal{X}\to\mathcal{Y})$ 
the incidence matrix associated with the edge $\mathcal{X}\to\mathcal{Y}$ of the 
renormalization graph, where $\lambda$ stands for one or more of the indices 
$\epsilon,\mu,\nu$ on which the matrix depends.  Here $\epsilon=\mathrm{sign}(i)$ 
and $\mu$ and $\nu$ are functions of $i$ and $j$ given in table \ref{tbl:scenarios}.  
For multi-edge paths, we will add a Roman numeral superscript to identify the 
appropriate scenario and make  the dependence on $i$ and $j$ unique.  
For example, the composite incidence matrix the edge sequence 
$\cB\to\cP\to\cP^*\to\cT_\mu\to\cD_\mu$ will be written as 
$\mathrm{M}^{\mathrm{III}}_{\epsilon,\mu}(\cB\to\cD)$, with the matrix product expansion
$$
\mathrm{M}^{\mathrm{III}}_{\epsilon,\mu}(\cB\to\cD)=\mathrm{M}_\epsilon(\cB\to\cP)\cdot\mathrm{M}(\cP\to\cP^*)\cdot\mathrm{M}_\mu(\cP^*\to\cT)\cdot\mathrm{M}_\mu(\cT\to\cD).
$$

%
%%%%%%%SUBSECTION
\subsection{Proof of $\mathcal{B}\to\mathcal{B^*}$ (scenario I)}\label{section:ProofOfScenarioI}
We begin our proof of theorem 1 by establishing statement $i)$ for $s\in\{0,\alpha\}$, 
statement $ii)$ for $i=\pm 1$, and the following proposition for the remaining $(i,j)$ of 
scenario I (see table \ref{tbl:scenarios} and figures \ref{fig:sAxis} and \ref{fig:4Scenarios}).

\begin{proposition}\label{prop:BB*} 
Let $(i,j)\in \mathrm{I}\cap\Z^2$, let $s \in I_{i,j}$, and let $\cB \sim \widehat{\cB}(1,s)$.  
Then $\cB\rightarrow \cB^*$ where 
\begin{equation}\label{eq:BB*}
\cB^* \sim \widehat{\cB}(l^*,r(s)l^*),\; \mbox{with}\; \; l^*=\beta^{|i|+|j|+2},\;\pi(\cB^*)=(-1)^{|i|+|j|}.
\end{equation}
The incidence matrices for this scenario are given in Appendix C.
\end{proposition}

\begin{figure}[h!]
\hfil\epsfig{file=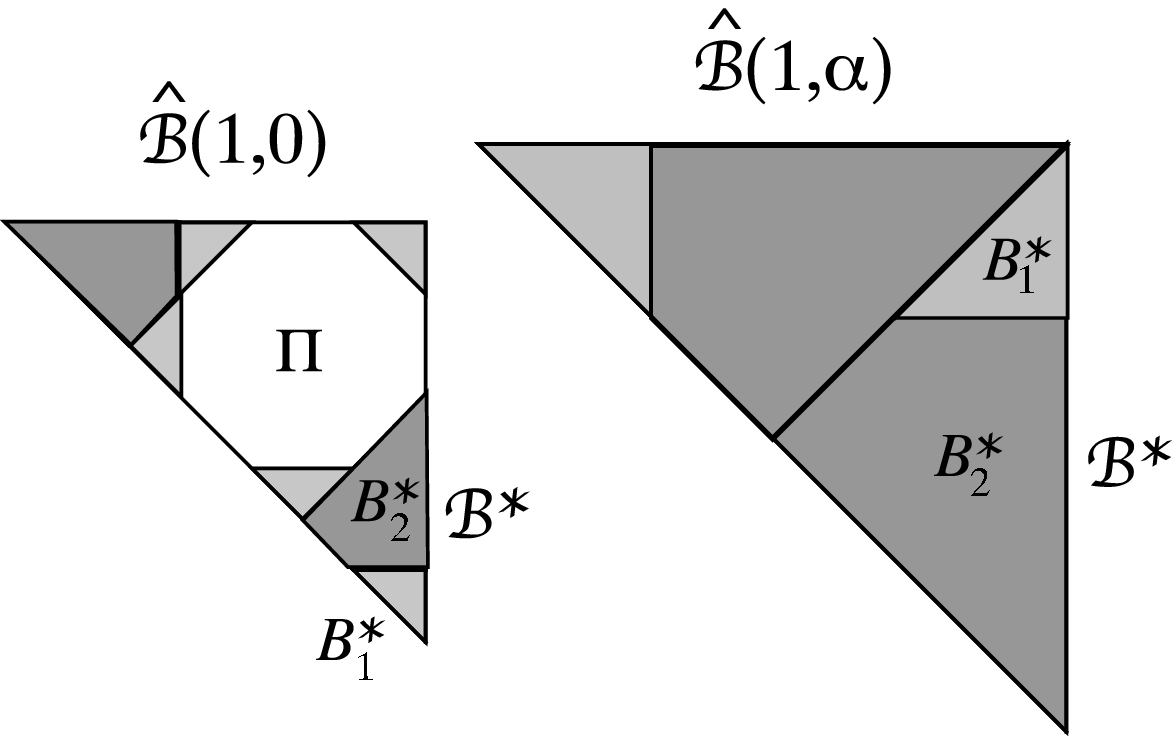,width=4 in}\hfil
\caption{\label{fig:tiling0} \small Return orbits for $\cB\rightarrow\cB^*$, $s=0$ (left) and $s=\alpha$ (right). }
\end{figure}

 \noindent{\sc Proof of theorem 1, statement $i)$.} For $s=0$ we assume,
without loss of generality, that $\cB =\widehat{\cB}(1,0)$.  
The data for this dressed domain and its atoms $B_1$ and $B_2$ are displayed in 
table \ref{tbl:Bdef} p.~\pageref{tbl:Bdef} with $l=1$ and $h=0$.   
By direct iteration of $\rho_{\stackrel{}{{\cB}}}$ on the tiles
\begin{eqnarray*}
 B^*_{1}&=&\mathtt{T}_{(0,-4\beta)}\mathtt{R}_7 Q_1(\alpha\beta^2),\\
 B^*_{2}&=&\mathtt{T}_{(0,-2\beta)}\mathtt{R}_6 Q_1(2 \beta),\\
 \Pi&=&\mathtt{T}_{(-\alpha\beta^2,-\alpha)} Q_5(\alpha\beta^2),
\end{eqnarray*} 
one verifies that the three orbits tile the span of $\cB$ (see figure \ref{fig:tiling0}), 
and  produce a return map which promotes $B_1$ to a positive-parity dressed domain 
$\cB^*$ congruent to $\widehat{\cB}(\beta,0)$. The incidence matrix is included in Appendix C.
%%%%%%%%% FIGURE 
%%%%%%%%% 

Turning now to the other endpoint, we assume $\cB =\widehat{\cB}(1,\alpha)$.  
The data for this dressed domain and its atoms $B_3$, $B_4$, and $B_5$ are displayed in 
table \ref{tbl:Bdef} with $l=1,\; h=\alpha$.   
By direct iteration of $\rho_{\stackrel{}{{\cB}}}$ on the tiles
\begin{eqnarray*}
 B^*_{1}&=&\mathtt{T}_{(0,-2\beta)}\mathtt{R}_5 Q_1(\alpha\beta),\\
 B^*_{2}&=&\mathtt{T}_{(0,-2\alpha)}\mathtt{R}_3 Q_7(2),
\end{eqnarray*} 
one verifies that the two orbits tile the spanning domain of $\cB$ (see figure \ref{fig:tiling0}), and  produce a return map which promotes $B_3$ to a negative-parity dressed domain $\cB^*$ congruent to $\widehat{\cB}(1,0)$.  The incidence matrix is given in Appendix C.

\noindent$\Box$\medskip

\noindent{\sc Proof of theorem 1, statement $ii)$, for $i=\pm 1$.}  
The parameter values $s=\Delta_{-1}=\Delta_0=\beta$ and $s=\Delta_1=1$ are distinguished 
from the other cases of scenario I by the presence of two higher-level base triangles 
with disjoint return orbits, both of which are needed to complete the tiling of the 
parent base triangle.  
The case $s=\beta$ is proved by direct iteration according to Tiling Plan 2, which is 
reported in Appendix B and illustrated in figure \ref{fig:scen1}; the treatment of $s=1$ 
is included in the Electronic Supplement. \/ $\Box$

\medskip
An important consequence of this is the splitting of the exceptional set into disjoint 
ergodic components. Such a behaviour, already observed in quadratic two-dimensional 
piecewise isometries \cite{KouptsovLowensteinVivaldi} here takes a very simple form.
Moreover, infinitely many examples of it appear in our family, corresponding to
the set of all accumulation points, at $\beta^k$ and  $\alpha-\beta^k,\; k=1,2,\ldots$.  
(The cases with $k>1$ belong to scenario II, to be treated later.)

%%%%%%%%% FIGURE 
\begin{sidewaysfigure}[h]
\hfil\epsfig{file=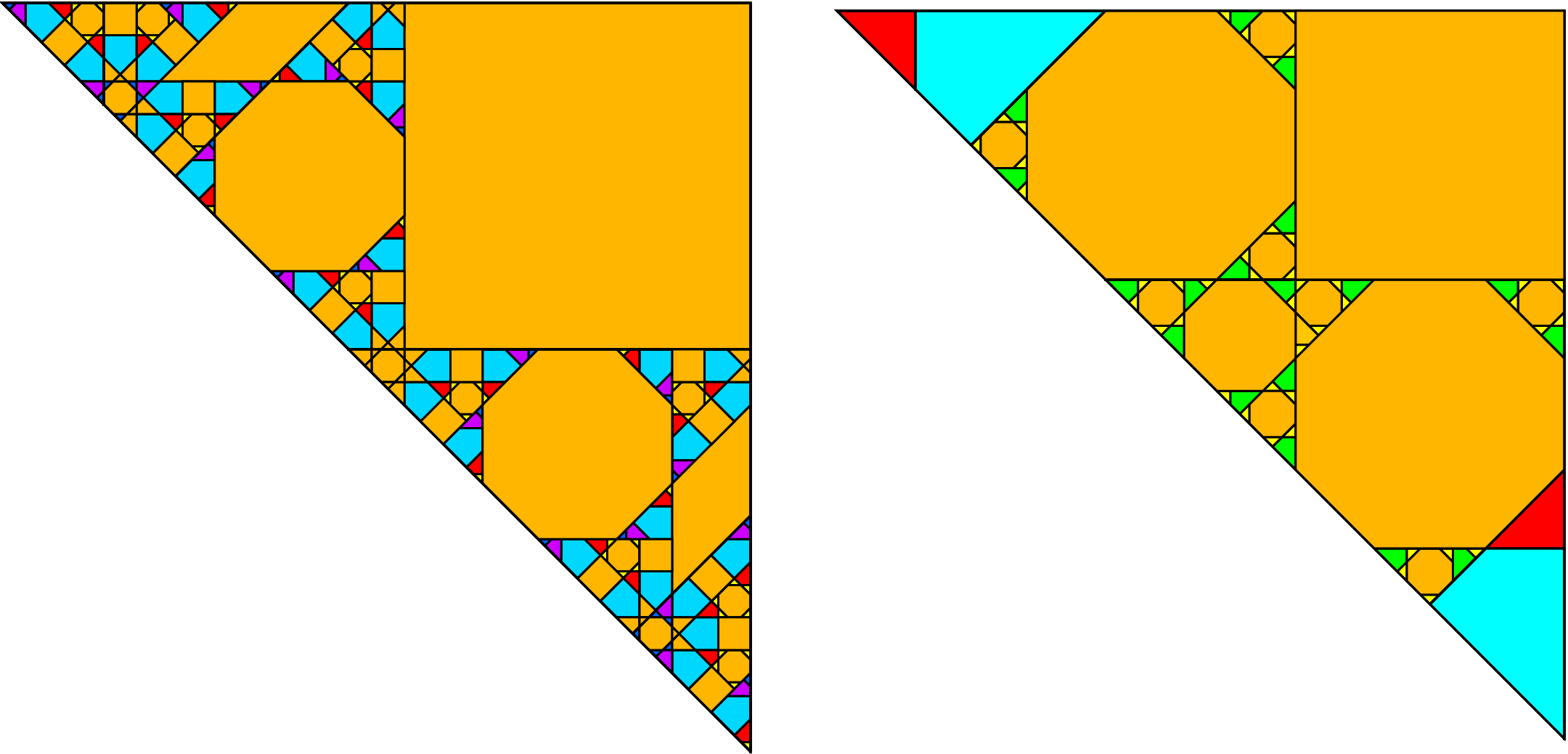,width=8 in}\hfil
\caption{\label{fig:scen1} \small Illustration of the complete tiling of $\mathrm{B}(1,s)$ 
by the return orbits of the higher-level base triangle(s) and a finite number of periodic domains.  
The tiling on the left is for $s=107/200 \in I_{0,-1}$, while that on the right is for the 
accumulation point $s=\beta$.}
\end{sidewaysfigure}

\noindent{\sc Proof of proposition \ref{prop:BB*}}.
For all finite $(i,j)$ of scenario I, we prove the renormalization $\mathcal{B}\to\mathcal{B}^*$ 
by direct iteration of the domain map $\rho_{\stackrel{}{\cB}}$.
Let us illustrate the salient features of the calculations in the case $(i,j)=(0,-1)$, 
corresponding to $s\in I_{0,-1}=(\beta+\beta^3,\alpha\beta]$.
The corresponding Tiling Plan 1, shown in Appendix B, 
has been validated using the procedure specified in section \ref{section:DirectIteration}.  
Extension of these results to the endpoint  $\alpha\beta $ is straightforward,
requiring us to ignore those tiles of $\cB^*$ and periodic domains of Tiling Plan 1
which degenerate to lower-dimensional objects, and allow for the possibility of 
redundant edge conditions in the specification of tile shapes.  
In the present example, only 6 of the original 14 tiles survive at the endpoint, 
with $\cB^*$ degenerating to a two-atom right triangle, as it should in accordance 
with the vanishing of $r(\alpha\beta)$. Even though $\cB^*$ for $s=\alpha\beta$ 
has only two atoms, their return paths are the same as in the open interval, 
and their incidence matrix is formed by the first two columns of the $5\times 5$ 
matrix ${\rm M}(0,-1)$ displayed in Appendix C.
\clearpage

The tiling plans for the remaining $(i,j)$-values of scenario $I$ are analogous. Equivalent data tables will be found in the Electronic Supplement, while the corresponding incidence
matrices are listed in Appendix C. 

%%%%SUBSECTION 10.3

\subsection{Proof of $\mathcal{B}\to\mathcal{P}$ (scenarios II and IV)}\label{section:ProofOfB-P}
Scenarios II and IV deal with the peripheral parts of the $s$-interval, namely 
$s\leqslant \alpha \beta^2$ and  $s\geqslant \alpha-\alpha\beta^2=2\alpha\beta$.
In these regions we find asymptotic phenomena, which develop at the accumulation
points of the singularities of the renormalization function $r$.
The analysis is divided into two cases (see figure \ref{fig:4Scenarios}). 
In Scenario II, the index $|i|$ diverges while $j$ remains within small bounds, so that 
$s$ approaches one of the accumulation points $s=0,\alpha$ of $f$, without 
approaching the first-order accumulation points $s=\beta^k$ or $\alpha-\beta^k$. 
Scenario IV deals with larger values of $|j|$, and includes all doubly asymptotic 
cases: $|i|,|j|\to\infty$. We shall establish the following result
(the set \textrm{IV} is defined in (\ref{eq:ShortHandScenarios}), 
p.~\pageref{eq:ShortHandScenarios}):

\begin{proposition}\label{prop:BPP*}
Let $(i,j)\subset  {\rm II}\cup {\rm IV}$, 
let $s \in I_{i,j}$, and let $\cB \sim \widehat{\cB}(1,s)$.  
Then 
\beq\label{eq:BP}
\cB\rightarrow \cP\sim \widehat{\cP}(\alpha, s),\; \mbox{with}\;\; \pi(\cP)=1.  
\eeq
Moreover,
\beq\label{eq:PP*}
\cP\rightarrow\cP^*\sim\left\{\begin{array}{l}
\widehat{\cP}(\beta^{|i|-1}\alpha,s),\quad \mbox{with}\quad\pi(\cP^*)=(-1)^{|i|-1},
            \quad \mbox{if}\;\; i\cdot j>0,\\ 
\widehat{\cP}(\beta^{|i|-2}\alpha,s),\quad \mbox{with}\quad\pi(\cP^*)=(-1)^{|i|},
            \quad \mbox{if}\quad i\cdot j\leq0.
\end{array}\right.
\eeq
The incidence matrix for the combined renormalization step $\cB\rightarrow \cP^*$ is 
given by (\ref{eq:matrixBP*}).
 \end{proposition}
\bigskip

\noindent {\sc Proof.} 
The proof of $(\ref{eq:BP})$, including the calculation of the incidence matrix, 
is performed separately for the $s$ intervals $(0,\alpha\beta^3)$,  
$[\alpha\beta^3,\alpha\beta^2]$, and  $[2\alpha\beta,\alpha)$.
\medskip

\noindent{\sc Proof of (\ref{eq:BP}) for $s\in (0,\alpha\beta^3]$}.
We prove  (\ref{eq:BP}) in two stages.  The first, non-branching, stage involves the disjoint 
return orbits of a fixed number of initial tiles.  
As $s$ ranges over the interval of interest, the orbits evolve continuously without 
bifurcations. The second stage of the proof deals with all of the $s$-dependent branching, 
which is entirely accounted for by the transfer-map dynamics of an arrowhead dressed 
domain $\cA=(\mathrm{A},\mathbf{A},\rho_{\stackrel{}{\mathcal{A}}})$ congruent
to $\widehat{\cA} (\alpha\beta^2,s)$.

%%%%%%%%% FIGURE 
\begin{figure}[h!]
\hfil\epsfig{file=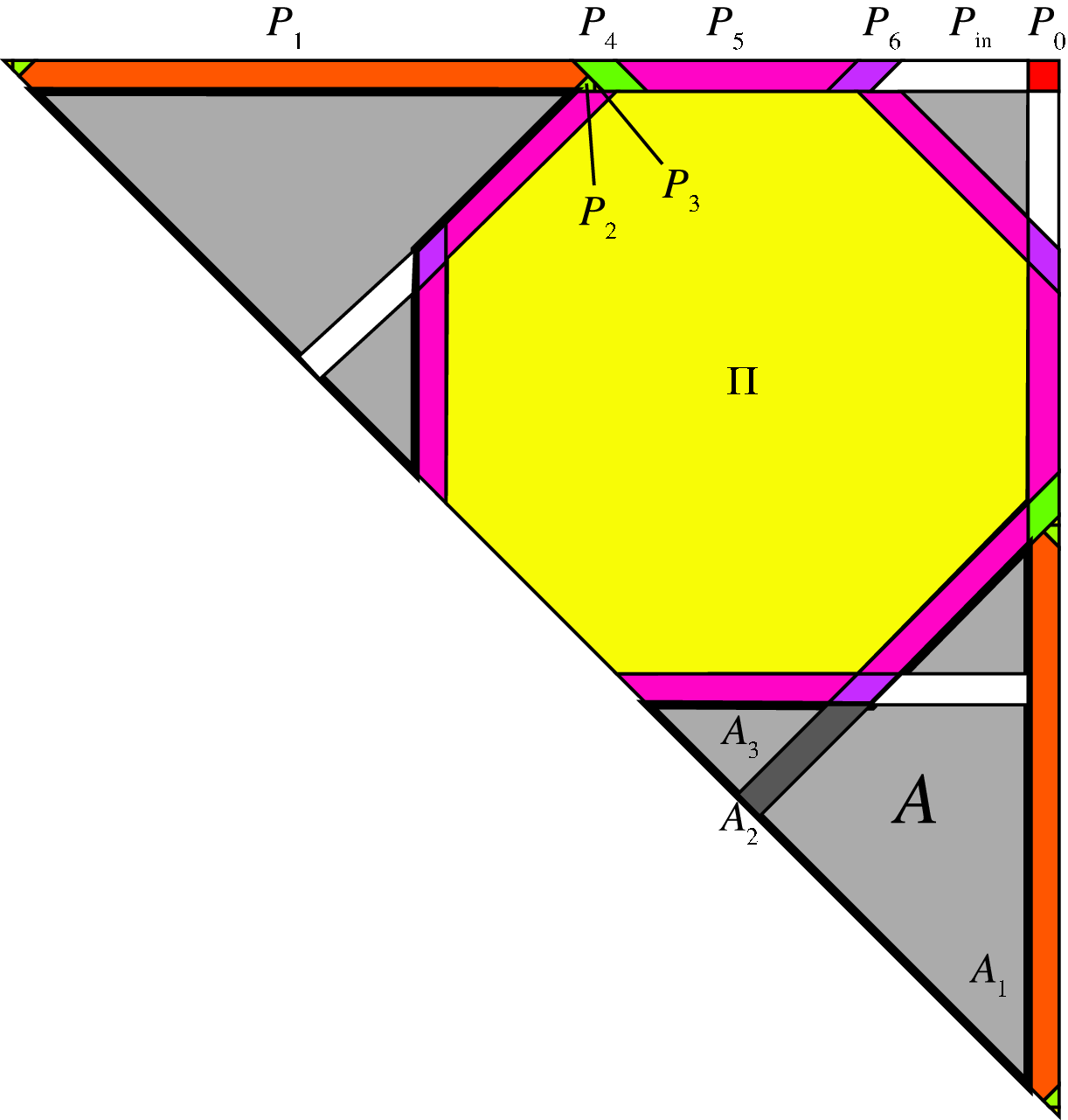,width=5 in}\hfil
\caption{\label{fig:A-7} \small Joint return-map partition of pencil $P$ and arrowhead $A$.  
Those atoms of the pencil's return map which reside in the quadrilateral $P_{\rm in}$ are mapped, 
via three iterations of $\rho_B$, into the entrance tile of the arrowhead $A$.  
The transfer map $\psi(r,h)$ for $A$ then maps them onto the exit tile (dark grey). 
From there, they are mapped back into $P$ by a single application of $\rho_B$.  }
\end{figure}
%%%%%%%%% 

The non-branching part of the proof establishes the first-return orbits of the disconnected 
domain ${\rm P}\cup{\rm A}$.  This set includes the persistent orbits  of $P_0,\ldots,P_6$, 
which begin and end in ${\rm P}$, never entering the arrowhead.  
In addition, we have $P_{\rm in}$ (the complement of $\cup_{i=0}^6 P_i$ 
in ${\rm P}$, ignoring boundaries), whose return orbit ends on the entrance tile of $\cA$.  
The set of atoms is rounded out by the three tiles of ${\bf A}$.  
Of these, $A_1$ and $A_3$ have orbits which return to ${\rm A}$ without 
entering ${\rm P}$, while that of $A_2$ ends up on ${\tt U}_0 P_{\rm in}$.  
All of these statements are proved by direct iteration, with the Tiling Plan 3.
We will use the notation $\rho_{\stackrel{}{\cP,\cA}}$ for the joint return map.

It remains to establish the $s$-dependent partition of $P_{\rm in}$ and piecewise isometric 
mapping of $P_{\rm in}$ under $\rho_{\stackrel{}{\cP}}$.  
The restriction of the latter to $P_{\rm in}$ is given by
\beq\label{eq:Pinmap}
\left.\rho_{\stackrel{}{\cP}}\right\vert_{P_{\rm in}} = 
\rho_{\stackrel{}{\cP,\cA}} \circ \psi_{\stackrel{}{\cA}} \circ  \rho_{\stackrel{}{\cP,\cA}}
\eeq 
where $\psi_\cA$ is the arrowhead transfer map.  
According to the Arrowhead Lemma (section \ref{section:Arrowheads}), 
the map $\psi_{\stackrel{}{\cA}}$ partitions the entrance tile of 
$\cA \sim \widehat{\cA} (\alpha\beta^2,s)$ into a strip congruent 
to $\widehat{\mathbf{S}}(\alpha\beta^2,s)$.
Since $ \rho_{\stackrel{}{\cP,\cA}}$ acts on each atom as an orientation-preserving isometry,  
the atoms contained in $P_{\rm in}$ inherit the strip structure.  
This is precisely what is needed to fill out the specification of the tiling of 
$\cP\sim\widehat{\cP}(\alpha,s)$.  
Finally, it is easy to verify that the composition of maps in (\ref{eq:Pinmap}) 
provides the correct mapping $\rho_{\stackrel{}{\cP}}$ of the atoms in $P_{\rm in}$. 

The temporal scaling information of the induction $\cB\rightarrow \cP$ is 
neatly summarized in an incidence matrix ${\rm M}_\epsilon(\cB\to\cP)$.  
For each of the persistent atoms $P_k,\;k=0,\ldots,6$, one lists in 
the relevant column of ${\rm M}_\epsilon(\cB\to\cP)$ the number of times the return orbit 
of $P_k$ visits each of the five atoms of ${\bf B}$.  
This information is tallied in the course of the CAP validation.    
For the remaining atoms, the same data set provides the tile counts for 
the partial orbits to and from the arrowhead, as well as the the tile counts for 
the return orbits of the arrowhead atoms. Combining this with the 
incidence matrix (\ref{eq:NEk}) for the arrowhead transfer map 
(Arrowhead Lemma, part $iii)$), we obtain the rest of 
${\rm M}_\epsilon(\cB\to\cP)$. The result is

\beq\label{eq:matrixBP}
{\rm M}_\epsilon(\cB\to\cP)=
\begin{array}{c|ccccccccc}
&P_0&P_1&P_2&P_3&P_4&P_5&P_6&\multicolumn{2}{c}{\overbrace{P_{2k-1}\;P_{2k}}^{4\leqslant k\leqslant L}}\\ \hline
B_1&0&1&1-\epsilon&1-\epsilon&0&0&0&a_k&c_k\\
B_2&0&0&0&0&1&6&3&b_k&d_k\\
B_3&1&1&1&1&1&1&1&1&1\\
B_4&0&0&1&0&0&0&0&0&0\\
B_5&0&0&0&1&0&0&0&0&0
\end{array}
\eeq
where $\epsilon=i/|i|$, $L=|i|+3$, and
$$
\begin{array}{ll}
a_k=-\frac{1}{2}+(-1)^k+\frac{1}{2}\cdot3^{k-2},&c_k=-\frac{1}{2}+\frac{1}{4}(-1)^k+\frac{1}{4}\cdot3^{k-2},\\
b_k=1-2(-1)^k+3^{k-2},&d_k=1-\frac{1}{2}(-1)^k+\frac{1}{2}\cdot3^{k-2}.
\end{array}
$$

The return times for the pencil's atoms, expressed in terms of iterations of $\rho_B$, 
are obtained by summing the respective columns of the incidence matrix.  
We have established the completeness of the tiling by return orbits of the pencil 
and arrowhead by calculating the total area of each return orbit (area of the source 
polygon multiplied by the return time), summing over all orbits, including the periodic 
one, and comparing with the area of ${\rm B}$. Note that it is unnecessary to consider 
anew the periodic orbits which pass through the arrowhead, since these are subsumed 
in the recursive tiling property of the arrowhead proved in section \ref{section:Arrowheads}.

\bigskip

\noindent{\sc Proof of (\ref{eq:BP}) and (\ref{eq:matrixBP}) for 
$s\in [\alpha\beta^3,\alpha\beta^2]$}.
On this interval, the pencil $\cP$ is minimal, with precisely 9 atoms.  
Here, the mediation of an arrowhead is not needed, and we prove the renormalization step 
and incidence matrix using the method of direct iteration.  
The statement and validation of the tiling scheme (from which the incidence matrix 
(\ref{eq:matrixBP}) can be verified) may be found in the Electronic 
Supplement \cite{ESupplement}.

\bigskip
 
\noindent {\sc Proof of (\ref{eq:BP}) for $s\in [2\alpha\beta,\alpha)$.} 
Consider now the relation between a value of $s$ in  $(0,\alpha\beta^2]$ and its 
mirror value $\tilde{s}=\alpha-s$. In the latter case, the return-map partition 
of $B_3(\tilde{s})$ (the third atom of $\widehat{\cB}(1,s)$ in the canonical ordering of table \ref{tbl:Bdef}) is very far from being pencil-like.  
On the other hand, one of the five atoms of $B_3(\tilde{s})$ is in fact 
a pencil $\tilde{P}$ with the same parameters and same return-map (up to placement) as $P=B_3(s)$.   

We first compare explicitly the return orbits of tiles $P_2$ and $P_3$ 
to those of $\tilde{P}_2$ and $\tilde{P}_3$.
They are different, but lead to the same image tiles, up to placement, 
after returning to the pencil.  

\begin{figure}[h]
\hfil\epsfig{file=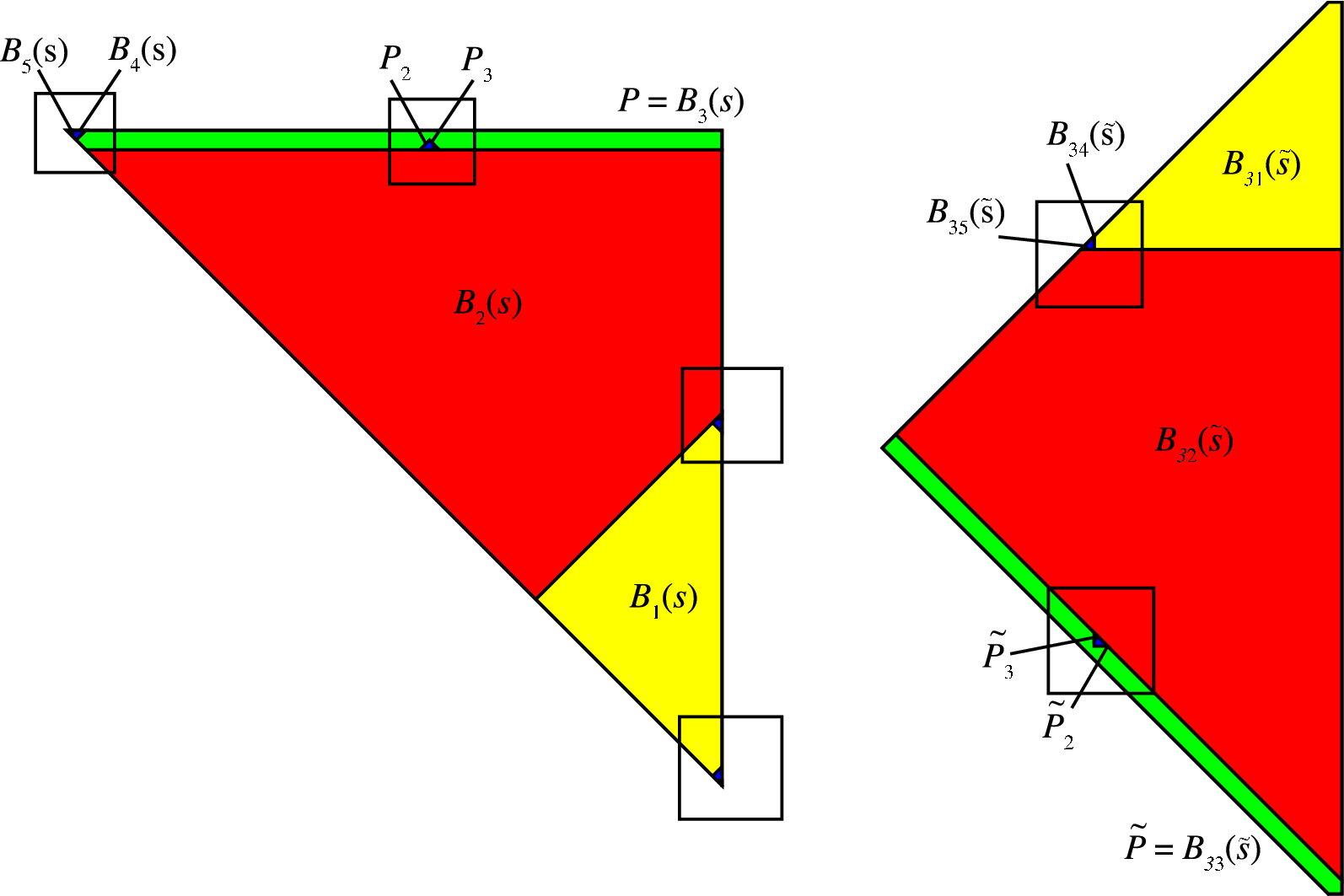,width=5in}\hfil
\caption{\label{fig:ADsymmetry}\small
Return orbits of pencil tiles \#2 and \#3, shown in boxes against the background of the return-map 
partitions of $B_1(s)$ (on left) and $B_{13}(\tilde{s})$, for a value of $\tilde{s}$ near 1.3.  
Note that once the two orbits are deleted, the effective domain maps needed to 
calculate the return orbits of the remaining pencil tiles are conjugate via a rotation and translation.}
\end{figure}

Now let us remove tiles \#2 and \#3 from the game and consider how to construct 
the return orbits of the remaining tiles of each pencil.  
These orbits will visit only the tiles $B_1(s)$ and $B_2(s)$ 
(respectively, $B_{31}(\tilde{s})$ and $B_{32}(\tilde{s})$,
the sub-tiles 1 and 2 of $B_3(\tilde{s})$) before returning to the 
pencil, since the remaining two tiles lie on the orbits of tiles \#2 and \#3.  
Now $B_1(s)$ and $B_{31}(\tilde{s})$ are congruent and the action of the 
respective mappings on them are the same, up to placement.  
If the overlaps of $B_2(s)$ with the orbits of \#2 and \#3 are deleted, 
then the truncated polygon is found to be congruent to $B_{32}(\tilde{s})$, 
and once again the mappings are the same, up to placement.  
Thus the return orbits of all the remaining tiles of the two pencils are 
identical, up to placement.  The above arguments are illustrated in 
figure \ref{fig:ADsymmetry} 
%%%%%%%%%%%%%%%%%%%%%%%%%%%%%%%%%%%%%%%%%%%%%%%% FIGURE 
%%%%%%%%%%%%%%%%%%%%%%%%%%%%%%%%%%%%%%%%%%%%

For the mirror intervals with $s\geqslant 2\alpha\beta$, the $\cB\rightarrow\cP$ incidence matrix is again given by (\ref{eq:matrixBP}), with $\epsilon=1$. 
\clearpage 

\noindent {\sc Proof of (\ref{eq:PP*}).} 
The proof is a direct application of the Pencil Shortening Lemma 
(section \pageref{section:Pencil}).
Multiplying the incidence matrices for $\cB\rightarrow\cP$ and  $\cP\rightarrow\cP^*$, we obtain, for $s\in (0,\alpha\beta^2]\cup [2\alpha\beta,\alpha)$,
$$
{\rm M}_\epsilon(\cB\to\cP^*)=\left(\begin{array}{ccccccccc}
0&-\frac{1}{2}&\frac{1}{2}- \epsilon(i)-2 |i|&-2&-\frac{1}{2}&-\frac{1}{2}&-\frac{1}{2}&-\frac{1}{2}&-\frac{1}{2}\\
0&1&4|i|-6&4&1&1&1&1&1\\
1&1&9&4&1&1&1&1&1\\
0&0&1&0&0&0&0&0&0\\
0&0&0&0&0&0&0&0&0
\end{array}\right) +
$$
\beq\label{eq:matrixBP*}
\left(\begin{array}{ccccccccc}
0&-1&-\frac{3}{4}&-\frac{3}{2}&-\frac{1}{4}&1&\frac{1}{4}&-1&-\frac{1}{4}\\
0&2&\frac{3}{2}&3&\frac{1}{2}&-2&-\frac{1}{2}&2&\frac{1}{2}\\
0&0&0&0&0&0&0&0&0\\
0&0&0&0&0&0&0&0&0\\
0&0&0&0&0&0&0&0&0
\end{array}\right)\, (-1)^i +
\eeq
$$
\left(\begin{array}{ccccccccc}
0&\frac{1}{6}&\frac{7}{12}&\frac{7}{18}&\frac{1}{12}&\frac{1}{2}&\frac{1}{4}&\frac{3}{2}&\frac{3}{4}\\
0&\frac{1}{3}&\frac{7}{6}&\frac{7}{9}&\frac{1}{6}&1&\frac{1}{2}&3&\frac{3}{2}\\
0&0&0&0&0&0&0&0&0\\
0&0&0&0&0&0&0&0&0\\
0&0&0&0&0&0&0&0&0
\end{array}\right)\,3^{|i|}
$$ 
%%%%SUBSECTION
\subsection{Proof of $\cP^*\rightarrow \cB^*$ (scenario II)}
Scenario II  includes both a countable set of accumulation points, namely $\Delta_i,\; |i|\geqslant 2$, 
as well as a subset of the continuity intervals, $I_{i,j}$, $i,j$ finite.  
In the former case, the renormalization is described by statement (ii) of theorem 1, 
while the latter case is handled by the following proposition:

\begin{proposition}\label{prop:P*B}
Let $(i,j)\in{\rm II}\cap\mathbb{Z}^2$, let $s\in I_{i,j}$, 
and let $\cP^*$ be as in (\ref{eq:PP*}).  
Then $\cP^*\rightarrow \cB^*$ where
\beq\label{eq:P*B*}
\cB^*\sim \widehat{\cB}(l^*,r(s )l^*), \quad\mbox{with}\quad l^*=\beta^{|i|+|j|+2},\;\pi(\cB^*)=(-1)^{|i|+|j|}.
\eeq
The incidence matrices for this renormalization step are listed in Appendix D.
\end{proposition}

\begin{figure}[h]
\hfil\epsfig{file=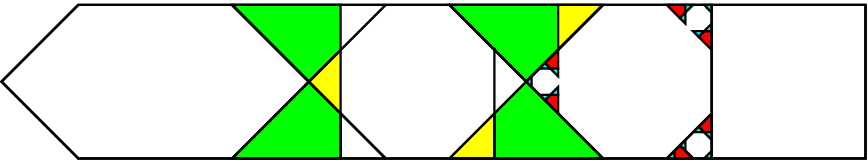,width=5in}\hfil
\caption{\label{fig:accumtiling}\small
Tiling of the minimal pencil by two disjoint return orbits of $s=0$ base 
triangles, plus several periodic orbits.}
\end{figure}

\noindent{\sc Proof.}
Thanks to the scale invariance of the dynamics, it is sufficient to restrict ourselves to the values of $s$
(both the accumulation points and the continuity intervals) where the pencil is minimal.  
These are handled on a case-by-case basis by direct iteration, exactly as was done for scenario I.  
Again, we relegate the tiling plans to the Electronic Supplement, with the exception of those for the 
accumulation points $s=\beta^k,\;  k=2,3,\ldots$, where tiling of the pencil requires 
two return orbits of higher-level base triangles (see Tiling Plan 4 and figure \ref{fig:accumtiling}).   
For the sake of the completeness of the scaling data, we list all the relevant incidence matrices 
in Appendices D and F.

%%%%%%%%%%%%%%%%%%%%%%%%%%%%%%%%%%%%%%%%%%%%%%%% FIGURE 
%%%%%%%%%%%%%%%%%%%%%%%%%%%%%%%%%%%%%%%%%%%%
%
%
%SUBSECTION  10.4
\subsection{Proof of $\mathcal{B}\to\mathcal{T}_\mu$ (scenarios III and IV)}

\begin{proposition}\label{prop:BT}
\noindent 
Let $(i,j)\in {\rm III}\cup {\rm IV}$ (see equation (\ref{eq:ShortHandScenarios})),
let $s \in I_{i,j}$, and let $\cB \sim \widehat{\cB}(1,s)$.  
Then $\cB\rightarrow \cT_{\mu(i,j)}$, where $\mu(i,j)$ is as in table \ref{tbl:scenarios}
p.~\pageref{tbl:scenarios}, and
\begin{eqnarray}\label{eq:BT}
\cT_{-}&\sim& \widehat{\cT}_{-}(\beta^{|i|+1},s-\beta^{|i|+1}),\;\mbox{with}\;\pi(\cT_{-})=(-1)^{|i|},\\
\cT_{+}&\sim& \widehat{\cT}_{+}(\beta^{|i|+2},\beta^{|i|}-s),\;\;\;\mbox{with}\;\;\pi(\cT_{+})=(-1)^{|i|+1}.\nonumber
\end{eqnarray}
The relevant incidence matrices are given in formulae (\ref{eq:matrixBT--})--(\ref{eq:matrixP*T+}).
\end{proposition}

\begin{figure}[h]
\hfil\epsfig{file=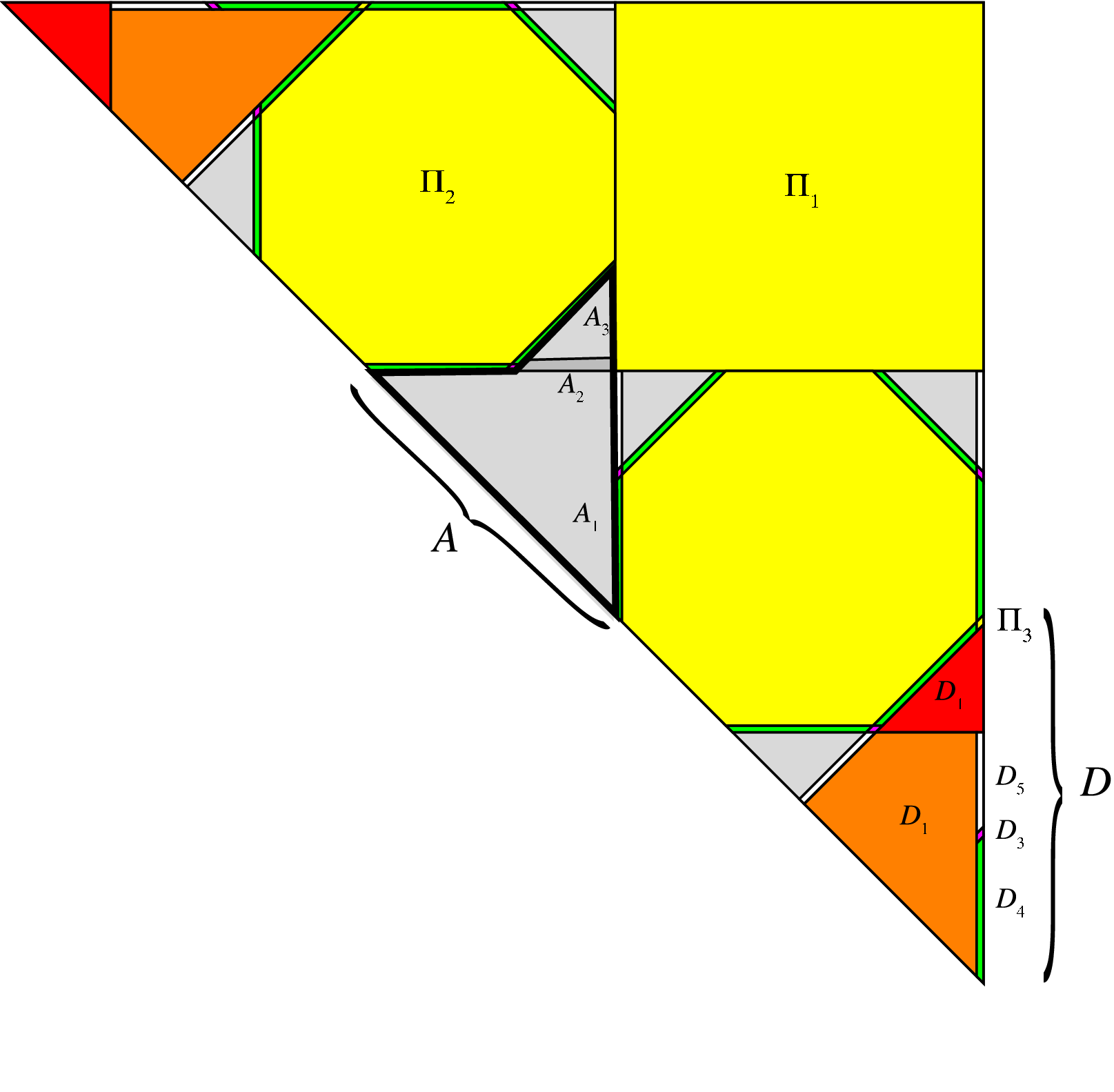,width=5 in}\hfil
\caption{\label{fig:nonbranchIII-} \small Tiling of the base for scenario III$_{-\nu}$,  
according to the Tiling Plan 7.}
\end{figure}

\noindent{\sc Proof.}
We start with scenario III$_{\mu\nu}$, with $s\in [\beta-\beta^6,\beta)\cup(\beta,\beta+\beta^4]$, 
corresponding to $i=0$ for $\cT_-$ and $i= -1$ for $\cT_+$ in (\ref{eq:BT}).
Our first task is to prove that the return map for $B_1$,the first atom of $\widehat{\cB}(1,s)$, promotes the 
triangle to a fringed triangle of type $\widehat{\cT}_-(\beta, s-\beta)$ 
(resp. $\widehat{\cT}_+(\beta^3, \beta-s))$. To help establish this, we introduce 
an auxiliary arrowhead (analogous to that introduced in section \ref{section:ProofOfB-P}) 
and prove by direct iteration that the induced domain map is the correct one.  
The non-branching part of the return-map proof, valid for all $s$ in the chosen 
interval and illustrated in figure \ref{fig:nonbranchIII-},  is accomplished with computer 
assistance, according to the Tiling Plans 5 and 6.  
There ${\rm T}^\pm= B_1$, and $T^-_1,\ldots,T^-_4$ (resp. $T^-_1,\ldots,T^-_7$) are 
atoms whose return orbits are non-branching.  The domain $T^\pm_{\rm in}$ is a `blank' polygon 
which is delivered by the computer program to the entrance tile of the arrowhead.  
Arrowhead dynamics endows the polygon with the ($s$-dependent) structure of a strip, 
and maps it to the exit domain $A_2$, from which it is delivered, isometrically to 
its final destination in ${\rm T}$ by the computer program. The completeness of the 
tiling is checked by summing the areas of all tiles contained in the computer-generated 
orbits of the source domains.  
%%%%%%%%% FIGURE 
%%%%%%%%%%%%%%

An analogous treatment for $s$ in the mirror interval $[1-\beta^4, 1+\beta^6]$ can be given, 
with the tiling domain ${\rm T}^\pm = B_5$, the fifth atom of $\widehat{\cB}(1,s)$. 
The tiling schemes, which are very similar to Tiling Plans 5 and 6, have been relegated 
to the Electronic Supplement.  
Once again, however, we list all the incidence matrices, 
using a superscript to indicate the sign of index $i$. 
Moreover, it is useful, from here on, to introduce a parameter
\begin{equation}\label{eq:K}
K=|j| - (\mu(i,j)+1)/2.
\end{equation}
In the current context, this will be used to keep track of the numbers of atoms in the various parametric dressed domains.

For negative $i$, we have the scenario III incidence matrices
\beq\label{eq:matrixBT--}
{\rm M}_{--}(\cB\to\cT)=
\begin{array}{c|cccc}
&T^-_1&T^-_2&\multicolumn{2}{c}{\overbrace{T^-_{2k-1}\;T^-_{2k}}^{2\leqslant k\leqslant K+2}}\\ \hline
B_1&1&1&1&1\\
B_2&0&0&a_k & c_k\\
B_3&0&0&b_k & d_k\\
B_4&0&1&2&2\\
B_5&1&0&0&0
\end{array}
\eeq
where
$$
\begin{array}{ll}
a_k=1+2(-1)^k+3^{k-1},&\qquad c_k=1+\frac{1}{2}(-1)^k+\frac{1}{2}\cdot3^{k-1},\\
b_k=\frac{1}{2}+3(-1)^k+\frac{1}{2}\cdot3^{k-1},&\qquad
d_k=\frac{1}{2}+\frac{3}{4}(-1)^k+\frac{1}{4}\cdot3^{k-1}.
\end{array}
$$

\beq\label{eq:matrixBT-+}
{\rm M}_{-+}(\cB\to\cT)=
\begin{array}{c|ccccccccc}
&T^+_1&T^+_2&T^+_3&T^+_4&T^+_5&T^+_6&T^+_7&\multicolumn{2}{c}{\overbrace{T^+_{2k}\;T^+_{2k+1}}^{4\leqslant k\leqslant K+4}}\\ \hline
B_1&1&1&1&1&1&1&1&1&1\\
B_2&0&0&0&10&18&2&5&a_k & c_k\\
B_3&0&0&1&7&10&3&4&b_k & d_k\\
B_4&0&1&0&0&0&0&0&0&0\\
B_5&1&0&0&0&0&0&0&0&0
\end{array}
\eeq
where 
$$
\begin{array}{ll}
a_k=17+2(-1)^k+3^{k-1},&\qquad 
c_k=17+\frac{1}{2}(-1)^k+\frac{1}{2}\cdot3^{k-1},\\
b_k=\frac{17}{2}+3(-1)^k+\frac{1}{2}3^{k-1},&\qquad 
d_k=\frac{17}{2}+\frac{3}{4}(-1)^k+\frac{1}{4}\cdot3^{k-1}.
\end{array}
$$

For positive $i$, the corresponding matrices are
\beq\label{eq:matrixBT+-}
{\rm M}_{+-}(\cB\to\cT)=
\begin{array}{c|cccc}
&T^-_1&T^-_2&\multicolumn{2}{c}{\overbrace{T^-_{2k-1}\;T^-_{2k}}^{2\leqslant k\leqslant K+2}}\\ \hline
B_1&1&0&0&0\\
B_2&0&1&2& 2\\
B_3&0&0&a_k & c_k\\
B_4&0&0&b_k & d_k\\
B_5&1&1&1&1
\end{array}
\eeq
where
$$
\begin{array}{ll}
a_k=\frac{1}{2}+5(-1)^k+\frac{1}{2}3^k,&\qquad 
c_k=\frac{1}{2}+\frac{5}{4} (-1)^k+\frac{1}{4}3^k,\\
b_k=1+2(-1)^k+3^{k-1},&\qquad
d_k=1+\frac{1}{2} (-1)^k+\frac{1}{2}\cdot3^{k-1}.
\end{array}
$$
\beq\label{eq:matrixBT++}
{\rm M}_{++}(\cB\to\cT)=
\begin{array}{c|ccccccccc}
&T^+_1&T^+_2&T^+_3&T^+_4&T^+_5&T^+_6&T^+_7&\multicolumn{2}{c}{\overbrace{T^+_{2k}\;T^+_{2k+1}}^{4\leqslant k\leqslant K+4}}\\ \hline
B_1&1&0&0&0&0&0&0&0&0\\
B_2&0&1&0&0&0&0&0&0&0\\
B_3&0&0&2&18&29&6&10&a_k & c_k\\
B_4&0&0&0&10&18&2&5 & b_k & d_k\\
B_5&1&1&1&1&1&1&1&1&1
\end{array}
\eeq
where 
$$
\begin{array}{ll}
a_k=\frac{53}{2}+5(-1)^k+\frac{1}{2}3^{k},&\qquad 
c_k=\frac{53}{2}+\frac{5}{4}(-1)^k+\frac{1}{4}3^{k},\\
b_k=17+2(-1)^k+3^{k-1},&\qquad 
d_k=17+\frac{1}{2}(-1)^k+\frac{1}{2}\cdot3^{k-1}.
\end{array}
$$

We now turn to scenario IV$_{\mu\nu}$, where the renormalization process has already 
reached the minimal pencil stage.  
The next step, in which fringed triangles are induced, is completely analogous to 
what we have just studied  for $\cB\rightarrow \cT_\mu$, again leading to (\ref{eq:BT}). 
The relevant Tiling Plans 7 and 8 are given in Appendix B 
for the two values of $\mu$.

The corresponding incidence matrices are:
\beq\label{eq:matrixP*T-}
{\rm M}_{-}(\cP^*\to\cT)=
\begin{array}{c|ccccccc}
&T^-_1&T^-_2&T^-_3&T^-_4&\multicolumn{2}{c}{\overbrace{T^-_{2k-1}\;T^-_{2k}}^{3\leqslant k\leqslant K+2}}\\ \hline
P^*_0&0&0&0&0&0&0\\
P^*_1&0&0&0&0&0&0\\
P^*_2&0&1&2&2&2&2\\
P^*_3&0&0&0&0&0&0\\
P^*_4&2&2&4&4&4&4\\\
P^*_5&1&1&1&1&1&1\\
P^*_6&0&0&6&3&a_k&d_k\\
P^*_7&0&0&0&0&b_k&e_k\\
P^*_8&0&0&0&0&c_k&f_k\\
\end{array}
\eeq
where 
$$
\begin{array}{ll}
a_k=1+2(-1)^k+3^{k-2},&\qquad d_k=1+\frac{1}{2}(-1)^k+\frac{1}{2}\cdot3^{k-2},\\
b_k=-\frac{1}{2}-(-1)^k+\frac{1}{2}\cdot 3^{k-3},&\qquad e_k=-\frac{1}{2}-\frac{1}{4}(-1)^k+\frac{1}{4}\cdot3^{k-3},\\
c_k=1+2(-1)^k+3^{k-3},&\qquad f_k=1+\frac{1}{2}(-1)^k+\frac{1}{2}\cdot3^{k-3},
\end{array}
$$
\beq\label{eq:matrixP*T+}
{\rm M}_{+}(\cP^*\to\cT)=
\begin{array}{c|cccccccccccc}
&T^+_1&T^+_2&T^+_3&T^+_4&T^+_5&T^+_6&T^+_7&\multicolumn{2}{c}{\overbrace{T^+_{2k}\;T^+_{2k+1}}^{4\leqslant k\leqslant K+4}}\\ \hline
P^*_0&0&0&0&0&0&0&0&0&0\\
P^*_1&0&0&0&0&0&0&0&0&0\\
P^*_2&0&1&0&0&0&0&0&0&0\\
P^*_3&0&0&0&0&0&0&0&0&0\\
P^*_4&2&2&0&0&0&0&0&0&0\\\
P^*_5&1&1&1&1&1&1&1&1&1\\
P^*_6&0&0&0&4&6&2&2&a_k&d_k\\
P^*_7&0&0&0&0&1&0&0&b_k&e_k\\
P^*_8&0&0&0&2&2&0&1&c_k&f_k\\
\end{array}
\eeq
where 
$$
\begin{array}{ll}
a_k=5+2(-1)^k+3^{k-2},&\qquad 
d_k=5+\frac{1}{2}(-1)^k+\frac{1}{2}\cdot3^{k-2},\\
b_k=\frac{3}{2}-(-1)^k+\frac{1}{2}\cdot 3^{k-3},&\qquad 
e_k=\frac{3}{2}-\frac{1}{4}(-1)^k+\frac{1}{4}\cdot3^{k-3},\\
c_k=1+2(-1)^k+3^{k-3},&\qquad 
f_k=1+\frac{1}{2}(-1)^k+\frac{1}{2}\cdot3^{k-3}.
\end{array}
$$

%
%%%SUBSECTION
\subsection{Proof of $\cT_\mu\rightarrow\cD_\mu\rightarrow \cD_\nu^*$ (scenarios III and IV)} 
The next phase of scenarios III$_{\mu\nu}$ and  IV$_{\mu\nu}$ is the transition from fringed 
triangle to double-strip, followed by shortening of the strip. 

\begin{proposition}
Let $(i,j)\in{\rm III}\cup {\rm IV}$, 
let $s \in I_{i,j}$, and let $\cT_{\mu(i,j)}$ be as in (\ref{eq:BT}).
Then $\cT_{\mu(i,j)}\rightarrow \cD_{\mu(i,j)}$, where
\begin{eqnarray}
\cD_{-}&\sim& \widehat{\cD}_{-}(\beta^{|i|+3},s-\beta^{|i|+1}),\;\mbox{with}\;\pi(\cD_{-})=(-1)^{|i|+1},\label{eq:TD}\\
\cD_{+}&\sim& \widehat{\cD}_{+}(\beta^{|i|+4},\beta^{|i|}-s),\;\;\;\mbox{with}\;\;\pi(\cD_{+})=(-1)^{|i|},\nonumber
\end{eqnarray}
and the corresponding incidence matrices are given in (\ref{eq:matrixT+D+}) and (\ref{eq:matrixT-D-}),
Moreover,
\beq\label{eq:DD*}
\cD_{\mu(i,j)}\rightarrow\cD^*_{\nu(j)}\sim\left\{\begin{array}{ll}
\widehat{\cD}_{\nu(j)}(\beta^{|i|+|j|},s-\beta^{|i|+1}),& \;\mbox{if}\; \mu(i,j)=-1,\\ 
\widehat{\cD}_{\nu(j)}(\beta^{|i|+|j|},\beta^{|i|}-s), &\;\mbox{if}\; \mu(i,j)=+1,
\end{array}\right.
\eeq
with $\pi(\cD^*_{\nu(j)})=(-1)^{|i|+|j|}$.
\end{proposition}
\bigskip

\noindent {\sc Proof of (\ref{eq:TD}).}
We begin with the case $\mu(i,j)= -1$. Our task is to calculate the return-map orbits of the 
trapezoidal atom $T^-_3$ of the fringed triangle (see figure \ref{fig:FringedTriangle}), 
showing that the corresponding partition is that of a double-strip with the stated parameters.  
Scale invariance allows us to fix the parameter $l$ of the fringed triangle at $l=\beta$ and 
let the parameter $h$ range over $(0,\beta^4]$.  
This treats simultaneously all integer values of the index $i$.

The non-branching part of the proof consists of the computer-assisted validation 
of tiling plan 9 in Appendix B, illustrated in figure \ref{fig:T-toD-}.  
%%%%%%%%% FIGURE 
\begin{figure}[h]
\hfil\epsfig{file=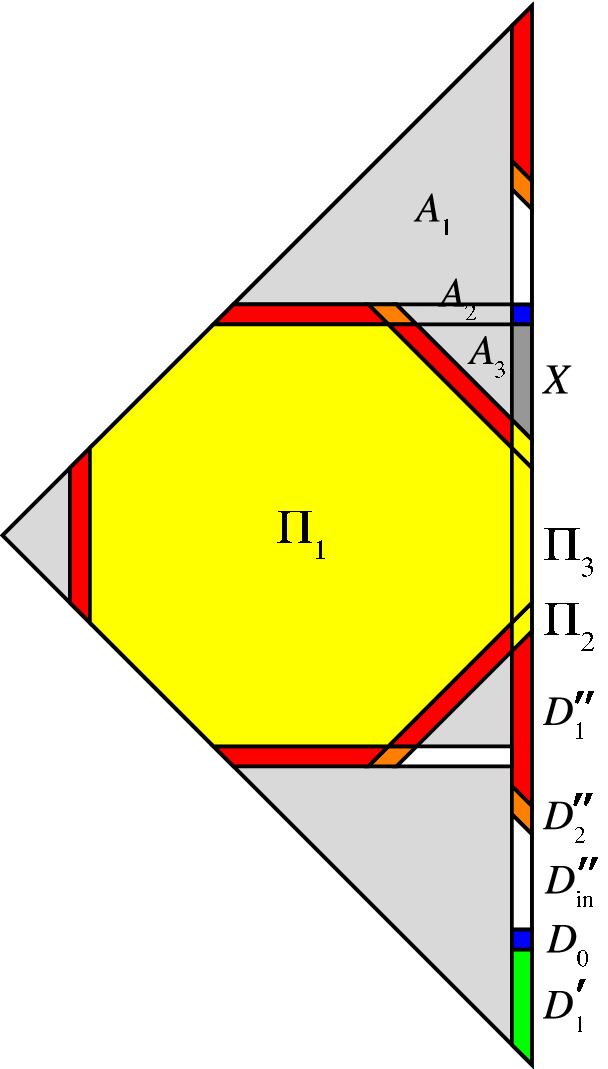,width=2.5 in}\hfil
\caption{\label{fig:T-toD-} \small Tiling of the triangle ${\rm T}_-$, 
according to the Tiling Plan 9.}
\end{figure}
%%%%%%%%%%%%%%
It is worth noting that the induced return map of the double strip corresponds 
to the variant $\nu= -1$ (since the double strip has negative parity and a net 
rotation ${\tt R}_2$ of atom $D_0$). 
  
Once again, arrowhead dynamics is responsible for all of the branching, 
but with an important difference from previous stages of the renormalization.
The double strip, labelled $\cD_-$, has two constituent strips, 
produced by separate $s$-dependent branching mechanisms.  
One of these is associated with the arrowhead, while the other is produced by 
the already established piecewise isometry of proposition \ref{prop:BT}.    
It should be emphasized that, as always, computer assistance is enlisted only 
for the non-branching orbit pieces which are present for all $s$ in the chosen interval.   
This includes the isometric mapping of $D'$ onto $X$, but not the mapping of $X$ 
back into ${\rm D}$ which completes the return orbit. 

The case $\mu(i,j)= +1$ is handled in an analogous manner, with the 
non-branching part of the proof consisting of computer-assisted validation of tiling plan 10.  
There we identify two special tiles congruent to $D'$, namely $X$ and $Y$.  
The isometric mappings from $D'$ to $X$ and from $Y$ to $D'\oslash {\rm D}$ are 
included in the tiling plan, while the mapping from $X$ to $Y$ is a direct application 
of the domain map of the fringed triangle.

As a by-product of the return-map calculations leading to (\ref{eq:TD}), we obtain 
the following incidence matrices:

\beq\label{eq:matrixT-D-}
{\rm M}_{-}(\cT\to\cD)=
\begin{array}{c|cccc}
&D_0&\multicolumn{2}{c}{\overbrace{D^{''}_{2k-1}\; D^{''}_{2k}}^{1\leqslant j\leqslant K}}&\overbrace{D'_l}^{1\leqslant l\leqslant 2K-2}\\ \hline
T^-_1&1&a_k&c_k&0\\
T^-_2&0&b_k&d_k&0\\
T^-_3&1&1&1&1\\
T^-_4\cdots T^-_6&0&0&0&0\\
\underbrace{T^-_{j+6}}_{1\leqslant j\leqslant 2K-2}&0&0&0&\delta_{j l}
\end{array}
\eeq
where 
$$
\begin{array}{ll}
a_k=\frac{3}{2}-(-1)^k+\frac{1}{2}3^{k-1},&\qquad c_k=\frac{3}{2}-\frac{1}{4}(-1)^k+\frac{1}{4}\cdot 3^k\\
 b_k=-1+2(-1)^k+3^{k-1}&\qquad d_k=-1+\frac{1}{2}(-1)^k+\frac{1}{2}\cdot3^{k-1}
\end{array}
$$

\beq\label{eq:matrixT+D+}
{\rm M}_{+}(\cT\to\cD)=
\begin{array}{c|cccc}
&D_0&\multicolumn{2}{c}{\overbrace{D^{''}_{2k-1}\; D^{''}_{2k}}^{1\leqslant j\leqslant K}}&\overbrace{D'_l}^{1\leqslant l\leqslant 2K-2}\\ \hline
T^+_1&0&a_k&c_k&0\\
T^+_2&12&b_k&d_k&12\\
T^+_3&0&0&0&0\\
T^+_4&0&2&2&0\\
T^+_5&1&0&0&0\\
T^+_6&0&0&0&0\\
T^+_7&0&0&0&0\\
T^+_8&1&1&1&1\\
T^+_9\cdots T^+_{11}&0&0&0&0\\
\underbrace{T^+_{j+11}}_{1\leqslant j\leqslant 2K-2}&0&0&0&\delta_{j l}
\end{array}
\eeq
where 
$$
\begin{array}{ll}
a_k=-\frac{1}{2}+(-1)^k+\frac{1}{2}3^{k+2},&\qquad c_k=-\frac{1}{2}+\frac{1}{4}(-1)^k+\frac{1}{4}3^{k+2},\\
b_k=13-2(-1)^k+3^{k+2},&\qquad d_k=13-\frac{1}{2}(-1)^k+\frac{1}{2}\cdot 3^{k+2}.
\end{array}
$$
\noindent {\sc Proof of (\ref{eq:DD*}).}
We apply the Double-Strip Shortening Lemma (section \pageref{section:DoubleStrip}) 
to reduce the number of atoms in the double strip from $4K-1$ to 11 in $k=K-3$ steps.  
Since $K=|j|-(\mu(i,j)+1)/2$, we have
$$
k=|j|-3-\frac{1}{2}(\mu(i,j)+1)=|j| -\left\{\begin{array}{ll}3&\mu(i,j)= -1,\\4&\mu(i,j)=+1.\end{array}\right. 
$$
Since as a result of the shortening, the first argument of $\cD_{\mu}(l,h)$ decreases by a factor $\beta^k$ 
and the parity is multiplied by $(-1)^k$, we see that (\ref{eq:DD*}) follows from (\ref{eq:TD}).
 
%
%%SUBSECTION 
\subsection{Proof of $\cD^*_\nu\to\cB^*$ (scenarios III and IV)}\label{section:ProofOfScenariosIIIandIV}
This is the final stage of scenarios III$_{\mu\nu}$ and IV$_{\mu\nu}$.

\begin{proposition}
Let $(i,j)\subset  {\rm III}\cup {\rm IV}$, 
let $s \in I_{i,j}$, and let $\cD^*_\nu$ be as in (\ref{eq:DD*}).  Then 
\beq\label{eq:D*B*}
\cD^*_\nu\rightarrow \cB^* \sim \widehat{\cB}(l^*,r(s) l^*), \;\mbox{with}\;\;l^*=\beta^{|i|+|j|+2},\; \pi(\cB^*)=(-1)^{|i|+|j|}.
\eeq
The corresponding incidence matrices are given in (\ref{eq:matrixD*B*}).
\end{proposition}

\begin{figure}[h]
\hfil\epsfig{file=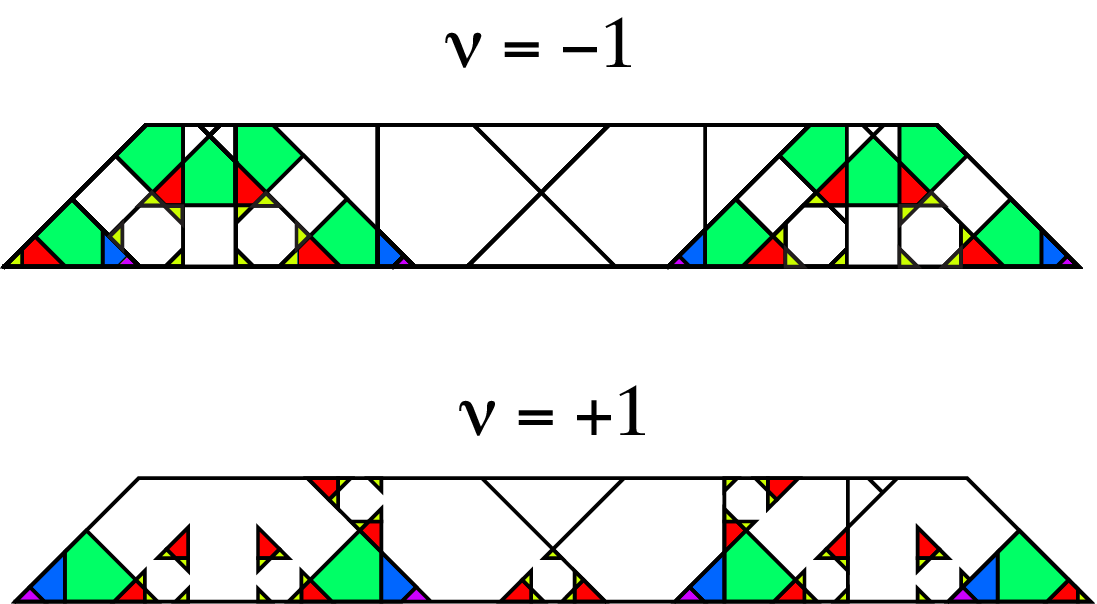,width=5 in}\hfil
\caption{\label{fig:D*B*} \small Tiling of the minimal double strip 
$\cD*_\nu$ for $s=87/200$ ($\nu=-1$) and $s=211/500$ ($\nu=+1$).}
\end{figure}

\noindent {\sc Proof}. 
Setting
$$
h=\left\{\begin{array}{ll}s-\beta^{|i|+1}&\mu(i,j)= -1,\\ \beta^{|i|} - s& \mu(i,j)=+1,\end{array}\right.
$$
we rescale both arguments of $\widehat{\cD}_{\nu}$ in (\ref{eq:DD*}) by a factor $\beta^{3-|i|-|j|}$ to obtain 
$$
\widehat{\cD}_{\nu}(\beta^3,h'),\quad h'\in [\beta^5,\beta^4],
$$
$$
h'=\beta^{3-|i|-|j|} h.
$$
For both values of $\nu$, we have shown by direct iteration (see tiling plans 11 and 12 and figure \ref{fig:D*B*}) that the induced return map on the tile $D'_1$ of the minimal double strip is that of a base triangle of the same parity, namely
$$
\widehat{\cD}_{\nu}(\beta^3, h')\rightarrow \cB_{\nu} \sim \widehat{B}(\beta^5,\beta^4-h').
$$
To complete the proof, we undo the scale transformation, multiplying both arguments of $\widehat{\cD}_{\nu}$ and
$\widehat{B}$ by the same factor $\beta^{|i|+|j|-3}$, and returning to the variable $s$.  
The result is equation (\ref{eq:D*B*}). \/ $\Box$
%%%%%%%%% FIGURE 
%%%%%%%%%%%%%%
\clearpage

The incidence matrices for $\cD^*_\pm \rightarrow \cB^*$ are:
\beq\label{eq:matrixD*B*} 
\begin{array}{ccc}
{\rm M}_{-}(\cD^*\to\cB^*)= & \qquad\qquad & {\rm M}_{+}(\cD^*\to\cB^*)=\\ \\
\begin{array}{c|ccccc}
&B^*_1&B^*_2&B^*_3&B^*_4&B^*_5\\ \hline
D^*_0&10&4&1&0&0\\
D^{*''}_1&10&4&2&2&2\\
D^{*''}_2&10&4&0&0&0\\
D^{*''}_3&0&0&0&0&0\\
D^{*''}_4&0&0&0&0&0\\
D^{*''}_5&0&0&0&0&1\\
D^{*''}_6&0&0&0&1&0\\
D^{*'}_1&1&1&1&1&1\\
D^{*'}_2&4&1&0&0&0\\
D^{*'}_3&0&0&0&0&0\\
D^{*'}_4&0&0&0&0&0
\end{array} &

\qquad\qquad &

\begin{array}{c|ccccc}
&B^*_1&B^*_2&B^*_3&B^*_4&B^*_5\\ \hline
D^*_0&5&2&2&0&0\\
D^{*''}_1&10&4&5&2&2\\
D^{*''}_2&0&0&0&0&0\\
D^{*''}_3&0&0&0&0&0\\
D^{*''}_4&0&0&0&0&0\\
D^{*''}_5&0&0&0&0&1\\
D^{*''}_6&0&0&0&1&0\\
D^{*'}_1&1&1&1&1&1\\
D^{*'}_2&4&1&2&0&0\\
D^{*'}_3&0&0&0&0&0\\
D^{*'}_4&0&0&0&0&0
\end{array}
\end{array}
\eeq

\subsection{Conclusion of the proof}

We have now proved all of the induction steps comprising the renormalization 
graph of figure \ref{fig:RenormGraph}.  
By composing the relevant return maps along the ten closed circuits of the 
graph, we can now assign to each of the $s$-parametric intervals $I_{i,j}$ a 
renormalization of the class of base triangles congruent to 
$\widehat{\cB}(l,h)$, with $s=h/l\in I_{i,j}$, renormalization 
functions $\kappa$ and $r$, and parities as prescribed by theorem 1.  

The proof of each induction step, via either the Shortening Lemmas of 
section \ref{section:PrototypeDressedDomains} or the propositions of 
section \ref{section:ProofOfTheorem1}, yield as by-products a verification 
of a uniform return path for each pair $(i,j)$, as well as explicit 
formulae for the incidence matrices.

By composing the tilings of the individual induction steps, one obtains 
for each $(i,j)$ a tiling of the parent base triangle.  
This tiling consists of the return orbits of the atoms of the higher-level 
base triangle, as well as those of finitely many periodic tiles explicitly 
identified in the proof. For any $(i,j)\in\Z^2$, the tile coordinates are 
all affine functions of $s$.  
The exhaustiveness of the tiling for each step, hence for the full 
renormalization, has been established by evaluation of area sums, 
included in the Electronic Supplement \cite{ESupplement}.

The proof of theorem 1 is now complete. $\Box$

%%%SECTION
\section{Incidence matrices}\label{section:IncidenceMatrices}

Theorem 1 contains a statement asserting the constancy of the incidence matrix 
on each parameter interval $I_{i,j}$. We did not provide explicit expressions
for the incidence matrices, due to the inherently complicated nature of the return-map dynamics. 
In section \ref{section:ProofOfTheorem1}, we saw how the ten qualitatively distinct 
routes from initial $\cB$ to final $\cB^*$ all led to a simple spatial scaling 
factor $\beta^{|i|+|j|+2}$ and a simple parity-changing factor $(-1)^{|i|+|j|}$.  
No comparable simplicity can be extracted from the incidence matrices sprinkled 
throughout the proof of section \ref{section:ProofOfTheorem1},  
as they are strongly scenario-dependent.
In particular, the incidence matrices vary from interval to interval in a complicated way, 
breaking the symmetry under $(i,j)\mapsto (-i,-j)$ enjoyed by the spatial scaling and 
parity-changing factors.

In section \ref{section:ProofOfTheorem1}, we computed the incidence matrices associated 
with all of the edges of the renormalization graph. For scenarios I and II, we displayed 
(in appendices C and D) the explicit $5\times 5$ matrices associated with a number of 
individual intervals and sequences of intervals. On the other hand, for scenarios 
III and IV, we presented the relevant incidences matrices for individual induction 
steps, but made no attempt at assembling these into all-embracing formulae.  
The goal of the present section is to derive such formulae and examine some of 
their implications.

\subsection{Composite formulae for scenarios III and IV}
We begin by deriving a concise representation for the incidence matrices associated with 
the composite renormalization step $\cB\to\cD_\mu$.  
Similar treatment will be given for $\cP^*\to\cD_\mu$.

\begin{proposition}
Let $(i,j)\in {\rm III}$,
let $s \in I_{i,j}$, and let ${\rm M} = \mathrm{M}^{\mathrm{III}}_{\epsilon,\mu}(\cB\to\cD)$.   Further, let $K$ be given by (\ref{eq:K}).
Then for  $k=1,\ldots,K,\; l=1,\ldots,K-1$, the column vectors ${\rm M}_{\mathbf{B},D}$ 
of the matrix ${\rm M}$ (with $\mathbf{B}=(B_1,\ldots,B_5)^T$) may be written as follows:
\begin{eqnarray}
{\rm M}_{{\bf B},D_0}&=&{\bf X}_0,\nonumber\\
{\rm M}_{{\bf B},D^{''}_{2k-1}}&=&{\bf X}^{''} +4\times (-1)^k\,{\bf Y}^{''}  +2\times 3^k\,{\bf Z}^{''},\nonumber\\
{\rm M}_{{\bf B},D^{''}_{2k}}&=&{\bf X}^{''} + (-1)^k\,{\bf Y}^{''}  +3^k\,{\bf Z}^{''},\label{eq:M_XYZ} \\
{\rm M}_{{\bf B},D^{'}_{2k-1}}&=&{\bf X}^{'} + 4\times (-1)^k\,{\bf Y}^{'}  +2\times 3^k\,{\bf Z}^{'} ,\nonumber \\
{\rm M}_{{\bf B},D^{'}_{2k}}&=&{\bf X}^{'} +(-1)^k\, {\bf Y}^{'}  + 3^k\,{\bf Z}^{'}\nonumber .
\end{eqnarray}
\end{proposition}

\noindent {\sc Proof.}  One verifies the formulae by explicitly calculating
$$
{\rm M}={\rm M}_{\epsilon,\mu}(\cB\to\cT)\cdot  {\rm M}_{\mu}(\cT\to\cD).
$$
The results are listed in the following tables:

$$
\renewcommand{\arraystretch}{1.2}
\begin{array}{|c|c|c|c|c|c|c||c|c|c|c|c|c|c|}\hline
\multicolumn{7}{|c||}{\epsilon=-1, \mu=-1}&\multicolumn{7}{|c|}{\epsilon=-1, \mu=+1}\\ \hline
{\bf X}_0 &{\bf X}^{''}&{\bf Y}^{''}&{\bf Z}^{''}&{\bf X}^{'}&{\bf Y}^{'}&{\bf Z}^{'}&{\bf X}_0 &{\bf X}^{''}&{\bf Y}^{''}&{\bf Z}^{''}&{\bf X}^{'}&{\bf Y}^{'}&{\bf Z}^{'}  \\ \hline
2&\frac{3}{2} & -\frac{1}{4}&\frac{3}{4}&2&0&0  & 14& \frac{31}{2}& -\frac{1}{4}&\frac{27}{4}&14&0&0\\
6&6&0&0&7&-\frac{1}{2}&\frac{9}{2} & 64&66&0&0&63&-\frac{1}{2}&\frac{81}{2}\\
5&5&0&0&\frac{11}{2}&-\frac{3}{4}&\frac{9}{4}  & 35&39&0&0&\frac{67}{2}&-\frac{3}{4}&\frac{81}{4}\\
2&1&-\frac{1}{2}&\frac{1}{2}&4&0&0  &  12&13&-\frac{1}{2}&\frac{9}{2}&12&0&0\\
1&\frac{3}{2}&\frac{1}{4}&\frac{1}{4}&0&0&0  &  0&-\frac{1}{2}&\frac{1}{4}&\frac{9}{4}&0&0&0\\
\hline
\end{array}
$$

$$
\renewcommand{\arraystretch}{1.2}
\begin{array}{|c|c|c|c|c|c|c||c|c|c|c|c|c|c|}\hline
\multicolumn{7}{|c||}{\epsilon=+1, \mu=-1}&\multicolumn{7}{|c|}{\epsilon=+1, \mu=+1}\\ \hline
{\bf X}_0 &{\bf X}^{''}&{\bf Y}^{''}&{\bf Z}^{''}&{\bf X}^{'}&{\bf Y}^{'}&{\bf Z}^{'}&{\bf X}_0 &{\bf X}^{''}&{\bf Y}^{''}&{\bf Z}^{''}&{\bf X}^{'}&{\bf Y}^{'}&{\bf Z}^{'}  \\ \hline
1&\frac{3}{2} & \frac{1}{4}&\frac{1}{4}&0&0&0  & 0&- \frac{1}{2}& \frac{1}{4}&\frac{9}{4}&0&0&0\\
2&1&-\frac{1}{2}&\frac{1}{2}&4&0&0  & 12&13&-\frac{1}{2}&\frac{9}{2}&12&0&0\\
10&10&0&0&\frac{21}{2}&-\frac{5}{4}&\frac{27}{4}  & 101&108&0&0&\frac{197}{2}&-\frac{5}{4}&\frac{243}{4}\\
6&6&0&0&7&-\frac{1}{2}&\frac{9}{2}  &  64&66&0&0&63&-\frac{1}{2}&\frac{81}{2}\\
2&\frac{3}{2}&-\frac{1}{4}&\frac{3}{4}&2&0&0  &  14&\frac{31}{2}&-\frac{1}{4}&\frac{27}{4}&14&0&0\\
\hline
\end{array}
$$

We now calculate, using (\ref{eq:shortModd})--(\ref{eq:shortNeven}), and (\ref{eq:M_XYZ}), the matrix $ \tilde{\rm M}$ and column vector ${\rm N}$ such that
$$
{\rm M}^{\rm III}_{\epsilon,\mu,\nu}(\cB\to\cD^*)={\rm M}^{\rm III}_{\epsilon,\mu}(\cB\to\cD)\cdot  {\rm M}_{\mu,\nu}(\cD\to\cD^*)= 
  \tilde{\rm M} +{\rm N}\cdot (\underbrace{1,1,\ldots,1}_{11}).
$$
where, for $K$ even, we have
$$
\begin{array}{lll}
\tilde{\rm M}_{{\bf B},D^*_0}&=&{\rm M}_{{\bf B},D_0}={\bf X}_0,\\
\noalign{\vskip 1pt}
\tilde{\rm M}_{{\bf B},D^{*''}_k}&=&{\rm M}_{{\bf B},D'_{2 K+k-8}}\\
&=&
     \left\{\begin{array}{ll}
     {\bf X}'+4\times (-1)^{\frac{k-7}{2}}\,{\bf Y}'+2\times 3^{K+\frac{k-7}{2}}\,{\bf Z}',& k=1,3,5,\\ 
      {\bf X}'+(-1)^{\frac{k-8}{2}}\,{\bf Y}'+3^{K+\frac{k-8}{2}}\,{\bf Z}',& k=2,4,6,
     \end{array}\right.\\
\noalign{\vskip 1pt}
\tilde{\rm M}_{{\bf B},D^{*'}_l}&=&{\rm M}_{{\bf B},D^{''}_{2 K+l-4}}\\
&=&
      \left\{\begin{array}{ll}
     {\bf X}^{''}+4\times(-1)^{\frac{l-7}{2}}\, {\bf Y}^{''}+ 2\times 3^{K+\frac{l-3}{2}}\,{\bf Z}^{''},& l=1,3,\\ 
      {\bf X}^{''}+(-1)^{\frac{l-8}{2}}{\bf Y}^{''}+ 3^{K+\frac{l-4}{2}}\,{\bf Z}^{''},&l=2,4,
     \end{array}\right.
\end{array}
$$
$$
{\rm N}_{\bf B}=\xi_{K+1}({\bf X}^{''}-4\,{\bf Y}^{''}) + 2 \,\eta_{K+1}{\bf Z}^{''} + \xi_{K-1}(2\,{\bf X}^{'}-8\, {\bf Y}^{'}) + 4\,\eta_{K-1}{\bf Z}^{'},
$$
where, using the summation formula for geometric series,
$$
\xi_K = \sum^{\frac{K-5}{2}}_{l=0} 4^l= -\frac{1}{3} +\frac{1}{24} \, 2^K,\qquad
\eta_K = \sum^{\frac{K-5}{2}}_{l=0} 4^l \, 3^{K-4-2 l} = \frac{1}{45}\,3^K - \frac{3}{40}\, 2^K.
$$

For $K$ odd, we find
$$
\begin{array}{lll}
\tilde{\rm M}_{{\bf B},D^*_0}&=&{\rm M}_{{\bf B},D_0}={\bf X}_0,\\
\noalign{\vskip 1pt}
\tilde{\rm M}_{{\bf B},D^{*''}_k}&=&{\rm M}_{{\bf B},D^{''}_{2 K+k-6}}\\
&=&
     \left\{\begin{array}{ll}
     {\bf X}^{''}-4\times(-1)^{\frac{k-5}{2}}\,{\bf Y}^{''}+2\times 3^{K+\frac{k-5}{2}}\,{\bf Z}^{''},& k=1,3,5,\\ 
      {\bf X}^{''}- (-1)^{\frac{k-6}{2}}\,{\bf Y}^{''}+ 3^{K+\frac{k-6}{2}}{\bf Z}^{''},& k=2,4,6,
     \end{array}\right.\\
\noalign{\vskip 1pt}
\tilde{\rm M}_{{\bf B},D^{*'}_l}&=&{\rm M}_{{\bf B},D^{'}_{2 K+l-6}}\\
&=&
      \left\{\begin{array}{ll}
     {\bf X}^{'}-4\times (-1)^{\frac{l-5}{2}}\,{\bf Y}^{'}+4\times 3^{K+\frac{l-5}{2}}\,{\bf Z}^{'}\,,& l=1,3,\\ 
      {\bf X}^{'}-(-1)^{\frac{l-5}{2}}\,{\bf Y}^{'}+3^{K+\frac{l-6}{2}}\,{\bf Z}^{'},& l=2,4,
     \end{array}\right.\\
\noalign{\vskip 1pt}
{\rm N}_{\bf B}&=&\xi_{K}(2\,{\bf X}^{'}-8\,{\bf Y}^{'}+{\bf X}^{''} 
          -4\,{\bf Y}^{''}) + \eta_{K}(4\,{\bf Z}^{'} + 2 \,{\bf Z}^{''}).
\end{array}
$$

The above formulae can also be used for calculating incidence matrices for $\cP^*\rightarrow \cD^*$, 
with the index set ${\bf B}$ replaced by ${\bf P}^*$ and the vectors ${\bf X}_0$, ${\bf X}^{'}_k$, etc.,
now having nine components, taken from the following tables. 
$$
\renewcommand{\arraystretch}{1.2}
\begin{array}{|c|c|c|c|c|c|c||c|c|c|c|c|c|c|}\hline
\multicolumn{7}{|c||}{\mu=+1}&\multicolumn{7}{|c|}{\mu=-1}\\ \hline
{\bf X}_0 &{\bf X}^{''}&{\bf Y}^{''}&{\bf Z}^{''}&{\bf X}^{'}&{\bf Y}^{'}&{\bf Z}^{'}&{\bf X}_0 &{\bf X}^{''}&{\bf Y}^{''}&{\bf Z}^{''}&{\bf X}^{'}&{\bf Y}^{'}&{\bf Z}^{'}  \\ \hline
0&0&0&0&0&0&0&0&0&0&0&0&0&0\\
0&0&0&0&0&0&0&0&0&0&0&0&0&0\\
12&13&-\frac{1}{2}&\frac{9}{2}&12&0&0  &  2&1&-\frac{1}{2}&\frac{1}{2}&4&0&0\\
0&0&0&0&0&0&0&0&0&0&0&0&0&0\\
24&25&-\frac{1}{2}&\frac{27}{2}&24&0&0  &  6&5&-\frac{1}{2}&\frac{3}{2}&8&0&0\\
14&\frac{31}{2}&-\frac{1}{4}&\frac{27}{4}&14&0&0  &  2&\frac{3}{2}&-\frac{1}{4}&\frac{3}{4}&2&0&0\\
22&24&0&0&21&-\frac{1}{2}&\frac{27}{2}  &  6&6&0&0&7&-\frac{1}{2}&\frac{3}{2}\\
3&2&0&0&\frac{7}{2}&\frac{1}{4}&\frac{9}{4}  &  0&0&0&0&-\frac{1}{2}&\frac{1}{4}&\frac{1}{4}\\
8&10&0&0&7&-\frac{1}{2}&\frac{9}{2}  &  0&0&0&0&1&-\frac{1}{2}&\frac{1}{2}\\
\hline
\end{array}
$$
Combining the above results with the matrices ${\rm M}_{\nu}(\cD^*\to \cB^*)$ of (\ref{eq:matrixD*B*}), 
we have calculated the composite incidence matrices for $\cB\to \cB^*$ (scenario III) and 
$\cP^*\to\cB^*$ (scenario IV), and listed them in Appendices E and F, respectively.

\subsection{Hausdorff dimensions for selected fixed points of $r$} 
We conclude this section with a brief discussion of the fractal properties of the 
exceptional set (complementary to all periodic orbits). 
This is an application of the incidence matrix formulae derived above. 
For reasons of space, our analysis will be limited to a single case:
scenario $\mathrm{IV}_{-+}$, with parameter intervals $I_{i,j}$, $i<0,j<-2$, $j 
\mbox{ even}$.  
The incidence matrix for this infinite family of intervals is
\beq\label{eq:Mij}
\mathrm{M}(i,j)=\mathrm{M}_{-}(\cB\to\cP^*)\cdot \mathrm{M}^{IV}_{- +}(\cP^*\to\cB^*),
\eeq
with the first factor (a function of $i$) taken from (\ref{eq:matrixBP*}),  
and the second (a function of $j$) from Appendix F.

We further restrict our attention to the simplest strictly renormalizable cases, 
namely the fixed points of $r(s)$,
$$
s_{\mathrm{fix}}(i,j)=\frac{\beta^{|i|+1}(1+\beta^{|j|})}{1+\beta^{|i|+|j|+2}},
$$
with the period-2 symbolic representation $(i,-j,i,-j,\ldots)$.  
For each such parameter value, the temporal scaling factor is just the largest eigenvalue 
$\tau(i,j)$ of the incidence matrix $\mathrm{M}(i,j)$. 
The spatial and temporal scale factor then, in standard fashion \cite{Falconer,Lowenstein}, 
determine the Hausdorff dimension of the exceptional set through
\beq
d_H(i,j)=-\frac{\log(\tau(i,j))}{\log(\beta^{|i|+|j|+2})}.
\eeq
The quintic eigenvalue equation for $d_H(i,j)$ has two trivial solutions, $0$ and $-1$, so that $\tau(i,j)$ 
and so the calculation of $\tau(i,j)$ reduces to finding the largest root of a cubic polynomial.

\begin{figure}[t]
\hfil\epsfig{file=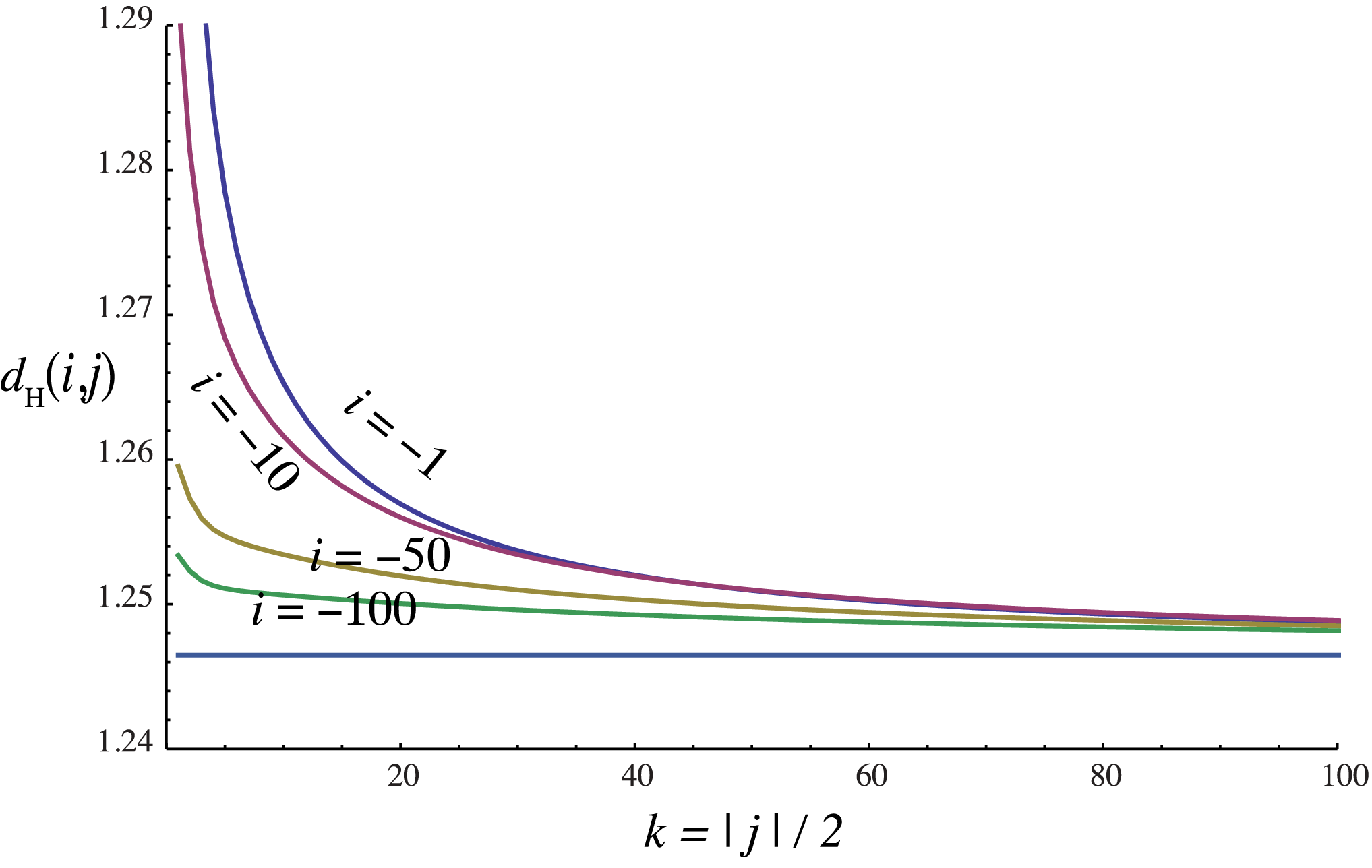,width=5 in}\hfil
\caption{\label{fig:dH} Plots of $d_H(i,j)$ versus $k=|j|/2$ for $i= -1$ (top curve),$-10,-50,$ and $-100$ and 
$2\leqslant k\leqslant 100$.  
Calculated points have been joined by straight segments to aid the eye. 
The horizontal asymptote corresponds to the value $d_\infty$.\small }  
\end{figure}

We have performed the numerical calculation of the dimension for four values of $i$ and $100$ 
values of $j$ in the designated range.  
The results, together with the limiting value for $j\rightarrow -\infty$, 
namely $d_{\infty}=\log{3}/\log{\omega}=1.24648\ldots$, are plotted in figure \ref{fig:dH}.  
Note that the dimension exceeds $d_{\infty}$ for all $i$ and $j$ in our index set, tending 
monotonically to $d_{\infty}$ when either $i$ or $j$ tends to $-\infty$.

%%%%%%%%% FIGURE 
%%%%%%%%%

%%%%%%%%%%%%%%%%%%%%%%%%%%

\clearpage
 
\section*{Appendix A: Catalogue of polygonal shapes}
\begin{table}[h]
\caption{\label{tbl:polytypes}\small Polygonal shapes encountered in the text. All of these are convex with the exception of $\widehat{\mathrm{A}}$, the arrowhead.}
$$
\begin{array}{|c|c|}
\hline
Q_1(l)&[(6,1,3),(l,0,0)]\\ \hline
Q_2(l)&[(4,6,0,2),(l,l,l,l)]\\ \hline
Q_3(l)&[(6,7,2,3),(l,l,l,l)]\\ \hline
Q_4(l_1,l_2)&[(6,1,2,3),(l_1,0,-l_2,0)]\\ \hline
Q_5(l)&[(0,1,\ldots,7),(l,l,\ldots,l)]\\ \hline
Q_6(l_1,l_2)&[(5,6,0,2,3),(0,l_2,l_1,l_2,0)]\\ \hline
Q_7(l)&[(5,7,0,2),(0,l,l,0)]\\ \hline
Q_8(l_1,l_2)&[(4,5,7,0,1,2),(0,0,l_1,l_2,l_2,l_1)]\\ \hline
Q_9(l_1,l_2)&[(1,2,3,5,6,7),(l_1,l_2,l_1,l_1,l_2,l_1)]\\ \hline
Q_{10}(l_1,l_2)&[(4,6,0,2),(l_2,l_1,l_2,l_1)]\\ \hline
Q_{11}(l_1,l_2,l_3)&[(1,2,3,5,6,7),(l_1,l_2,l_1,l_1,l_2,l_1)]\\ \hline
Q_{12}(l_1,l_2)&[(6,0,2,3),(0,l_1,l_2,0)]\\ \hline
Q_{13}(l_1,l_2)&[(5,6,0,2),(0,l_2,l_1,0)]\\ \hline
Q_{14}(l_1,l_2)&[(0,1,\ldots,7),(l_1,l_1,l_2,l_2,l_1,l_1,l_2,l_2)]\\ \hline
Q_{15}(l_1,l_2)&[(0,1,\ldots,7),(l_1,l_2,l_1,l_2,l_1,l_2,l_1,l_2)]\\ \hline
\widehat{\mathrm{A}}(l)&[(5,0,2),(l,l,0)]\cup[(5,0,3),(l,l,0)]\\ \hline
\end{array}
$$
\end{table}
%%%%%%%%% FIGURE 
\begin{figure}[h]
\hfil\epsfig{file=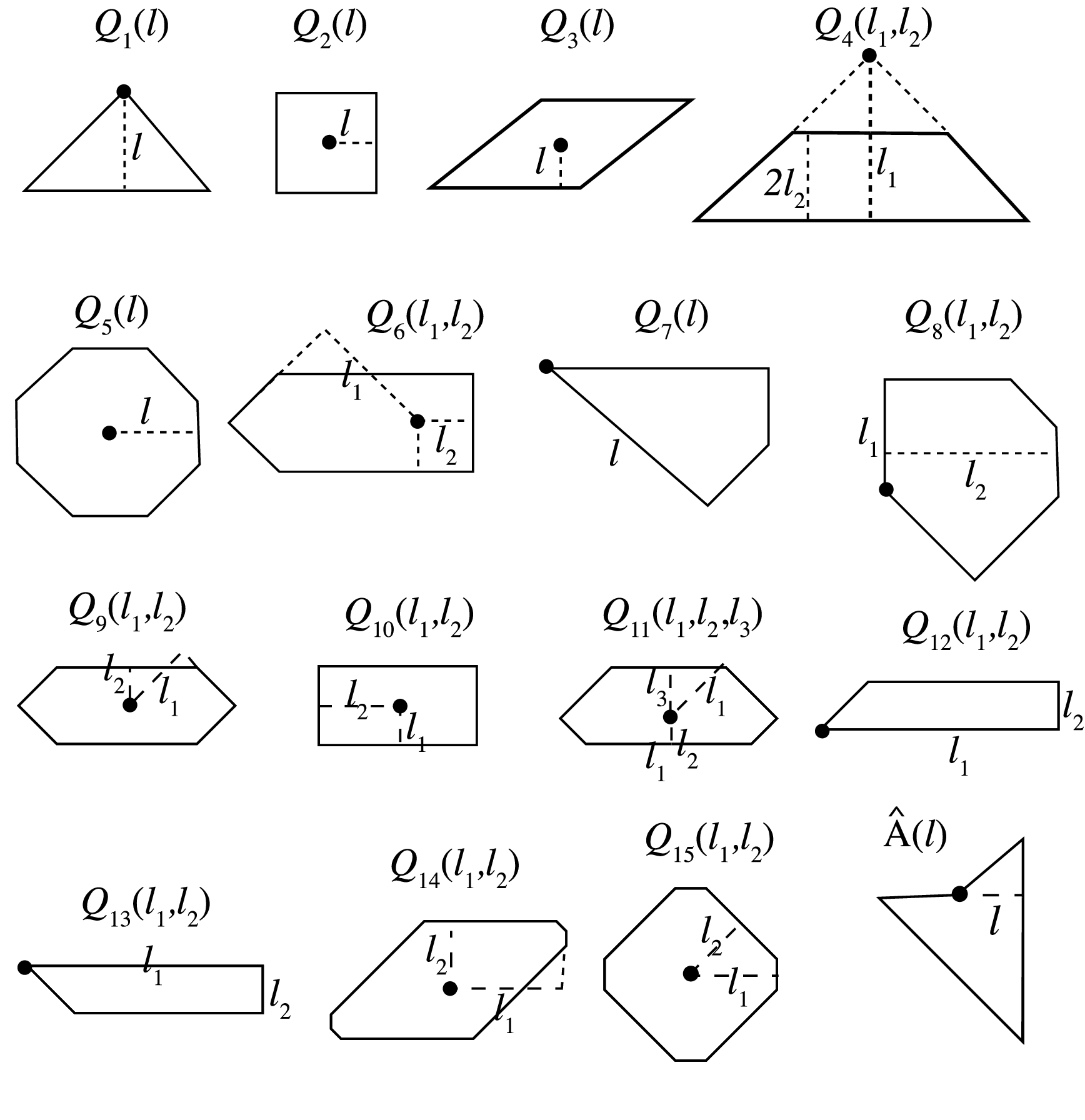,width=5 in}\hfil
\caption{\label{fig:polytypes} \small Polygonal shapes listed in Table \ref{tbl:polytypes}, 
with anchor points emphasized.}
\end{figure}
%%%%%%%%%%%%%%
\clearpage

\section*{Appendix B: Tiling plans}
Key:
\begin{itemize}
\item  Heading, `Parent dressed domain': \ \  PWI for constructing orbits by direct iteration.

\item Heading, `Tiling domains':\ \  spans of listed source tiles, expressed in terms of the representative polygons of Appendix A.

\item Col. 1, `Tile':\ \  name of source tile whose orbit is to be calculated.

\item Cols 2, 3,  `$Q_{\#}$', `Parameters': \ \  
source tile representative polygon and parameters (see Appendix A).

\item Cols 4, 5,  `$\mathtt{R}_{\#}$', `Translation':\ \ index $n$ of rotation $\mathtt{R}_{n}$, and translation $(d_x,d_y)$ of source tile relative to representative polygon.

\item Col 6, 7:  `tile', `$\mathtt{R}_{\#}$':\ \  destination tile and index $n$ of net rotation $\mathtt{R}_{n}$ relative to source tile. These specify the net isometry to be verified by direct iteration.
\end{itemize}

%%%%%TILING PLAN 1
\begin{table}[h]
$$
\begin{array}{|c|}\hline
\\
\mbox{TILING PLAN 1 :}  \quad \cB\rightarrow \cB^*, \quad s\in I_{0,-1}=(\beta+\beta^3,\alpha\beta)\\  
\\ \hline
\begin{array}{l}
\\
\mbox{Parent dressed domain:}\qquad \widehat{\cB}(1,s) \\ 
\\
\mbox{Tiling domain:} \quad {\rm B}^*= {\tt T}_{(-2 \beta^2+\beta s,-4 \beta- s)}\,\, {\tt R}_2 \, \,Q_1(2\beta^2-\beta s)\\ 
\end{array}
\\
\\
\mbox{\sc Data for Direct Iteration}\\
\begin{array}{c|c|c|c|c|c|c}\hline
\multicolumn{3}{c|}{\mbox{Source Polygon}}&\multicolumn{2}{|c|}{\mbox{Initial Placement}}&\multicolumn{2}{|c}{\mbox{Destination}}\\ \hline
\mbox{Tile}&Q_{\#}&\mbox{Parameters}&{\tt R}_{\#}&\mbox{Translation}&\mbox{Tile}&{\tt R}_{\#}\\ \hline\hline
  B^*_1&1&-2 \alpha\beta^3+\beta s&7&(0,-2-4\beta^3)& B^*_1\oslash B^*&6 \\ \hline
 B^*_2& 7&-4\beta^2+\alpha s&6&(0,14-12\alpha+\alpha s)& B^*_2\oslash B^* &5 \\ \hline
 B^*_3& 6&-\alpha\beta^3+\beta s,\alpha\beta-\ s &5&(0,12-10\alpha)& B^*_3\oslash B^* &6 \\ \hline
 B^*_4& 7&2\beta-\alpha s&1&(-2\beta+\alpha s,12-10\alpha)& B^*_4\oslash B^*&7 \\ \hline
 B^*_5& 1&\alpha\beta^2-\beta s &2&(-\alpha\beta^2+\beta s, 14-11\alpha- s)& B^*_5\oslash B^*&0 \\ \hline \hline
 \Pi_1& 2& s &0&(-  s,-  s)& \Pi_1 &2\\ \hline
 \Pi_2& 5&2\beta-\ s&0&(-2\beta-s, - s)& \Pi_2 &5\\ \hline
 \Pi_3&3&-\beta+ s&0&(-1- s,\beta-s)& \Pi_3 &4\\ \hline
 \Pi_4& 5&\alpha\beta- s&0&(-\alpha\beta- s,-\alpha\beta-s)& \Pi_4 &3\\ \hline
 \Pi_5& 5&\alpha\beta- s&0&(\alpha\beta^2- s,-\alpha- s)& \Pi_5 &3\\ \hline
 \Pi_6& 7&2\beta-\alpha s&5&(4\beta-(2+\alpha) s,-2\beta-2 s)& \Pi_6&0\\ \hline
 \Pi_7& 1&\alpha\beta^2-\beta s&6&(-6+5\alpha-\omega s,-\alpha- s)& \Pi_ 7&0\\ \hline
 \Pi_8& 5&-2\alpha\beta^2+ s &0&(2\alpha\beta^2- s,10-8\alpha- s)& \Pi_ 8&3\\ \hline
 \Pi_9&10&\alpha\beta- s,-2\alpha\beta^2+ s
&1&(\beta- s,-3\beta- s)&\Pi_9&4\\ 
\end{array}\\ \hline
\end{array}
$$
\end{table}
%

%TILING PLAN 2
\begin{table}[h]
$$
\begin{array}{|c|}\hline
\\
\mbox{TILING PLAN 2 :} \quad \cB\rightarrow \cB^*, \quad s =\beta\\  
\\ \hline
\begin{array}{l}
\\
\mbox{Parent dressed domain:}\qquad \widehat{\cB}(1, s) \\ 
\\
\mbox{Tiling domains:} \quad\begin{array}{l}
{\rm B}^*= {\tt T}_{(-\beta,1-2\alpha)}\,\, {\tt R}_2 \, \,Q_1( \beta )\\ 
{\rm B}^\dagger={\tt T}_{(-4\beta,-2\alpha\beta^2) }\,\, {\tt R}_7 \, \,{Q}_1(\beta^3 )
\end{array}
\end{array}
\\
\\
\mbox{\sc Data for Direct Iteration}\\
\begin{array}{c|c|c|c|c|c|c}\hline

\multicolumn{3}{c|}{\mbox{Source Polygon}}&\multicolumn{2}{|c|}{\mbox{Initial Placement}}&\multicolumn{2}{|c}{\mbox{Destination}}\\ \hline
\mbox{Tile}&Q_{\#}&\mbox{Parameters}&{\tt R}_{\#}&\mbox{Translation}&\mbox{Tile}&{\tt R}_{\#}\\ \hline\hline
B^*_1& 1 & \beta^2& 5 &  (0,-4\beta)& B^*_1\oslash B^*  &2\\ \hline
B^*_2& 7 &\alpha\beta & 3 &  (0,2-3\alpha) & B^*_2\oslash B^* &3\\ \hline
B^\dagger_1&1 &\beta^4& 2 & (-13+8\alpha,15-11\alpha)& B^\dagger_1\oslash B^\dagger &2\\ \hline
 B^\dagger_2&7 & \alpha\beta^3& 0 & (-3\alpha\beta,-2\alpha\beta^2) &B^\dagger_2\oslash B^\dagger &3\\ \hline \hline
 \Pi_1&2 & \beta & 0 &  \left(-\beta,-\beta \right) &\Pi_1&2\\ \hline
 \Pi_2&5 & \beta  & 0 &  (-3\beta,-\beta &\Pi_2&5\\ \hline
 \Pi_3&5 & \beta ^3& 0 & (11-9\alpha,-\beta)  &\Pi_3 &5\\ \hline
 \Pi_4&5 &\beta^2& 0 & (-1,-1) &\Pi_4&3
 \end{array}\\ \hline
\end{array}
$$
\end{table}
%

%%%%%%%%TILING PLAN 3
\begin{table}[h]
$$
\begin{array}{|c|}\hline
\\
\mbox{TILING PLAN 3 :} \quad \cB\rightarrow \cP, \quad s\in(0,\alpha\beta^2)\\  
\\ \hline
\begin{array}{l}
\\
\mbox{Parent dressed domain:}\qquad \widehat{\cB}(1, s) \\ 
\\
\mbox{Tiling domains:} \quad\begin{array}{l}
{\rm P}= {\tt T}_{(-2+\beta s,- s)}\,\,{Q}_6(  2 -\beta s , s  )\\ 
{\rm A}={\tt T}_{(-2\beta^2-\alpha\beta s,-2\alpha\beta-2 s)}\,\,\widehat{\rm A}(2 \beta^2-\alpha s) 
\end{array}

\end{array}
\\
\\
\mbox{\sc Data for Direct Iteration}\\
\begin{array}{c|c|c|c|c|c|c}\hline
\multicolumn{3}{c|}{\mbox{Source Polygon}}&\multicolumn{2}{|c|}{\mbox{Initial Placement}}&\multicolumn{2}{|c}{\mbox{Destination}}\\ \hline
\mbox{Tile}&Q_{\#}&\mbox{Parameters}&{\tt R}_{\#}&\mbox{Translation}&\mbox{Tile}&{\tt R}_{\#}\\ \hline\hline
 P_0&2 &  s& 0 &(-s,-  s) &P_0\oslash P&2\\ \hline
 P_1&9 & \beta- s,s  & 0 &(-\alpha- s,-s) &P_1\oslash P&4\\ \hline
 P_2&7 & \alpha s  & 1 &(-2\beta-\alpha \omega s,-2 s)&P_2\oslash P&7\\ \hline
 P_3& 1 &\beta s & 3 &(-2\beta-2 s,-2 s)&P_5\oslash P&6\\ \hline
 P_4&3 & 
  s & 3 &(-2\beta- s,- s)&P_4\oslash P&5\\ \hline
 P_5&4 & 
\alpha\beta^2-\beta s,\alpha\beta^2-\omega s & 4 &(-\alpha\beta- s,-\alpha\beta^2+\beta s)&P_3\oslash P&4\\ \hline
P_6& 3 & s& 0 &(-2\beta^2- s,- s) &P_6\oslash P&3\\ \hline
 P_{\rm in}&12 & 2\beta^2-\alpha\ s,2 s & 0 &(-2\beta^2-\alpha\beta s,-2 s)&A_2\oslash {\rm A}&0\\ \hline\hline
 A_1&1 & \alpha\beta-\omega s & 2 &(-\alpha\beta+\beta s,-\alpha- s)&A_1\oslash {\rm A}&5\\ \hline
A_2& 13 & 2\beta^2-\alpha s,2 s & 5 &(-2\beta^2-\alpha\beta s,-2\alpha\beta-2 s)&P_{\rm in}\oslash P&3\\ \hline
 A_3&1 &\alpha\beta^2-\omega s & 4 &(-\alpha\beta- s,-\alpha +\beta s)&A_3\oslash {\rm A}&1\\ \hline\hline
 \Pi&5 &\alpha\beta- s & 0 &(-\alpha\beta- s,-\alpha\beta- s)&\Pi&3
\end{array}\\ \hline
\end{array}
$$
\end{table}
%%%%%%%%%%%%%
%
%%%%%TILING PLAN 4
\begin{table}[h]
$$
\begin{array}{|c|}\hline
\\
\mbox{TILING PLAN 4 : }\quad \cP\rightarrow \cB, \quad s=\beta^2\\  
\\ \hline
\\
\mbox{Parent dressed domain:}\qquad \widehat{\cP}(\alpha,\beta^2) \\ 
\\
\mbox{Tiling domains:} \quad \begin{array}{l}
{\rm B}^*={\tt T}_{(-\alpha\beta,0)}{\tt R}_4 Q_1(\beta^2) \\ 
{\rm B}^\dagger={\tt T}_{(-\alpha\beta^2,-2\beta^3)} Q_1(\beta^5)
\end{array}
\\
\\
\mbox{\sc Data for Direct Iteration}\\
\begin{array}{c|c|c|c|c|c|c}\hline
\multicolumn{3}{c|}{\mbox{Source Polygon}}&\multicolumn{2}{|c|}{\mbox{Initial Placement}}&\multicolumn{2}{|c}{\mbox{Destination}}\\ \hline
\mbox{Tile}&Q_{\#}&\mbox{Parameters}&{\tt R}_{\#}&\mbox{Translation}&\mbox{Tile}&{\tt R}_{\#}\\ \hline\hline
 B^*_1&1 &\beta^3 & 1 & (-3\beta^2,\beta^2) &B^*_1\oslash B^*&6\\ \hline
B^*_2& 7 & \alpha\beta ^2 & 0 & (-5+3\alpha,\beta^2) &B^*_2\oslash B^*&5\\ \hline
B^\dagger_1&1 & \beta^5 & 5 &(-\beta^5-\alpha\beta^2,-\beta^2) &B^\dagger_1\oslash B^\dagger&6\\ \hline
B^\dagger_2&7 &\alpha\beta^4&4 &  (-3\beta^3,-\beta^2) &B^\dagger_2\oslash B^\dagger&5\\ \hline
 \Pi_1&2 & \beta^2 & 0 &(0,0)&\Pi_1&2\\ \hline
\Pi_2& 9 & \alpha\beta^2,\beta^2& 0 & (-\alpha,0))&\Pi_2&4\\ \hline
 \Pi_3&1 & \beta^3 & 3 &(-1,-\beta^2)&\Pi_3&0\\ \hline
 \Pi_4&5 & \beta^2 & 0 &(-2\beta,0)&\Pi_4&5\\ \hline
 \Pi_5&5 &\beta^2& 0 & (-2\beta^2,0)
 &\Pi_5&3\\\hline
 \Pi_6&5 & \beta^4& 0 & (\alpha\beta^4-\alpha\beta,0)   &\Pi_6&3
 \end{array}\\ \hline
 \end{array}
$$
\end{table}

%%%%%%%%%TILING PLAN 5
\begin{table}[h]
$$
\begin{array}{|c|}\hline
\\
\mbox{TILING PLAN 5 :} \quad \cB\rightarrow \cT_-, \quad s\in (\beta,\beta+\beta^4]\\  
\\ \hline
\begin{array}{l}
\\
\mbox{Parent dressed domain:}\qquad \widehat{\cB}(1,s) \\ 
\\
\mbox{Tiling domains:} \quad \begin{array}{l} 
\mathrm{T}_- = {\tt T}_{(-\alpha\beta+\beta s,-\alpha-s)}{\tt R}_2 Q_1(\alpha\beta-\beta s)\\ 
\mathrm{A} = {\tt T}_{(-2\beta-\alpha\beta s,-2 s)}\widehat{A}(2\beta-\alpha s)\
\end{array}
\end{array}
\\
\\
\mbox{\sc Data for Direct Iteration}\\
\begin{array}{c|c|c|c|c|c|c}\hline
\multicolumn{3}{c|}{\mbox{Source Polygon}}&\multicolumn{2}{|c|}{\mbox{Initial Placement}}&\multicolumn{2}{|c}{\mbox{Destination}}\\ \hline
\mbox{Tile}&Q_{\#}&\mbox{Parameters}&{\tt R}_{\#}&\mbox{Translation}&\mbox{Tile}&{\tt R}_{\#}\\ \hline\hline 
T^-_{1}&1&\beta  s&5& (0,-2\beta-2 s) &T^-_{1}\oslash{\rm T}_-&2\\ \hline
T^-_{2}&7&2-\omega s&3& (2\beta-2 s,-2\alpha+\alpha s)  &T^-_{2}\oslash{\rm T}_-&3\\ \hline
T^-_{3}&3&-\beta+s & 1&(\beta-s,7-6\alpha-s) &T^-_{3}\oslash{\rm T}_-&1\\ \hline
T^-_{4}&4&2\beta^2-\beta s, 2\alpha\beta-\omega s&2& (-2\beta^2+\beta s,-4 \beta-s) &T^-_{4}\oslash{\rm T}_-&0\\ \hline
T^-_{{\rm in}}& 13&
2\beta-\alpha s,-2\beta+2 s& 2&(2\beta-2 s,-4\beta-\alpha\beta s)&A_2\oslash\mathrm{A}&3\\ \hline
A_1& 1&\alpha-\omega s&7& (-2 s ,-2 s)&A_1\oslash\mathrm{A}&3\\ \hline
A_2&12&2\beta-\alpha s,-2\beta+2 s& 0&(-2\beta-\alpha\beta s,-2 s)  &T^-_{{\rm in}}\oslash\mathrm{T}_-&6\\ \hline
A_3& 1&2\alpha\beta-\omega s &5& (-2 s  ,-2\beta) &\widehat{A}_3\oslash\mathrm{A}&7\\ \hline
\Pi_1&2&s&0& (-s  ,-s) &\Pi_1&2\\ \hline
 \Pi_2&5&2\beta-s &0& (-2 \beta-s,-s  ) &\Pi_2&5\\ \hline
 \Pi_3&3&-\beta+s&0& (-1-s,\beta-s) &\Pi_3&4
\end{array}\\ \hline
\end{array}
$$
\end{table}
%%%%%%%%%%%%%

%%%%%%%%%TILING PLAN 6
\begin{table}[h]
$$
\begin{array}{|c|}\hline
\\
\mbox{TILING PLAN 6 :} \quad \cB\rightarrow \cT_+, \quad s\in [\beta-\beta^6,\beta)\\  
\\ \hline
\begin{array}{l}
\\
\mbox{Parent dressed domain:}\qquad \widehat{\cB}(1,s) \\ 
\\
\mbox{Tiling domains:} \quad \begin{array}{l} 
\mathrm{T}_+={\tt T}_{(-\alpha\beta+\beta s,-\alpha-s)}{\tt R}_2 Q_1(\alpha\beta-\beta s)\\ 
\mathrm{A}={\tt T}_{(-2\beta^2-\alpha s,-2 s)}{\tt R}_3\widehat{A}(-2\alpha\beta^2+\alpha s)\
\end{array}
\end{array}
\\
\\
\mbox{\sc Data for Direct Iteration}\\
\begin{array}{c|c|c|c|c|c|c}\hline
\multicolumn{3}{c|}{\mbox{Source Polygon}}&\multicolumn{2}{|c|}{\mbox{Initial Placement}}&\multicolumn{2}{|c}{\mbox{Destination}}\\ \hline
\mbox{Tile}&Q_{\#}&\mbox{Parameters}&{\tt R}_{\#}&\mbox{Translation}&\mbox{Tile}&{\tt R}_{\#}\\ \hline\hline
 T^+_{1}&1 & \beta s & 5 & (0,-2\beta-2 s) &T^+_{1}\oslash \mathrm{T}_+&2\\ \hline
  T^+_{2}&7 &\alpha s  & 3 & (0,-2\beta-\alpha\omega s) &T^+_{2}\oslash \mathrm{T}_+&3\\ \hline
  T^+_{3}&4 & \beta s,-2\beta+\omega s & 7 & (0,-2) &T^+_{3}\oslash \mathrm{T}_+&6\\ \hline
  T^+_{4}&3 &\beta-s  & 7 & \beta^2-s,-3+\alpha-s)&T^+_{4}\oslash \mathrm{T}_+&1\\ \hline
  T^+_{5}&2 & \beta-s & 1 & (\alpha\beta^3-s,-10+6\alpha-s) &T^+_{5}\oslash \mathrm{T}_+&4\\ \hline
   T^+_{6}&10&\beta-s,-2\beta^2+s & 7 & (\alpha\beta^3/2-s,-5\alpha\beta/2-s ) &T^+_{6}\oslash {\rm T}_+&6\\ \hline
  T^+_{7}&2 & \beta-s & 1 & (-s,-\alpha-s) &T^+_{7}\oslash \mathrm{T}_+&4\\ \hline
  T^+_{\rm in}&12 & -2\alpha\beta^2+\alpha s,2\beta-2s & 3 & (2\beta^2-\alpha s,-2\alpha\beta-2 s) &A_2\oslash \mathrm{A}&0\\ \hline
 A_1&1 & -2\beta+\omega s & 5 & (-2\beta,-2s ) &A_1\oslash {\rm A}&5\\ \hline
 A_2 &13 &  -2\alpha\beta^2+\alpha s,2\beta-2s  & 0 & (-2\beta^2-\alpha s,-2 s ) &T^+_{\rm in}\oslash \mathrm{T}_+&5\\ \hline
A_3&1 & \beta^3-1+\omega s & 7 & (-2\beta,-2\beta) &A_3\oslash {\rm A}&1\\ \hline
 \Pi_1&2 & s   & 0 & (-s  ,-s  ) &\Pi_1&2\\ \hline
\Pi_2&5 & s  & 0 & (-2\beta-s,-s ) &\Pi_2&5\\ \hline
 \Pi_3&5 & \alpha\beta-s& 0 & (-\alpha\beta-s,-\alpha\beta-s) &\Pi_3&3
\end{array}\\ \hline
\end{array}
$$
\end{table}

%%%%%%%%%TILING PLAN 7
\begin{table}[h]
$$
\begin{array}{|c|}\hline
\\
\mbox{TILING PLAN 7 :} \quad \cP^*\rightarrow \cT_-, \quad s\in (\beta^2,\beta^2+\beta^6]\\  
\\ \hline
\begin{array}{l}
\\
\mbox{Parent dressed domain:}\qquad \widehat{\cP}(\alpha,s) \\ 
\\
\mbox{Tiling domains:} \quad \begin{array}{l} 
\mathrm{T}_- = {\tt T}_{(-\alpha\beta,-\alpha\beta^2+\alpha s)}{\tt R}_4Q_1(\alpha\beta^2-\beta s)\\ 
\mathrm{A} = {\tt T}_{(-2\beta^4- s,-2\alpha\beta^4-\beta s)}{\tt R}_6\widehat{A}(\alpha\beta^5+\alpha\beta^2-\alpha s)
\end{array}
\end{array}
\\
\\
\mbox{\sc Data for Direct Iteration}\\\begin{array}{c|c|c|c|c|c|c}\hline
\multicolumn{3}{c|}{\mbox{Source Polygon}}&\multicolumn{2}{|c|}{\mbox{Initial Placement}}&\multicolumn{2}{|c}{\mbox{Destination}}\\ \hline
\mbox{Tile}&Q_{\#}&\mbox{Parameters}&{\tt R}_{\#}&\mbox{Translation}&\mbox{Tile}&{\tt R}_{\#}\\ \hline\hline
T^-_1&1&\beta s&1&(-2\beta^2- s,s)&T^-_1\oslash {\rm T}_-&6\\ \hline
T^-_2&7&2\beta-\alpha\omega s&0&(-2\alpha\beta+\omega s,2\beta^2- s)&T^-_2\oslash {\rm T}_-&5\\ \hline
T^-_3&4&2\beta^3-\beta s&4&(-3\beta^2,-2\beta^3+\alpha s)&T^-_3\oslash {\rm T}_-&0\\ \hline
T^-_4&3&-\beta^2+s&0&(-\beta^4-\alpha\beta,\beta^2)&T^-_4\oslash {\rm T}_-&7\\ \hline
T^-_5&4&\begin{array}{c}\alpha\beta^3-\beta s,\\ \beta^4+\beta-\omega s\end{array}&0&(-\alpha\beta,\beta^4+\beta-\alpha s)&T^-_5\oslash {\rm T}_-&0\\ \hline
T^-_6&3&-\beta^2+s&7&(\beta^4+\alpha\beta,\beta^2)&T^-_6\oslash {\rm T}_-&1\\ \hline
T^-_7&4&\begin{array}{c}2\alpha\beta^4-\beta s,\\ \beta^5+\beta-\omega s\end{array}&4&(\alpha\beta^4+\beta,-2\alpha\beta^4+\alpha s)&T^-_7\oslash {\rm T}_-&0\\ \hline
T^-_8&3&-\beta^2+s&0&(-4\beta^4-\beta,\beta^2)&T^-_8\oslash {\rm T}_-&7\\ \hline
T^-_{\rm in}&12&\begin{array}{c}\alpha\beta^5+\alpha\beta^2-\alpha s,\\ -2\beta^2+2 s\end{array}&0&(-\beta^5+\beta+\beta s,2\beta^2-s)&A_2\oslash {\rm A}&6\\ \hline
A_1&1&\beta^4+\beta-\omega s&0&(-\alpha\beta^2,\alpha\beta^3-\alpha s)&A_1\oslash {\rm A}&5\\ \hline
A_2&13&\begin{array}{c}\alpha\beta^5+\alpha\beta^2-\alpha s,\\ -2\beta^2+2 s\end{array}&3&(-2\beta^4- s,-2\alpha\beta^4-\beta s)&T^-_{\rm in}\oslash {\rm T}_-&1\\ \hline
A_3&1&\beta^5+\beta-\omega s&2&(-2\alpha\beta^2+\alpha s,-2\beta^3)&A_3\oslash {\rm A}&1\\ \hline
\Pi_1&2&s&0&(0,0)&\Pi_1&2\\ \hline
\Pi_2&3&-\beta^2+s&7&(-3\beta,\beta^2)&\Pi_2&4\\ \hline
\Pi_3&9&\beta-s,s&0&(-\alpha,0)&\Pi_3&4\\ \hline
\Pi_4&1&\beta s&3&(-2\beta- s,-s)&\Pi_4&0\\ \hline
\Pi_5&5&2\beta^2-s&0&(-2\beta^2,0)&\Pi_5&3\\ \hline
\Pi_6&5&s&0&(-2\beta-s,-s)&\Pi_6&5\\ \hline
\Pi_7&5&2\alpha\beta^3-s&0&(-2\alpha\beta^3,-2\beta^3)&\Pi_7&3\\ \hline
\end{array}\\ \hline
\end{array}
$$
\end{table}

%%%%%%%%%TILING PLAN 8
\begin{table}[h]
$$
\begin{array}{|c|}\hline
\\
\mbox{TILING PLAN 8 :}  \quad \cP^*\rightarrow \cT_+, \quad s\in [\beta^2-\beta^7,\beta^2)\\  
\\ \hline
\begin{array}{l}
\\
\mbox{Parent dressed domain:}\qquad \widehat{\cP}(\alpha,s) \\ 
\\
\mbox{Tiling domains:} \quad \begin{array}{l} 
{\rm T}_+={\tt T}_{(-\alpha\beta,-\alpha\beta^2+\alpha s)}{\tt R}_4 Q_1(\alpha\beta^2-\beta s)\\ 
{\rm A} ={\tt T}_{(-4\alpha\beta^2+s,\alpha\beta^5-\alpha\beta^2+\beta s)}{\tt R}_6\widehat{A}(\alpha\beta^5-\alpha\beta^2+\alpha s)
\end{array}
\end{array}
\\
\\
\mbox{\sc Data for Direct Iteration}\\
\begin{array}{c|c|c|c|c|c|c}\hline
\multicolumn{3}{c|}{\mbox{Source Polygon}}&\multicolumn{2}{|c|}{\mbox{Initial Placement}}&\multicolumn{2}{|c}{\mbox{Destination}}\\ \hline
\mbox{Tile}&Q_{\#}&\mbox{Parameters}&{\tt R}_{\#}&\mbox{Translation}&\mbox{Tile}&{\tt R}_{\#}\\ \hline\hline
T^+_1&1&\beta s&1&(-2\beta^2- s,s)&T^+_1\oslash {\rm T}_+&6\\ \hline
T^+_2&7&\alpha s&0&(-2\beta^2- \omega s,s)&T^+_2\oslash {\rm T}_+&5\\ \hline
T^+_3&4&\beta s,-2\beta^2+s&7&(-2\beta+s,s)&T^+_3\oslash {\rm T}_+&2\\ \hline
T^+_4&3&\beta^2- s&2&(2\beta^2-1,\beta^3)&T^+_4\oslash {\rm T}_+&7\\ \hline
T^+_5&2&\beta^2- s&1&(-3\beta^3-\beta,\alpha\beta^4)&T^+_5\oslash {\rm T}_+&4\\
T^+_6&10&-2\beta^3+ s,\beta^2- s&1&(-5\alpha\beta^2/2,\alpha\beta^4/2)&T^+_6\oslash {\rm T}_+&2\\ \hline
T^+_7&2&\beta^2- s&1&(-\alpha\beta,0)&T^+_7\oslash {\rm T}_+&4\\ \hline
T^+_8&4&\begin{array}{c}-2\beta^4+\beta s,\\ -4\alpha\beta^3+s\end{array}&3&\begin{array}{c}(-2\alpha\beta^2- s,\\ 9\beta^3-\beta- s)\end{array}&T^+_8\oslash {\rm T}_+&2\\ \hline
T^+_9&3&\beta^2- s&7&(-9\beta^3,9\beta^3-\alpha\beta)&T^+_9\oslash {\rm T}_+&1\\ \hline
T^+_{\rm in}&12&\begin{array}{c}-9\beta^3+\beta+\alpha s,\\ 2\beta^2-2 s\end{array}&7&\begin{array}{c}(-10\beta^3+\beta s,\\-4\beta^4+s)\end{array}&A_2\oslash {\rm A}&7\\ \hline
A_1&1&\beta^4-\beta+s&0&(-\alpha\beta^2,\beta^4-\beta+\alpha s)&A_1\oslash {\rm A}&5\\ \hline
A_2&13&\begin{array}{c}-9\beta^3+\beta+\alpha s,\\2\beta^2-2 s\end{array}&3&\begin{array}{c}(-4\alpha\beta^3+s,\\ -9\beta^4+\beta+\beta s)\end{array}&T^+_{\rm in}\oslash {\rm T}_+&2\\ \hline
A_3&1&-4\alpha\beta^3+\omega s&2&(-\alpha s,-2\beta^3)&A_3\oslash {\rm A}&1\\ \hline
\Pi_1&2& s&0&(0,0)&\Pi_1&2\\ \hline
\Pi_2&9&\beta- s,\beta s&0&(-\alpha,0)&\Pi_2&4\\ \hline
\Pi_3&1&\beta s&3&(-2\beta- s,- s)&\Pi_3&0\\ \hline
\Pi_4&5& s&0&(-2\beta^2,0)&\Pi_4&3\\ \hline
\Pi_5&5& s&0&(-2\beta,0)&\Pi_5&5\\ \hline
\Pi_6&5&-2\beta^3+ s&0&(-2\alpha\beta^3,-2\beta^3)&\Pi_6&3
\end{array}\\ \hline
\end{array}
$$
\end{table}

%%%%%%%%%%%%%TILING PLAN 9
\begin{table}[h]
$$
\begin{array}{|c|}\hline
\\
\mbox{TILING PLAN  9 :}  \quad \cT_-\rightarrow \cD_-, \quad s-\beta\in(0,\beta^4]\\  
\\ \hline
\begin{array}{l}
\\
\mbox{Parent dressed domain:}\qquad \widehat{\cT}_-(\beta,h),\quad h=s-\beta \\ 
\\
\mbox{Tiling domains:} \quad \begin{array}{l} 
{\rm D}={\tt T}_{(-\beta^2+\alpha h,-\beta+ h)}{\tt R}_2 Q_4(\beta^2-\beta h,\beta^2-\omega h)\\ 
{\rm X}={\tt T}_{(h,-\alpha\beta^3+\alpha h)}{\tt R}_ 2Q_{12}(\alpha\beta^3-\alpha h,2h)\\
{\rm A}={\tt T}_{(-\alpha\beta^3+\beta h,2 h)}{\tt R}_3\widehat{A}(\alpha\beta^3-\alpha h)
\end{array}
\end{array}
\\
\\
\mbox{\sc Data for Direct Iteration}\\
\begin{array}{c|c|c|c|c|c|c}\hline
\multicolumn{3}{c|}{\mbox{Source Polygon}}&\multicolumn{2}{|c|}{\mbox{Initial Placement}}&\multicolumn{2}{|c}{\mbox{Destination}}\\ \hline
\mbox{Tile}&Q_{\#}&\mbox{Parameters}&{\tt R}_{\#}&\mbox{Translation}&\mbox{Tile}&{\tt R}_{\#}\\ \hline\hline
D_1''&4&\beta^3-\beta h,\beta^3-\omega h&2&(-\beta^3+\alpha h,11-8 \alpha + h)&D_1''\oslash {\rm D}&0\\ \hline
D_2''&3&h&2&(0,18-13 \alpha+h)&D_2''\oslash {\rm D}&7\\ \hline
D_{\rm in}''&12&\alpha\beta^3-\alpha h,2 h&6&(-h,18-13 \alpha+\alpha\beta h)&A_2\oslash {\rm A}&5\\ \hline
D_0&2&h&0&(0,-2\alpha\beta^2+h)&D_0\oslash {\rm D}&2\\ \hline
D_1'&12&\alpha\beta^3-\alpha h,2 h&2&(\alpha\beta^3-\alpha h,2 h)&X&0\\ \hline
A_1&1&\beta^2-\omega h&5&(-h,2 h)&A_1\oslash {\rm A}&5\\ \hline
A_2&13&\alpha\beta^3-\alpha h,2 h&0&(-\alpha\beta^3+\beta h)&D_{\rm in}''\oslash {\rm D}&2\\ \hline
A_3&1&\beta^3-\omega h&7&(-h,0)&A_3\oslash {\rm A}&1\\ \hline
\Pi_1&5&\beta^2-h&0&(-\beta^2,-\beta^2+h)&\Pi_1&3\\ \hline
\Pi_2&3&h&1&(0,-\alpha\beta^2+h)&\Pi_2&0\\ \hline
\Pi_3&4&\beta^3-\beta h,\beta^3-\omega h&6&(\beta^3-\alpha h,-\beta^2+h)&\Pi_3&0\\ \hline
\end{array}\\ \hline
\end{array}
$$
\end{table}
%%%%%%%%%%TILING PLAN 10
\begin{table}[h]
$$
\begin{array}{|c|}\hline
\\
\mbox{TILING PLAN 10 :}  \quad \cT_+\rightarrow \cD_+, \quad \beta-s\in(0,\beta^6]\\  
\\ \hline
\begin{array}{l}
\\
\mbox{Parent dressed domain:}\qquad \widehat{\cT}_+(\beta^3,h),\quad h=\beta-s \\ 
\\
\mbox{Tiling domains:} \quad \begin{array}{l} 
{\rm D}={\tt T}_{(\beta^4+ h, -\beta^3+h)}\,\mathtt{R}_3 \,Q_4(\beta^4-\beta h,\beta^4-\omega h)\\
{\rm A}={\tt T}_{(13-9\alpha+\beta h, 13-9\alpha-h)}\,\widehat{\rm A}(\alpha\beta^4-\alpha h)\\
{\rm X}={\tt T}_{(\beta^5+\beta h, -\beta^5+h)}\,\mathtt{R}_3 \,Q_{12}(\alpha\beta^5-\alpha h, 2 h)\\
{\rm Y}={\tt T}_{(\beta^5+\beta h, 3 \beta^3-h)}\,\mathtt{R}_5\, Q_{13}(\alpha\beta^5-\alpha h, 2 h)
\end{array}
\end{array}
\\
\\
\mbox{\sc Data for Direct Iteration}\\
\begin{array}{c|c|c|c|c|c|c}\hline
\multicolumn{3}{c|}{\mbox{Source Polygon}}&\multicolumn{2}{|c|}{\mbox{Initial Placement}}&\multicolumn{2}{|c}{\mbox{Destination}}\\ \hline
\mbox{Tile}&Q_{\#}&\mbox{Parameters}&{\tt R}_{\#}&\mbox{Translation}&\mbox{Tile}&{\tt R}_{\#}\\ \hline\hline
D_{\rm in}''&13&\alpha\beta^4-\alpha h,2h&7&(\beta^4+h,-\beta^4+\beta h)&A_2\oslash {\rm A}&6\\ \hline
D_{\rm in}'&12&-\alpha\beta^5-\alpha h,2 h&3&(\beta^3+\beta h,-\beta^3+h)&X&0\\ \hline
D_0&2&h&1&(2 \beta^4,-2 \beta^4)&D_0 \oslash {\rm D}&6\\ \hline
{\rm Y}&13&\alpha\beta^5-\alpha h,2h&5&(\beta^5+\beta h,3\beta^3-h)&D_{\rm in}'\oslash {\rm D}&2\\ \hline
A_1&1&\beta^3-\omega h&7&(-11+8\alpha-h,13-9\alpha-h)&A_1\oslash {\rm A}&3\\ \hline
A_2&12&\alpha\beta^4-\alpha h,2h&0&(13-9\alpha+\beta h,13-9\alpha-h)&D_{\rm in}''\oslash {\rm D}&3\\ \hline
A_3&1&\beta^4-\omega h&5&(-11+8\alpha-h,13-9\alpha+h)&A_3\oslash {\rm A}&7\\ \hline
\Pi_1&5&\beta^2-h&0&(2\beta^3,\alpha\beta^3)&\Pi_1&3\\ \hline
\Pi_2&5&\beta^3-h&0&(\alpha\beta^2,-2\beta^3)&\Pi_2&5\\ \hline
\Pi_3&2&h&1&(-\alpha\beta^3,\alpha\beta^3)&\Pi_3&4\\ \hline
\Pi_4&4&\beta^2-\beta h,\beta^2-\omega h&7&(-11+8\alpha-h,-\beta^3-h)&\Pi_4&0\\ \hline
\Pi_5&3&h&3&(\beta^3,-\beta^3)&\Pi_5&4\\ \hline
\Pi_6&3&h&2&(\beta^4,-\beta^4)&\Pi_6&0\\ \hline
\Pi_7&4&\beta^5-\beta h,\beta^5-\omega h&7&(\beta^4-h,-\beta^5-h)&\Pi_7&0\\ \hline
\Pi_8&3&h&3&(-\beta^3,\beta^3)&\Pi_8&0\\ \hline
\Pi_9&4&\beta^4-\beta h,\beta^4-\omega h&7&(-\beta^4-h,\beta^3-h)&\Pi_9&0\\
\end{array}\\ \hline
\end{array}
$$
\end{table}
 
  %%%%%%%%%TILING PLAN 11
\begin{table}[h]
$$
\begin{array}{|c|}\hline
\\
\mbox{TILING PLAN 11 :}  \quad \cD^*_-\rightarrow \cB^*, \quad s-\beta\in [\beta^5,\beta^4)\\  
\\ \hline
\begin{array}{l}
\\
\mbox{Parent dressed domain:}\qquad \widehat{\cD}_-(\beta^3,h),\quad h=s-\beta \\ 
\\
\mbox{Tiling domain:} \quad {\rm B}^*={\tt T}_{(\beta^3,\beta^4-\alpha h)} Q_1(\beta^4-\beta h)
\end{array}
\\
\\
\mbox{\sc Data for Direct Iteration}\\\begin{array}{c|c|c|c|c|c|c}\hline
\multicolumn{3}{c|}{\mbox{Source Polygon}}&\multicolumn{2}{|c|}{\mbox{Initial Placement}}&\multicolumn{2}{|c}{\mbox{Destination}}\\ \hline
\mbox{Tile}&Q_{\#}&\mbox{Parameters}&{\tt R}_{\#}&\mbox{Translation}&\mbox{Tile}&{\tt R}_{\#}\\ \hline\hline

B^*_1&1 &-\beta^6+\beta h&3&(-130+92\alpha-h,-h)&B^*_1\oslash {\rm B}^*&2\\ \hline
B^*_2& 7 &-\alpha\beta^5+\alpha h&1&-72+51\alpha-\omega h,-h)&B^*_2\oslash {\rm B}^*&3\\ \hline
B^*_3& 6 &\beta^5+\beta h,\beta^4-h&1&(2\alpha\beta^4-h,-h)&B^*_3\oslash {\rm B}^*&2\\ \hline
B^*_4& 7 &\alpha\beta^4-\alpha h&6&(2\alpha\beta^4-h,\alpha\beta^4-\omega h)&B^*_4\oslash {\rm B}^*&1\\ \hline
B^*_5& 1 &\beta^5-\beta h&0&(-65+46 \alpha,\beta^5-\alpha h)&B^*_5\oslash {\rm B}^*&0\\ \hline
\Pi_1& 1 &\beta^5-\beta h&4&\beta^4,-\beta^5-\alpha h)&\Pi_1&0\\ \hline
\Pi_2 &7 &\alpha\beta^4-\alpha h&2&h,-\alpha\beta^4+\omega h)&\Pi_2&0\\ \hline
\Pi_3& 1 &\beta^4-\beta h&4&-\beta^3,-\beta^4+\alpha h)&\Pi_3&0\\ \hline
\Pi_4&5 &-\beta^5+h&0&(\beta^5,-\beta^5)&\Pi_4&5\\ \hline
\Pi_5&10 &\beta^4-h,-\beta^5+h&1&(0,0)&\Pi_5&4\\ \hline
\Pi_6& 1 &-\beta^4+\omega h&7&(-2\beta^3+h,h)&\Pi_6&0\\ \hline
\Pi_7&8 &\alpha\beta^4-\alpha\beta h,2 h&5&(-\alpha\beta^3+\beta h,h)&\Pi_7&0\\ \hline
\end{array}\\ \hline
\end{array}
$$
\end{table}
%%%%%%%%%%%%%
  %%%%%%%%%TILING PLAN 12
\begin{table}[h]
$$
\begin{array}{|c|}\hline
\\
\mbox{TILING PLAN 12 :}  \quad \cD_+\rightarrow \cB^*, \quad s-\beta\in [\beta^6,\beta^5]\\  
\\ \hline
\begin{array}{l}
\\
\mbox{Parent dressed domain:}\qquad \widehat{\cD}(\beta^4,h),\quad h=s-\beta\\ 
\\
\mbox{Tiling domain:} \quad  {\rm B}^*={\tt T}_{(\beta^4,\beta^5-\alpha h)} Q_1(\beta^5-\beta h)
\end{array}
\\
\\
\mbox{\sc Data for Direct Iteration}\\\begin{array}{c|c|c|c|c|c|c}\hline
\multicolumn{3}{c|}{\mbox{Source Polygon}}&\multicolumn{2}{|c|}{\mbox{Initial Placement}}&\multicolumn{2}{|c}{\mbox{Destination}}\\ \hline
\mbox{Tile}&Q_{\#}&\mbox{Parameters}&{\tt R}_{\#}&\mbox{Translation}&\mbox{Tile}&{\tt R}_{\#}\\ \hline\hline
 B^*_1&1 &-\beta^7+\beta h&3&(314-222\alpha-h,-h)& B^*_1\oslash {\rm B}^*&2\\ \hline
B^*_2& 7 &-\alpha\beta^6+\alpha h&1&(174-123\alpha-\omega h,-h)&B^*_2\oslash {\rm B}^*&3\\ \hline
B^*_3& 6 &\beta^6+\beta h,\beta^5-h&1&2\alpha\beta^5-h,-h)&B^*_3\oslash {\rm B}^*&2\\ \hline
B^*_4& 7 &\alpha\beta^5-\alpha h&6&(2\alpha\beta^5-h,\alpha\beta^5-\omega h)&B^*_4\oslash {\rm B}^*&1\\ \hline
B^*_5& 1 &\beta^6-\beta h&0&(157-111\alpha,\beta^6-\alpha h)&B^*_5\oslash {\rm B}^*&0\\ \hline
\Pi_1& 1 &\beta^5-\beta h&4&-\beta^4,-\beta^5+\alpha h&\Pi_1&0\\ \hline
\Pi_2 &8 &\alpha\beta^5-\alpha\beta h,2h&6&-\alpha\beta^5-\beta h,h)&\Pi_2&0\\ \hline
\Pi_3& 7 &\alpha\beta^5-\alpha h&2&(h,-\alpha\beta^5+\omega h)&\Pi_3&0\\ \hline
\Pi_4&1 &\beta^6-\beta h&4&(\beta^5,-\beta^6+\alpha h)&\Pi_4&0\\ \hline
\Pi_5&5 &-\beta^6+h&0&(-\beta^6,\beta^6)&\Pi_5&5\\ \hline
\Pi_6& 6 &-\beta^6+\omega h,\beta^5-h&5&(h,h)&\Pi_6&0\\ \hline
\end{array}\\ \hline
\end{array}
$$
\end{table}
%%%%%%%%%%%%%

%%%%%%%%%%%%%%%%%%%%%%%%%%

\clearpage

\section*{Appendix C: Incidence matrices $\mathrm{M}(i,j)$ for scenario I}
$$
\begin{array}{c}{\rm M}(-1,3)=\\\left(\begin{array}{ccccc}271&121&157&71&71\\650&260&326&251&374\\370&148&187&134&197\\210&93&120&54&54\\36&18&24&12&12\end{array}\right)\end{array}\;
\begin{array}{c}{\rm M}(1,0)=\\\left(\begin{array}{ccccc}0&0&0&0&0\\6&3&0&2&2\\29&29&11&31&43\\18&18&8&19&28\\7&4&0&4&4\end{array}\right)\end{array}\;
\begin{array}{c}{\rm M}(-1,1)=\\\left(\begin{array}{ccccc}23&11&15&7&7\\90&36&46&31&46\\50&20&27&16&25\\18&9&12&6&6\\0&0&0&0&0\end{array}\right)\end{array}
$$

$$
\begin{array}{c}{\rm M}(1,-1)=\\\left(\begin{array}{ccccc}0&0&0&0&0\\18&9&12&6&6\\145&58&76&48&72\\90&36&46&31&46\\23&11&15&7&7\end{array}\right)\end{array}\;
\begin{array}{c}{\rm M}(-1,2)=\\\left(\begin{array}{ccccc}110&44&14&36&57\\400&169&64&66&66\\224&95&35&39&39\\90&36&12&27&42\\0&0&0&6&12\end{array}\right)\end{array}\;
\begin{array}{c}{\rm M}(1,-2)=\\\left(\begin{array}{ccccc}0&0&0&6&12\\90&36&12&27&42\\644&272&101&108&108\\400&169&64&66&66\\110&44&14&36&57\end{array}\right)\end{array}
$$
$$
\begin{array}{c}{\rm M}(-1,0)=\\\left(\begin{array}{ccccc}7&4&0&4&4\\18&18&8&19&28\\10&10&3&10&13\\6&3&0&2&2\\0&0&0&0&0\end{array}\right)\end{array}\;
\begin{array}{c}{\rm M}(1,-3)=\\\left(\begin{array}{ccccc}36&18&24&12&12\\210&93&120&54&54\\1045&418&526&390&576\\650&260&326&251&374\\271&121&157&71&71\end{array}\right)\end{array}\;
\begin{array}{c}{\rm M}(0,1)=\\\left(\begin{array}{ccccc}5&2&2&2&2\\10&4&6&3&6\\26&14&18&10&10\\18&9&12&6&6\\10&4&5&4&7\end{array}\right)\end{array}\;
$$
$$
\begin{array}{c}{\rm M}^*(-1,\infty)=\\ \left(\begin{array}{cc} 1 & 1 \\ 0 & 0 \\ 0 & 0 \\ 0 & 1 \\ 1 & 0 \end{array}\right)\end{array}\;
\begin{array}{c}{\rm M}^\dagger(-1,\infty)=\\ \left(\begin{array}{cc} 0 & 0 \\28 & 13 \\ 15 & 6 \\ 0 & 0 \\ 0 & 0 \\\end{array}\right)\end{array}\;
\begin{array}{c}{\rm M}^*(1,-\infty)=\\ \left(\begin{array}{cc} 1 & 0 \\ 0 & 1 \\ 0 & 0 \\ 0 & 0 \\ 1 & 1 \end{array}\right)\end{array}\;
\begin{array}{c}{\rm M}^\dagger(1,-\infty)=\\ \left(\begin{array}{cc} 0 & 0 \\ 0 & 0 \\ 43 & 19 \\ 28 & 13 \\ 0 & 0 \\\end{array}\right)\end{array}\;
$$
$$
\begin{array}{c}{\rm M}(0,-2)=\\\left(\begin{array}{ccccc}48&21&9&9&9\\90&36&12&27&42\\70&28&10&18&27\\48&21&8&10&10\\20&8&3&2&2\end{array}\right)\end{array}\;
\begin{array}{c}{\rm M}(0,2)=\\\left(\begin{array}{ccccc}20&8&3&2&2\\48&21&8&10&10\\140&56&20&66&42\\90&36&12&42&27\\48&21&9&9&9\end{array}\right)\end{array}\;
\begin{array}{c}{\rm M}(0,-1)=\\\left(\begin{array}{ccccc}10&4&5&4&7\\18&9&12&6&6\\13&7&9&5&5\\10&4&6&3&6\\5&2&2&2&2\end{array}\right)\end{array}\;
$$
$$
\begin{array}{c}{\rm M}(0,0)_A=\\\left(\begin{array}{ccccc}2&2&1&0&0\\6&3&0&2&2\\5&2&0&4&10\\2&2&0&3&6\\1&1&1&1&1\end{array}\right)\end{array}\;
\begin{array}{c}{\rm M}(0,0)_B=\\\left(\begin{array}{ccccc}1&1&1&1&1\\2&2&0&3&6\\10&4&0&2&5\\6&3&0&2&2\\2&2&1&0&0\end{array}\right)\end{array}\;
\begin{array}{c}{\rm M}(-\infty,0)=\\ \left(\begin{array}{cc}1&1\\4&1\end{array}\right)\end{array}\;
\begin{array}{c}{\rm M}(+\infty,0)=\\ \left(\begin{array}{cc}1&1\\0&1\\1&0\end{array}\right)\end{array}\;
$$

\clearpage
\section*{Appendix D: Incidence matrices $\mathrm{M}(i,j)$ for scenario II}
$$
\begin{array}{c}{\rm M}(-2,0)=\\
\left(
\begin{array}{ccccc}
 0 & 0 & 0 & 0 & 0 \\
 0 & 0 & 0 & 0 & 0 \\
 6 & 3 & 0 & 2 & 2 \\
 0 & 0 & 0 & 0 & 0 \\
 12 & 6 & 0 & 4 & 4 \\
 7 & 4 & 0 & 4 & 4 \\
 6 & 6 & 2 & 4 & 4 \\
 1 & 1 & 1 & 1 & 1 \\
 2 & 2 & 0 & 3 & 6 \\
\end{array}
\right) \end{array}\;
\begin{array}{c}{\rm M}(-2,1)=\\ \left(
\begin{array}{ccccc}
 0 & 0 & 0 & 0 & 0 \\
 0 & 0 & 0 & 0 & 0 \\
 6 & 6 & 12 & 9 & 18 \\
 0 & 0 & 0 & 0 & 0 \\
 12 & 12 & 24 & 18 & 36 \\
 7 & 7 & 15 & 11 & 23 \\
 16 & 10 & 16 & 12 & 30 \\
 2 & 2 & 2 & 2 & 5 \\
 6 & 3 & 6 & 4 & 10 \\
\end{array}
\right) \end{array}\;
\begin{array}{c}{\rm M}(-2,2)=\\ 
\left(
\begin{array}{ccccc}
 0 & 0 & 0 & 0 & 0 \\
  0 & 0 & 0 & 0 & 0 \\
 90 & 36 & 12 & 27 & 42 \\
 0 & 0 & 0 & 0 & 0 \\
180 & 72 & 24 & 66 & 108 \\
110 & 44 & 14 & 36 & 57 \\
136 & 58 & 22 & 24 & 24 \\
20 & 8 & 3 & 2 & 2 \\
48 & 21 & 8 & 10 & 10 \\
\end{array}
\right)\end{array}
$$

$$
\begin{array}{c}{\rm M}(-2,3)=\\\left(
\begin{array}{ccccc}
 0 & 0 & 0 & 0 & 0 \\
 0 & 0 & 0 & 0 & 0 \\
210 & 93 & 120 & 54 & 54 \\
 0 & 0 & 0 & 0 & 0 \\
 492 & 222 & 288 & 132 & 132 \\
 271 & 121 & 157 & 71 & 71 \\
 230 & 92 & 116 & 86 & 128 \\
 25 & 10 & 12 & 12 & 18 \\
 90 & 36 & 46 & 31 & 46 \\

\end{array}
\right) \end{array}\;
\begin{array}{c}{\rm M}^*(-2,\infty)=\\
\left(
\begin{array}{cc}
 0 & 0 \\
 0 & 0 \\
 0 & 1 \\
 0 & 0 \\
2 & 2 \\
 1 & 1 \\
 0 & 0 \\
 0 & 0 \\
 0 & 0 \\
\end{array}
\right)
 \end{array}\;
 \begin{array}{c}{\rm M}^\dagger(-2,\infty)=\\
 \left(
\begin{array}{cc}
 0 & 0 \\
 0 & 0 \\
 0 & 0 \\
 0 & 0 \\
 0 & 0 \\
 0 & 0 \\
10 & 4 \\
 1 & 1 \\
 4 & 1 \\
\end{array}
\right)\end{array}
$$
$$
\begin{array}{c}{\rm M}(-1,-2)=\\
\left(
\begin{array}{ccccc}
 0 & 0 & 0 & 0 & 0 \\
 0 & 0 & 0 & 0 & 0 \\
 48 & 21 & 8 & 10 & 10 \\
 0 & 0 & 0 & 0 & 0 \\
 136 & 58 & 22 & 24 & 24 \\
 48 & 21 & 9 & 9 & 9  \\
  90 &36 & 12 & 18 & 24 \\
   0 & 0 & 0 & 0 & 0 \\
  0& 0 & 0 & 3 & 6  \\
\end{array}
\right) \end{array}\;
\begin{array}{c}{\rm M}(-1,-1)=\\
\left(
\begin{array}{ccccc}
 0 & 0 & 0 & 0 & 0 \\
 0 & 0 & 0 & 0 & 0 \\
 10 & 4 & 6 & 3 & 6 \\
 0 & 0 & 0 & 0 & 0 \\
 30 & 12 & 16 & 10 & 16 \\
 10 & 4 & 5 & 4 & 7 \\
 18 & 9 & 12 & 6 & 6 \\
 0 & 0 & 0 & 0 & 0 \\
 0 & 0 & 0 & 0 & 0 \\
\end{array}
\right) \end{array}
$$

\clearpage
\section*{Appendix E: Incidence matrices for scenario III}

\begin{sidewaystable}
\renewcommand{\arraystretch}{1.5}
\renewcommand{\arraycolsep}{3 pt}
$$
\begin{array}{l}
{\tt M}_{--+}^{\rm III}=\;
\left(
\begin{array}{ccccc}
 -\frac{19}{3} & -\frac{25}{12} & -\frac{1}{6} & -\frac{1}{6} & -\frac{1}{6} \\
 -25 & -10 & -2 & -\frac{1}{2} & 1 \\
 -\frac{205}{6} & -\frac{41}{3} & -\frac{7}{3} & -\frac{13}{12} & \frac{7}{6} \\
 -\frac{58}{3} & -\frac{41}{6} & -\frac{5}{3} & \frac{1}{3} & \frac{1}{3} \\
 \frac{29}{3} & \frac{41}{12} & \frac{5}{6} & -\frac{1}{6} & -\frac{1}{6} \\
\end{array}
\right)
\;+\;
2^{|j|} \left(
\begin{array}{ccccc}
 \frac{7}{3} & \frac{14}{15} & \frac{4}{15} & \frac{4}{15} & \frac{4}{15} \\
 7 & \frac{14}{5} & \frac{4}{5} & \frac{4}{5} & \frac{4}{5} \\
 \frac{91}{6} & \frac{91}{15} & \frac{26}{15} & \frac{26}{15} & \frac{26}{15} \\
 \frac{28}{3} & \frac{56}{15} & \frac{16}{15} & \frac{16}{15} & \frac{16}{15} \\
 -\frac{7}{6} & -\frac{7}{15} & -\frac{2}{15} & -\frac{2}{15} & -\frac{2}{15} \\
\end{array}
\right)
\;+\;
3^{|j|} \left(
\begin{array}{ccccc}
 5 & \frac{43}{20} & \frac{9}{10} & \frac{9}{10} & \frac{9}{10} \\
 \frac{29}{3} & \frac{58}{15} & \frac{6}{5} & \frac{27}{10} & \frac{21}{5} \\
 \frac{29}{6} & \frac{29}{15} & \frac{3}{5} & \frac{27}{20} & \frac{21}{10} \\
 \frac{10}{3} & \frac{43}{30} & \frac{3}{5} & \frac{3}{5} & \frac{3}{5} \\
 \frac{5}{3} & \frac{43}{60} & \frac{3}{10} & \frac{3}{10} & \frac{3}{10} \\
\end{array}
\right)\\
\\
{\tt M}_{---}^{\rm III}\;=\;
\left(
\begin{array}{ccccc}
 \frac{5}{3} & \frac{2}{3} & \frac{5}{6} & \frac{1}{12} & \frac{5}{6} \\
 -31 & -\frac{23}{2} & -14 & -5 & -5 \\
 -\frac{229}{6} & -\frac{167}{12} & -\frac{52}{3} & -\frac{35}{6} & -\frac{35}{6} \\
 -\frac{40}{3} & -\frac{16}{3} & -\frac{17}{3} & -\frac{19}{6} & -\frac{5}{3} \\
 \frac{20}{3} & \frac{8}{3} & \frac{17}{6} & \frac{19}{12} & \frac{5}{6} \\
\end{array}
\right)
\; + \;
2^{|j|} \left(
\begin{array}{ccccc}
 \frac{19}{6} & \frac{19}{15} & \frac{19}{12} & \frac{19}{30} & \frac{19}{30} \\
 -7 & -\frac{14}{5} & -\frac{7}{2} & -\frac{7}{5} & -\frac{7}{5} \\
 \frac{29}{6} & \frac{29}{15} & \frac{29}{12} & \frac{29}{30} & \frac{29}{30} \\
 \frac{23}{3} & \frac{46}{15} & \frac{23}{6} & \frac{23}{15} & \frac{23}{15} \\
 -\frac{1}{3} & -\frac{2}{15} & -\frac{1}{6} & -\frac{1}{15} & -\frac{1}{15} \\
\end{array}
\right)
\; + \;
3^{|j|} \left(
\begin{array}{ccccc}
 \frac{7}{3} & \frac{14}{15} & \frac{7}{6} & \frac{73}{60} & \frac{59}{30} \\
 11 & \frac{47}{10} & 6 & \frac{13}{5} & \frac{13}{5} \\
 \frac{11}{2} & \frac{47}{20} & 3 & \frac{13}{10} & \frac{13}{10} \\
 \frac{14}{9} & \frac{28}{45} & \frac{7}{9} & \frac{73}{90} & \frac{59}{45} \\
 \frac{7}{9} & \frac{14}{45} & \frac{7}{18} & \frac{73}{180} & \frac{59}{90} \\
\end{array}
\right)
\\
\\
{\tt M}_{-++}^{\rm III}\;=\;
\left(
\begin{array}{ccccc}
 \frac{40}{3} & \frac{16}{3} & \frac{5}{3} & \frac{41}{12} & \frac{25}{6} \\
 -\frac{23}{3} & -\frac{13}{6} & -\frac{1}{3} & \frac{5}{3} & \frac{5}{3} \\
 -\frac{19}{6} & \frac{1}{12} & \frac{1}{6} & \frac{25}{6} & \frac{25}{6} \\
 10 & 4 & 2 & \frac{7}{2} & 5 \\
 -5 & -2 & -1 & -\frac{7}{4} & -\frac{5}{2} \\
\end{array}
\right)
\; + \;
2^{|j|-1} \left(
\begin{array}{ccccc}
 \frac{707}{24} & \frac{707}{60} & \frac{101}{30} & \frac{101}{30} & \frac{101}{30}
   \\
 -\frac{1673}{12} & -\frac{1673}{30} & -\frac{239}{15} & -\frac{239}{15} &
   -\frac{239}{15} \\
 -\frac{1183}{24} & -\frac{1183}{60} & -\frac{169}{30} & -\frac{169}{30} &
   -\frac{169}{30} \\
 \frac{133}{4} & \frac{133}{10} & \frac{19}{5} & \frac{19}{5} & \frac{19}{5} \\
 -14 & -\frac{28}{5} & -\frac{8}{5} & -\frac{8}{5} & -\frac{8}{5} \\
\end{array}
\right)
\; + \;
3^{|j|-1} \left(
\begin{array}{ccccc}
 33 & \frac{66}{5} & \frac{21}{5} & \frac{219}{20} & \frac{177}{10} \\
 153 & \frac{639}{10} & \frac{117}{5} & \frac{117}{5} & \frac{117}{5} \\
 \frac{153}{2} & \frac{639}{20} & \frac{117}{10} & \frac{117}{10} & \frac{117}{10}
   \\
 22 & \frac{44}{5} & \frac{14}{5} & \frac{73}{10} & \frac{59}{5} \\
 11 & \frac{22}{5} & \frac{7}{5} & \frac{73}{20} & \frac{59}{10} \\
\end{array}
\right)
\\
\\
{\tt M}_{-+-}^{\rm III}=
\left(
\begin{array}{ccccc}
 -\frac{43}{6} & -\frac{29}{12} & -\frac{7}{3} & -\frac{5}{6} & -\frac{5}{6} \\
 -\frac{20}{3} & -\frac{8}{3} & -\frac{7}{3} & -\frac{11}{6} & -\frac{1}{3} \\
 -\frac{35}{3} & -\frac{14}{3} & -\frac{23}{6} & -\frac{37}{12} & -\frac{5}{6} \\
 -11 & -\frac{7}{2} & -4 & -1 & -1 \\
 \frac{11}{2} & \frac{7}{4} & 2 & \frac{1}{2} & \frac{1}{2} \\
\end{array}
\right)
\; + \;
2^{|j|-1} \left(
\begin{array}{ccccc}
 -\frac{4}{3} & -\frac{8}{15} & -\frac{2}{3} & -\frac{4}{15} & -\frac{4}{15} \\
 \frac{128}{3} & \frac{256}{15} & \frac{64}{3} & \frac{128}{15} & \frac{128}{15} \\
 \frac{104}{3} & \frac{208}{15} & \frac{52}{3} & \frac{104}{15} & \frac{104}{15} \\
 8 & \frac{16}{5} & 4 & \frac{8}{5} & \frac{8}{5} \\
 -16 & -\frac{32}{5} & -8 & -\frac{16}{5} & -\frac{16}{5} \\
\end{array}
\right)
\; + \;
3^{|j|-1} \left(
\begin{array}{ccccc}
 \frac{63}{2} & \frac{279}{20} & 18 & \frac{81}{10} & \frac{81}{10} \\
 54 & \frac{108}{5} & 27 & \frac{243}{10} & \frac{189}{5} \\
 27 & \frac{54}{5} & \frac{27}{2} & \frac{243}{20} & \frac{189}{10} \\
 21 & \frac{93}{10} & 12 & \frac{27}{5} & \frac{27}{5} \\
 \frac{21}{2} & \frac{93}{20} & 6 & \frac{27}{10} & \frac{27}{10} \\
\end{array}
\right)
\end{array}
$$
\end{sidewaystable}

\begin{sidewaystable}
\renewcommand{\arraystretch}{1.5}
\renewcommand{\arraycolsep}{3 pt}
$$
\begin{array}{l}
{\tt M}_{+-+}^{\rm III}=
\left(
\begin{array}{ccccc}
 \frac{29}{3} & \frac{41}{12} & \frac{5}{6} & -\frac{1}{6} & -\frac{1}{6} \\
 -\frac{58}{3} & -\frac{41}{6} & -\frac{5}{3} & \frac{1}{3} & \frac{1}{3} \\
 -\frac{335}{6} & -\frac{67}{3} & -\frac{11}{3} & -\frac{23}{12} & \frac{11}{6} \\
 -25 & -10 & -2 & -\frac{1}{2} & 1 \\
 -\frac{19}{3} & -\frac{25}{12} & -\frac{1}{6} & -\frac{1}{6} & -\frac{1}{6} \\
\end{array}
\right)
\; + \;
2^{|j|} \left(
\begin{array}{ccccc}
 -\frac{7}{6} & -\frac{7}{15} & -\frac{2}{15} & -\frac{2}{15} & -\frac{2}{15} \\
 \frac{28}{3} & \frac{56}{15} & \frac{16}{15} & \frac{16}{15} & \frac{16}{15} \\
 \frac{49}{3} & \frac{98}{15} & \frac{28}{15} & \frac{28}{15} & \frac{28}{15} \\
 7 & \frac{14}{5} & \frac{4}{5} & \frac{4}{5} & \frac{4}{5} \\
 \frac{7}{3} & \frac{14}{15} & \frac{4}{15} & \frac{4}{15} & \frac{4}{15} \\
\end{array}
\right)\; + \;3^{|j|} \left(
\begin{array}{ccccc}
 \frac{5}{3} & \frac{43}{60} & \frac{3}{10} & \frac{3}{10} & \frac{3}{10} \\
 \frac{10}{3} & \frac{43}{30} & \frac{3}{5} & \frac{3}{5} & \frac{3}{5} \\
 \frac{29}{2} & \frac{29}{5} & \frac{9}{5} & \frac{81}{20} & \frac{63}{10} \\
 \frac{29}{3} & \frac{58}{15} & \frac{6}{5} & \frac{27}{10} & \frac{21}{5} \\
 5 & \frac{43}{20} & \frac{9}{10} & \frac{9}{10} & \frac{9}{10} \\
\end{array}
\right)\\
\\
{\tt M}_{+--}^{\rm III}=
\left(
\begin{array}{ccccc}
 \frac{20}{3} & \frac{8}{3} & \frac{17}{6} & \frac{19}{12} & \frac{5}{6} \\
 -\frac{40}{3} & -\frac{16}{3} & -\frac{17}{3} & -\frac{19}{6} & -\frac{5}{3} \\
 -\frac{365}{6} & -\frac{265}{12} & -\frac{83}{3} & -\frac{55}{6} & -\frac{55}{6} \\
 -31 & -\frac{23}{2} & -14 & -5 & -5 \\
 \frac{5}{3} & \frac{2}{3} & \frac{5}{6} & \frac{1}{12} & \frac{5}{6} \\
\end{array}
\right)
\; + \;
2^{|j|} \left(
\begin{array}{ccccc}
 -\frac{1}{3} & -\frac{2}{15} & -\frac{1}{6} & -\frac{1}{15} & -\frac{1}{15} \\
 \frac{23}{3} & \frac{46}{15} & \frac{23}{6} & \frac{23}{15} & \frac{23}{15} \\
 -\frac{19}{3} & -\frac{38}{15} & -\frac{19}{6} & -\frac{19}{15} & -\frac{19}{15} \\
 -7 & -\frac{14}{5} & -\frac{7}{2} & -\frac{7}{5} & -\frac{7}{5} \\
 \frac{19}{6} & \frac{19}{15} & \frac{19}{12} & \frac{19}{30} & \frac{19}{30} \\
\end{array}
\right)\; + \;3^{|j|} \left(
\begin{array}{ccccc}
 \frac{7}{9} & \frac{14}{45} & \frac{7}{18} & \frac{73}{180} & \frac{59}{90} \\
 \frac{14}{9} & \frac{28}{45} & \frac{7}{9} & \frac{73}{90} & \frac{59}{45} \\
 \frac{33}{2} & \frac{141}{20} & 9 & \frac{39}{10} & \frac{39}{10} \\
 11 & \frac{47}{10} & 6 & \frac{13}{5} & \frac{13}{5} \\
 \frac{7}{3} & \frac{14}{15} & \frac{7}{6} & \frac{73}{60} & \frac{59}{30} \\
\end{array}
\right)
\\
\\
{\tt M}_{+++}^{\rm III}=
\left(
\begin{array}{ccccc}
 -5 & -2 & -1 & -\frac{7}{4} & -\frac{5}{2} \\
 10 & 4 & 2 & \frac{7}{2} & 5 \\
 -\frac{5}{2} & \frac{5}{4} & \frac{1}{2} & \frac{15}{2} & \frac{15}{2} \\
 -\frac{23}{3} & -\frac{13}{6} & -\frac{1}{3} & \frac{5}{3} & \frac{5}{3} \\
 \frac{40}{3} & \frac{16}{3} & \frac{5}{3} & \frac{41}{12} & \frac{25}{6} \\
\end{array}
\right)\; + \;
2^{|j|-1} \left(
\begin{array}{ccccc}
 -14 & -\frac{28}{5} & -\frac{8}{5} & -\frac{8}{5} & -\frac{8}{5} \\
 \frac{133}{4} & \frac{133}{10} & \frac{19}{5} & \frac{19}{5} & \frac{19}{5} \\
 -\frac{357}{2} & -\frac{357}{5} & -\frac{102}{5} & -\frac{102}{5} & -\frac{102}{5}
   \\
 -\frac{1673}{12} & -\frac{1673}{30} & -\frac{239}{15} & -\frac{239}{15} &
   -\frac{239}{15} \\
 \frac{707}{24} & \frac{707}{60} & \frac{101}{30} & \frac{101}{30} & \frac{101}{30}
   \\
\end{array}
\right)\; + \;3^{|j|-1} \left(
\begin{array}{ccccc}
 11 & \frac{22}{5} & \frac{7}{5} & \frac{73}{20} & \frac{59}{10} \\
 22 & \frac{44}{5} & \frac{14}{5} & \frac{73}{10} & \frac{59}{5} \\
 \frac{459}{2} & \frac{1917}{20} & \frac{351}{10} & \frac{351}{10} & \frac{351}{10}
   \\
 153 & \frac{639}{10} & \frac{117}{5} & \frac{117}{5} & \frac{117}{5} \\
 33 & \frac{66}{5} & \frac{21}{5} & \frac{219}{20} & \frac{177}{10} \\
\end{array}
\right)
\\
\\
{\tt M}_{++-}^{\rm III}=
\left(
\begin{array}{ccccc}
 \frac{11}{2} & \frac{7}{4} & 2 & \frac{1}{2} & \frac{1}{2} \\
 -11 & -\frac{7}{2} & -4 & -1 & -1 \\
 -20 & -8 & -\frac{13}{2} & -\frac{21}{4} & -\frac{3}{2} \\
 -\frac{20}{3} & -\frac{8}{3} & -\frac{7}{3} & -\frac{11}{6} & -\frac{1}{3} \\
 -\frac{43}{6} & -\frac{29}{12} & -\frac{7}{3} & -\frac{5}{6} & -\frac{5}{6} \\
\end{array}
\right)
\; + \;
2^{|j|-1} \left(
\begin{array}{ccccc}
 -16 & -\frac{32}{5} & -8 & -\frac{16}{5} & -\frac{16}{5} \\
 8 & \frac{16}{5} & 4 & \frac{8}{5} & \frac{8}{5} \\
 84 & \frac{168}{5} & 42 & \frac{84}{5} & \frac{84}{5} \\
 \frac{128}{3} & \frac{256}{15} & \frac{64}{3} & \frac{128}{15} & \frac{128}{15} \\
 -\frac{4}{3} & -\frac{8}{15} & -\frac{2}{3} & -\frac{4}{15} & -\frac{4}{15} \\
\end{array}
\right)\; + \;3^{|j|-1} \left(
\begin{array}{ccccc}
 \frac{21}{2} & \frac{93}{20} & 6 & \frac{27}{10} & \frac{27}{10} \\
 21 & \frac{93}{10} & 12 & \frac{27}{5} & \frac{27}{5} \\
 81 & \frac{162}{5} & \frac{81}{2} & \frac{729}{20} & \frac{567}{10} \\
 54 & \frac{108}{5} & 27 & \frac{243}{10} & \frac{189}{5} \\
 \frac{63}{2} & \frac{279}{20} & 18 & \frac{81}{10} & \frac{81}{10} \\
\end{array}
\right)
\end{array}
$$
\end{sidewaystable}

\clearpage
\section*{Appendix F: Incidence matrices for scenario IV}

\begin{sidewaystable}
\renewcommand{\arraystretch}{1.5}
\renewcommand{\arraycolsep}{3 pt}
$$
\begin{array}{l}
{\tt M}_{-+}^{\rm IV}=
\left(
\begin{array}{ccccc}
 0 & 0 & 0 & 0 & 0 \\
 0 & 0 & 0 & 0 & 0 \\
 -\frac{58}{3} & -\frac{41}{6} & -\frac{5}{3} & \frac{1}{3} &
   \frac{1}{3} \\
 0 & 0 & 0 & 0 & 0 \\
 -\frac{58}{3} & -\frac{41}{6} & -\frac{5}{3} & \frac{1}{3} &
   \frac{1}{3} \\
 -\frac{19}{3} & -\frac{25}{12} & -\frac{1}{6} & -\frac{1}{6} &
   -\frac{1}{6} \\
 -25 & -10 & -2 & -\frac{1}{2} & 1 \\
 \frac{25}{2} & 5 & 1 & \frac{1}{4} & -\frac{1}{2} \\
 -25 & -10 & -2 & -\frac{1}{2} & 1 \\
\end{array}
\right)
\; + \;
2^{|j|} \left(
\begin{array}{ccccc}
 0 & 0 & 0 & 0 & 0 \\
 0 & 0 & 0 & 0 & 0 \\
 \frac{28}{3} & \frac{56}{15} & \frac{16}{15} & \frac{16}{15} &
   \frac{16}{15} \\
 0 & 0 & 0 & 0 & 0 \\
 \frac{49}{3} & \frac{98}{15} & \frac{28}{15} & \frac{28}{15} &
   \frac{28}{15} \\
 \frac{7}{3} & \frac{14}{15} & \frac{4}{15} & \frac{4}{15} &
   \frac{4}{15} \\
 \frac{91}{4} & \frac{91}{10} & \frac{13}{5} & \frac{13}{5} &
   \frac{13}{5} \\
 -\frac{7}{2} & -\frac{7}{5} & -\frac{2}{5} & -\frac{2}{5} &
   -\frac{2}{5} \\
 \frac{7}{4} & \frac{7}{10} & \frac{1}{5} & \frac{1}{5} & \frac{1}{5}
   \\
\end{array}
\right)
\; + \;
3^{|j|} \left(
\begin{array}{ccccc}
 0 & 0 & 0 & 0 & 0 \\
 0 & 0 & 0 & 0 & 0 \\
 \frac{10}{3} & \frac{43}{30} & \frac{3}{5} & \frac{3}{5} &
   \frac{3}{5} \\
 0 & 0 & 0 & 0 & 0 \\
 10 & \frac{43}{10} & \frac{9}{5} & \frac{9}{5} & \frac{9}{5} \\
 5 & \frac{43}{20} & \frac{9}{10} & \frac{9}{10} & \frac{9}{10} \\
 \frac{29}{9} & \frac{58}{45} & \frac{2}{5} & \frac{9}{10} &
   \frac{7}{5} \\
 \frac{29}{54} & \frac{29}{135} & \frac{1}{15} & \frac{3}{20} &
   \frac{7}{30} \\
 \frac{29}{27} & \frac{58}{135} & \frac{2}{15} & \frac{3}{10} &
   \frac{7}{15} \\
\end{array}
\right)
\\
\\
{\tt M}_{--}^{\rm IV}=
\left(
\begin{array}{ccccc}
 0 & 0 & 0 & 0 & 0 \\
 0 & 0 & 0 & 0 & 0 \\
 -\frac{40}{3} & -\frac{16}{3} & -\frac{17}{3} & -\frac{19}{6} &
   -\frac{5}{3} \\
 0 & 0 & 0 & 0 & 0 \\
 -\frac{40}{3} & -\frac{16}{3} & -\frac{17}{3} & -\frac{19}{6} &
   -\frac{5}{3} \\
 \frac{5}{3} & \frac{2}{3} & \frac{5}{6} & \frac{1}{12} & \frac{5}{6}
   \\
 -31 & -\frac{23}{2} & -14 & -5 & -5 \\
 \frac{31}{2} & \frac{23}{4} & 7 & \frac{5}{2} & \frac{5}{2} \\
 -31 & -\frac{23}{2} & -14 & -5 & -5 \\
\end{array}
\right)
\; + \;
2^{|j|} \left(
\begin{array}{ccccc}
 0 & 0 & 0 & 0 & 0 \\
 0 & 0 & 0 & 0 & 0 \\
 \frac{23}{3} & \frac{46}{15} & \frac{23}{6} & \frac{23}{15} &
   \frac{23}{15} \\
 0 & 0 & 0 & 0 & 0 \\
 \frac{44}{3} & \frac{88}{15} & \frac{22}{3} & \frac{44}{15} &
   \frac{44}{15} \\
 \frac{19}{6} & \frac{19}{15} & \frac{19}{12} & \frac{19}{30} &
   \frac{19}{30} \\
 11 & \frac{22}{5} & \frac{11}{2} & \frac{11}{5} & \frac{11}{5} \\
 -4 & -\frac{8}{5} & -2 & -\frac{4}{5} & -\frac{4}{5} \\
 2 & \frac{4}{5} & 1 & \frac{2}{5} & \frac{2}{5} \\
\end{array}
\right)
\; + \;
3^{|j|} \left(
\begin{array}{ccccc}
 0 & 0 & 0 & 0 & 0 \\
 0 & 0 & 0 & 0 & 0 \\
 \frac{14}{9} & \frac{28}{45} & \frac{7}{9} & \frac{73}{90} &
   \frac{59}{45} \\
 0 & 0 & 0 & 0 & 0 \\
 \frac{14}{3} & \frac{28}{15} & \frac{7}{3} & \frac{73}{30} &
   \frac{59}{15} \\
 \frac{7}{3} & \frac{14}{15} & \frac{7}{6} & \frac{73}{60} &
   \frac{59}{30} \\
 \frac{11}{3} & \frac{47}{30} & 2 & \frac{13}{15} & \frac{13}{15} \\
 \frac{11}{18} & \frac{47}{180} & \frac{1}{3} & \frac{13}{90} &
   \frac{13}{90} \\
 \frac{11}{9} & \frac{47}{90} & \frac{2}{3} & \frac{13}{45} &
   \frac{13}{45} \\
\end{array}
\right)
\end{array}
$$
\end{sidewaystable}

\begin{sidewaystable}
\renewcommand{\arraystretch}{1.5}
\renewcommand{\arraycolsep}{3 pt}
$$
\begin{array}{l}
{\tt M}_{++}^{\rm IV}=
\left(
\begin{array}{ccccc}
 0 & 0 & 0 & 0 & 0 \\
 0 & 0 & 0 & 0 & 0 \\
 10 & 4 & 2 & \frac{7}{2} & 5 \\
 0 & 0 & 0 & 0 & 0 \\
 10 & 4 & 2 & \frac{7}{2} & 5 \\
 \frac{40}{3} & \frac{16}{3} & \frac{5}{3} & \frac{41}{12} &
   \frac{25}{6} \\
 -\frac{23}{3} & -\frac{13}{6} & -\frac{1}{3} & \frac{5}{3} &
   \frac{5}{3} \\
 \frac{23}{6} & \frac{13}{12} & \frac{1}{6} & -\frac{5}{6} &
   -\frac{5}{6} \\
 -\frac{23}{3} & -\frac{13}{6} & -\frac{1}{3} & \frac{5}{3} &
   \frac{5}{3} \\
\end{array}
\right)
\; + \;
2^{|j|-1} \left(
\begin{array}{ccccc}
 0 & 0 & 0 & 0 & 0 \\
 0 & 0 & 0 & 0 & 0 \\
 \frac{133}{4} & \frac{133}{10} & \frac{19}{5} & \frac{19}{5} &
   \frac{19}{5} \\
 0 & 0 & 0 & 0 & 0 \\
 \frac{77}{2} & \frac{77}{5} & \frac{22}{5} & \frac{22}{5} &
   \frac{22}{5} \\
 \frac{707}{24} & \frac{707}{60} & \frac{101}{30} & \frac{101}{30} &
   \frac{101}{30} \\
 -\frac{119}{3} & -\frac{238}{15} & -\frac{68}{15} & -\frac{68}{15} &
   -\frac{68}{15} \\
 -\frac{161}{12} & -\frac{161}{30} & -\frac{23}{15} & -\frac{23}{15} &
   -\frac{23}{15} \\
 -\frac{77}{12} & -\frac{77}{30} & -\frac{11}{15} & -\frac{11}{15} &
   -\frac{11}{15} \\
\end{array}
\right)
\; + \;
3^{|j|-1} \left(
\begin{array}{ccccc}
 0 & 0 & 0 & 0 & 0 \\
 0 & 0 & 0 & 0 & 0 \\
 22 & \frac{44}{5} & \frac{14}{5} & \frac{73}{10} & \frac{59}{5} \\
 0 & 0 & 0 & 0 & 0 \\
 66 & \frac{132}{5} & \frac{42}{5} & \frac{219}{10} & \frac{177}{5} \\
 33 & \frac{66}{5} & \frac{21}{5} & \frac{219}{20} & \frac{177}{10} \\
 51 & \frac{213}{10} & \frac{39}{5} & \frac{39}{5} & \frac{39}{5} \\
 \frac{17}{2} & \frac{71}{20} & \frac{13}{10} & \frac{13}{10} &
   \frac{13}{10} \\
 17 & \frac{71}{10} & \frac{13}{5} & \frac{13}{5} & \frac{13}{5} \\
\end{array}
\right)
\\
\\
{\tt M}_{+-}^{\rm IV}=
\left(
\begin{array}{ccccc}
 0 & 0 & 0 & 0 & 0 \\
 0 & 0 & 0 & 0 & 0 \\
 -11 & -\frac{7}{2} & -4 & -1 & -1 \\
 0 & 0 & 0 & 0 & 0 \\
 -11 & -\frac{7}{2} & -4 & -1 & -1 \\
 -\frac{43}{6} & -\frac{29}{12} & -\frac{7}{3} & -\frac{5}{6} &
   -\frac{5}{6} \\
 -\frac{20}{3} & -\frac{8}{3} & -\frac{7}{3} & -\frac{11}{6} &
   -\frac{1}{3} \\
 \frac{10}{3} & \frac{4}{3} & \frac{7}{6} & \frac{11}{12} &
   \frac{1}{6} \\
 -\frac{20}{3} & -\frac{8}{3} & -\frac{7}{3} & -\frac{11}{6} &
   -\frac{1}{3} \\
\end{array}
\right)
\; + \;
2^{|j|-1} \left(
\begin{array}{ccccc}
 0 & 0 & 0 & 0 & 0 \\
 0 & 0 & 0 & 0 & 0 \\
 8 & \frac{16}{5} & 4 & \frac{8}{5} & \frac{8}{5} \\
 0 & 0 & 0 & 0 & 0 \\
 -16 & -\frac{32}{5} & -8 & -\frac{16}{5} & -\frac{16}{5} \\
 -\frac{4}{3} & -\frac{8}{15} & -\frac{2}{3} & -\frac{4}{15} &
   -\frac{4}{15} \\
 \frac{56}{3} & \frac{112}{15} & \frac{28}{3} & \frac{56}{15} &
   \frac{56}{15} \\
 -\frac{4}{3} & -\frac{8}{15} & -\frac{2}{3} & -\frac{4}{15} &
   -\frac{4}{15} \\
 \frac{32}{3} & \frac{64}{15} & \frac{16}{3} & \frac{32}{15} &
   \frac{32}{15} \\
\end{array}
\right)
\; + \;
3^{|j|-1} \left(
\begin{array}{ccccc}
 0 & 0 & 0 & 0 & 0 \\
 0 & 0 & 0 & 0 & 0 \\
 21 & \frac{93}{10} & 12 & \frac{27}{5} & \frac{27}{5} \\
 0 & 0 & 0 & 0 & 0 \\
 63 & \frac{279}{10} & 36 & \frac{81}{5} & \frac{81}{5} \\
 \frac{63}{2} & \frac{279}{20} & 18 & \frac{81}{10} & \frac{81}{10} \\
 18 & \frac{36}{5} & 9 & \frac{81}{10} & \frac{63}{5} \\
 3 & \frac{6}{5} & \frac{3}{2} & \frac{27}{20} & \frac{21}{10} \\
 6 & \frac{12}{5} & 3 & \frac{27}{10} & \frac{21}{5} \\
\end{array}
\right)
\end{array}
$$
\end{sidewaystable}

\clearpage
%---------------------------------------------------------------------------------------------

\end{document}